\documentclass[10pt]{article}

\usepackage[T1]{fontenc}
\usepackage{amsmath,amsfonts,amssymb,amsthm,bbm}
\usepackage{graphicx,color}
\usepackage{epsfig,psfrag, fancyhdr}

\newtheorem{dfn}{Definition}[section]
\newtheorem{prop}[dfn]{Proposition}
\newtheorem{thm}[dfn]{Theorem}
\newtheorem{lem}[dfn]{Lemma}
\newtheorem{cor}[dfn]{Corollary}
\theoremstyle{remark}
\newtheorem{rem}[dfn]{Remark}

\DeclareMathOperator{\cyl}{cyl}
\DeclareMathOperator{\hyp}{hyp}
\DeclareMathOperator{\card}{card}
\DeclareMathOperator{\diam}{diam}
\DeclareMathOperator{\tr}{tr}

\newenvironment{dem}{\vskip 2mm\noindent {\it Proof} :}
                    {\hfill $\square$ \vskip 2mm \noindent}
 
\newcommand{\eps}{\varepsilon}

\def\PP{\mathbb{P}}
\def\RR{\mathbb{R}}
\def\EE{\mathbb{E}}
\def\NN{\mathbb{N}}
\def\ZZ{\mathbb{Z}}
\def\E{\mathcal{E}}
\def\H{\mathcal{H}}

\def\II{\mbox{ 1\hskip -.29em I}} 

\begin{document}

\thispagestyle{empty}

\title{\huge Lower large deviations and laws of large numbers for maximal flows through a box in first
  passage percolation}

\date{}
\maketitle

\begin{center}
\vskip-1cm {\Large Rapha\"el Rossignol\footnote{Rapha\"el Rossignol was supported by the Swiss National
    Science Foundation grants 200021-1036251/1 and 200020-112316/1.}}\\
{\it Universit\'e Paris Sud, Laboratoire de Math\'ematiques, b\^atiment 425,
91405 Orsay Cedex, France}\\
{\it E-mail:} raphael.rossignol@math.u-psud.fr\\
and\\
\vskip1cm {\Large Marie Th\'eret}\\
{\it \'Ecole Normale Sup\'erieure, D\'epartement Math\'ematiques et
Applications, 45 rue d'Ulm, 75230 Paris Cedex 05, France}\\
{\it E-mail:} marie.theret@ens.fr
\end{center}








\noindent
{\bf Abstract.} We consider the standard first passage percolation model in $\mathbb{Z}^d$ for
$d\geq 2$. We are interested in two quantities, the maximal flow $\tau$
between the lower half and the upper half of the box, and the maximal flow
$\phi$ between the top and the bottom of the box. A standard subadditive
argument yields the law of large numbers for $\tau$ in rational
directions. Kesten and Zhang have proved the law of large numbers for $\tau$ and $\phi$
when the sides of the box are parallel to the coordinate hyperplanes: the two variables grow linearly with the surface $s$ of the basis of the box, with the
same deterministic speed. We study
the probabilities that the rescaled variables $\tau /s$ and $\phi /s$ are
abnormally small. For $\tau$, the box can have any orientation,
whereas for $\phi$, we require either that the
box is sufficiently flat, or that its sides are parallel to the coordinate
hyperplanes. We show that these
probabilities decay exponentially fast with $s$, when $s$ grows to
infinity. Moreover, we prove an associated large deviation principle of speed
$s$ for $\tau /s$ and $\phi /s$, and we improve the conditions required to obtain the law of large numbers for these variables.\\





\section{Introduction}
The model of maximal flow in a randomly porous medium with independent and
identically distributed capacities has been initially studied by Kesten
(see \cite{Kesten:flows}), who introduced it as a ``higher dimensional version of
First Passage Percolation''. The purpose of this model is to understand the behaviour of the maximum
amount of flow that can cross the medium from one part to another. 

All the precise definitions will be given in section \ref{sec:def}, but let
us be a little more accurate. The random medium is represented by the lattice
$\ZZ^d$. We see each edge as a microscopic pipe which the fluid
can flow through. To each edge $e$, we attach a nonnegative capacity $t(e)$ which
represents the amount of fluid (or the amount of fluid per unit of time) that
can effectively go through the edge $e$. Capacities are then supposed to be
random, identically and independently distributed with common distribution
function $F$. Let $A$ be some hyperrectangle in $\RR^d$ (i.e. a box of
dimension $d-1$) and $n$ be an integer. The portion of media
that we will look at is a box $B_n$ of basis $nA$ and of height $2h(n)$, which
 $nA$ splits into two equal parts. The boundary of $B_n$ is thus split
 into two parts, $A_n^1$ and $A_n^2$. We define two flows through
 $B_n$: the maximal flow
$\tau_n$ for which the fluid can enter the box through $A_n^1$ and leave it
through $A_n^2$, and the maximal flow $\phi_n$ for which the fluid enters $B_n$ only through its
bottom side and leaves it through its top side. Existing results for $\phi_n$
and $\tau_n$ are
essentially of two
types: laws of large numbers and large deviation results. Subadditivity
implies a law of large numbers for $\tau_n$, when $B_n$ is oriented to
a rational direction (as defined in \cite{Boivin}). It is important to note that all the
results  concerning $\phi_n$ we present now were obtained for ``straight'' hyperrectangles $A$, i.e. hyperrectangles of
the form $\prod_{i=1}^{d-1} [0,a_i] \times \{0\}$. Due to the symmetries
of the lattice $\ZZ^d$, this simplifies considerably the task.
Kesten proved a law of large
numbers for $\phi_n$ in straight cylinders in $\ZZ^3$ (see
\cite{Kesten:flows}), under various conditions on the height $h(n)$, the value of $F(0)$ and
an exponential moment condition on $F$. In a remarkable paper, Zhang (see \cite{Zhang07})
recently optimized  Kesten's
condition on $F(0)$ and extended the result to $\ZZ^d$, $d\geq 2$ (see Theorem \ref{cvphi} below). Th\'eret proved a large deviation principle for $\phi_n$ at volume order for upper deviations (see
\cite{TheretUpper}). Lower large deviations for $\phi_n$ far from its
  asymptotic behaviour were investigated for Bernoulli
capacities in \cite{Chayes}, and for general functions in \cite{Theret:small},
and are shown to be of surface order, although a full large deviation
principle was not proved. 

The main results of this paper are the lower large deviation principles for $\tau_n$
and $\phi_n$ under various conditions, and the improvement of the
 moment conditions required to state the law of large numbers for these variables. More precisely, we shall show lower large deviation principles at the surface order for
$\tau_n$ for general $A$ and height $h(n)$, and for $\phi_n$ when $h(n)$ is
small compared to $n$ (see Theorem \ref{thmpgd} and Corollary \ref{corthmpgd}). We also show a
lower large deviation principle at the surface order for $\phi_n$ when $\log
h(n)$ is small compared to $n^{d-1}$ and when $A$ is straight (see Theorem
\ref{thmpgdphi}). Unfortunately, when $d\geq 3$, we are not able to prove the lower large
deviation principle for $\phi_n$ through general hyperrectangles and heights (see Remark \ref{rem:phiechoue}). Incidentally, we prove deviation results which are interesting on
their own for $\phi_n$ and $\tau_n$, for
general hyperrectangles $A$ (see Theorem \ref{devinf}, Theorem
\ref{cordevinf} and Theorem \ref{devinfphi}). A consequence of
these deviation results is a law of large numbers for $\tau_n$ in any
fixed direction, even irrational, under an optimal moment
condition (see Theorem \ref{thm:llntau}).  We also obtain a law of large numbers
for $\phi_n$ in straight boxes under an optimal condition on the
height of the box, and a weak moment condition. We stress the fact
that we do not use any subadditive ergodic theorem for the law of 
large numbers for $\tau_n$, since in our general setting, subadditivity of $\tau_n$
is lost in irrational directions. Instead, we use the ``almost
subadditivity'' of $\tau_n$ combined with a lower deviation inequality.

The paper is organized as follows. In section \ref{sec:def}, we give the precise
definitions and notations. In section \ref{sec:results}, we state the important
background we shall rely on and the main results of the paper. In section
\ref{sec:lowerdev}, we prove  the deviation results for $\tau_n$, and
for $\phi_n$ in flat cylinders, the proof of the corresponding result for
  $\phi_n$ in straight boxes being completed at the end of the
  paper. We also obtain also the law of large numbers for $\tau_n$ in this
section. Section \ref{sec:pgdtau} is devoted to the large deviation principle
for $\tau_n$, and its corollary, the large deviation principle for $\phi_n$ in
flat boxes. Finally, we prove the law of large numbers, the order of the lower
  large deviations and the large deviation
  principle for $\phi_n$ in straight boxes in section \ref{sec:pgdphi}.

\section{Definitions and notations}
\label{sec:def}

The most important notations are gathered in
this section.
\subsection{Maximal flow on a graph}
\label{subsec:maxflow}

First, let us define the notion of a flow on a finite
unoriented graph $G=(V,\E)$ with set of vertices $V$ and set of edges
$\E$.  We write $x\sim y$ when $x$ and $y$ are two neighbouring
vertices in $G$. Let
$t=(t(e))_{e\in \E}$ be a collection of non-negative real numbers,
which are called \emph{capacities}. It means that $t(e)$ is the maximal amount of
    fluid that can go through the edge $e$ per unit of time. To each edge $e$, one may associate two oriented edges, and we
shall denote by $\overrightarrow{\E}$ the set of all these oriented edges. Let $Y$ and $Z$ be
two finite, disjoint, non-empty sets of vertices of $G$: $Y$ denotes the source of the
network, and $Z$ the sink. A function $\theta$
on $\overrightarrow{\E}$ is called a flow from $Y$ to $Z$ with strength
$\|\theta\|$ and capacities $t$ if it is
antisymmetric, i.e.
$\theta_{\overrightarrow{xy}}=-\theta_{\overrightarrow{yx}}$, if it satisfies the
node law at each vertex $x$ of $V\smallsetminus (Y\cup Z)$:
$$\sum_{y\sim x}\theta_{\overrightarrow{xy}}=0\;,$$
if it satisfies the capacity constraints:
$$\forall e\in \E,\;|\theta(e)|\leq t(e)\;,$$
and if the ``flow in'' at $Y$ and the ``flow out'' at $Z$ equal $\|\theta\|$:
$$\|\theta\|=\sum_{y\in Y}\sum_{\substack{x\sim y\\ x\not \in
    Y}}\theta(\overrightarrow{yx})=\sum_{z\in Z}\sum_{\substack{x\sim
    z\\ x\not \in Z}}\theta(\overrightarrow{xz})\;.$$
The \emph{maximal flow from $Y$ to $Z$}, denoted by $\phi_t(G,Y,Z)$, is defined as
the maximum strength of all flows from $Y$ to $Z$ with capacities
$t$. \emph{ We stress the fact that $\phi_t(G,Y,Z)$ is non-negative for any
  choice of $G$, $Y$ and $Z$.} We shall in general omit the subscript $t$ when it is understood
from the context. The \emph{max-flow min-cut theorem} (see \cite{Bollobas98} for instance)
asserts that the maximal flow from $Y$ to $Z$ equals the minimal
capacity of a cut between $Y$ and $Z$. Precisely, let us say that
$E\subset\E$ is a cut between $Y$ and $Z$ in $G$ if every path from
$Y$ to $Z$ borrows at least one edge of $E$. Define $V(E)=\sum_{e\in
  E}t(e)$ to be the capacity of a cut $E$. Then,
\begin{equation}
\label{eq:maxflowmincut}
\phi_t(G,Y,Z)=\min\{V(E)\mbox{ s.t. }E\mbox{ is a cut between
}Y\mbox{ and }Z \mbox{ in } G\}\;.
\end{equation}

\subsection{On the cubic lattice}
\label{subsec:cubiclattice}
We use many notations introduced in \cite{Kesten:StFlour} and
\cite{Kesten:flows}. Let $d\geq2$. We consider the graph $(\mathbb{Z}^{d},
\mathbb E ^{d})$ having for vertices $\mathbb Z ^{d}$ and for edges
$\mathbb E ^{d}$, the set of pairs of nearest neighbours for the standard
$L^{1}$ norm: $\| z\|_1 = \sum_{i=1}^d |z_i|$ for $z=(z_1,...,z_d)\in \RR^d$. To each edge $e$ in $\mathbb{E}^{d}$ we assign a \emph{random
capacity} $t(e)$ with values in $\mathbb{R}^{+}$. We suppose that \emph{the family
$(t(e), e \in \mathbb{E}^{d})$ is independent and identically distributed,
with a common distribution function $F$:} this is the standard model of
first passage percolation on the graph $(\mathbb{Z}^d, \mathbb{E}^d)$. More
formally, we take the product
measure $\mathbb {P}$ on $\Omega= \prod_{e\in \mathbb{E}^{d}} [0, \infty[$,
    and we write its expectation $\mathbb{E}$. 

For a subset $X$ of $\mathbb{R}^d$, we denote by \emph{$\mathcal{H}^s (X)$ the
$s$-dimensional Hausdorff measure of $X$} (we will use $s=d-1$ and $s=d-2$).
Let $A$ {\bf$ \subset \RR^d$} be a non-degenerate \emph{hyperrectangle} (for the
usual scalar product), i.e., a box of dimension $d-1$ in
$\mathbb{R}^d$. \emph{All hyperrectangles will be supposed to be closed and
  non-degenerate in $\mathbb{R}^d$}. Thus, every hyperrectangle $A$ we will
consider is the image by an isometry of $\RR^d$ of a set of the form
$\prod_{i=1}^{d-1}[0, k_i]\times \{0\}$ for strictly positive real numbers
$k_i$. With this notation, we define \emph{the smallest length of
  $A$}, denoted by $l_{min}(A)$ as:
$$l_{min}(A)=\min_{i=1\ldots d-1}k_i\;,$$
i.e. the smallest length of a side of $A$. We denote by
$\vec{v}$ one of 
the two vectors of unit euclidean norm, orthogonal to $\hyp (A)$, the hyperplane spanned by
$A$. For $h$ a positive real number, we denote by \emph{$\cyl(A,h)$ the cylinder of basis $A$ and height $2h$},
i.e., the set 
$$ \cyl (A,h) \,=\, \{x+t \vec{v} \,|\, x\in A \,,\,  t\in
[-h,h]    \}\,.$$
The set $\cyl(A,h) \smallsetminus \hyp (A)$ has two connected
components, which we denote by $\mathcal{C}_1(A,h)$ and
$\mathcal{C}_2(A,h)$. For $i=1,2$, let $A_i^h$ be
the set of the points in $\mathcal{C}_i(A,h) \cap \mathbb{Z}^d$ which have
a nearest neighbour in $\mathbb{Z}^d \smallsetminus \cyl(A,h)$:
$$A_i^h\,=\,\{x\in \mathcal{C}_i(A,h) \cap
\mathbb{Z}^d \,|\, \exists y \in \mathbb{Z}^d \smallsetminus \cyl(A,h) \,,\,
\|x-y\|_{1} =1 \}\,.$$
Let $T(A,h)$ (respectively $B(A,h)$) be the top
(respectively the bottom) of $\cyl(A,h)$, i.e.,
$$ T(A,h) \,=\, \{ x\in \cyl(A,h) \,|\, \exists y\notin \cyl(A,h)\,,\,\,
\langle x,y\rangle \in \mathbb{E}^d \,\, and\,\,\langle x,y\rangle
\,\,intersects\,\, A+h\vec{v}  \}  $$
and
$$  B(A,h) \,=\, \{ x\in \cyl(A,h) \,|\, \exists y\notin \cyl(A,h)\,,\,\,
\langle x,y\rangle \in \mathbb{E}^d \,\, and\,\,\langle x,y\rangle
\,\,intersects\,\, A-h\vec{v}  \} \,.$$
The notation $\langle x,y\rangle$ corresponds to the edge of endpoints
$x$ and $y$. We define also the \emph{$r$-neighbourhood $\mathcal{V} (H,r)$ of a subset $H$ of $\mathbb{R}^d$} as
$$\mathcal{V}(H,r) \,=\, \{ x \in \mathbb{R}^d \,|\, d(x,H)<r\}\,,$$
where the distance is the euclidean one, i.e.  $d(x,H) =
\inf \{\|x-y\|_2 \,|\, y\in H \}$ and  $\| z\|_2 = \sqrt{\sum_{i=1}^d
  z_i^2}$ for $z=(z_1,...,z_d)\in \RR^d$.\\

\noindent\emph{The main characters.} For a given realization $(t(e),e\in \mathbb{E}^{d})$ we define $\tau (A,h)$ by:
$$ \tau(A,h) \,=\,  \phi_t(\cyl(A,h)\cap \ZZ^d,A_1^h,A_2^h ) \;,$$
where $\phi_t$ is defined in section \ref{subsec:maxflow} and
$\cyl(A,h)\cap \ZZ^d$ denotes the induced subgraph of $\ZZ^d$ with set
of vertices $\cyl(A,h)\cap \ZZ^d$, equipped with capacities
$t$. This definition makes sense if $A_1^h$ and $A_2^h$ are
  non-empty, otherwise we put $\tau(A,h)=0$. Notice that if $h>2\sqrt{d} $
  and $l_{min}(A)>\sqrt{d}$, then $A_1^h$ and $A_2^h$ are
  non-empty. Similarly, we define the variable $\phi(A,h)$ by:
$$ \phi(A,h) \,=\, \phi_t(\cyl(A,h)\cap \ZZ^d,B(A,h),T(A,h)) \;.$$
Finally, \emph{$p_c(d)$ denotes the critical parameter for
the Bernoulli bond percolation on $\mathbb{Z}^d$}.

\section{Background and main results}
\label{sec:results}
\subsection{Background}

The following result allows do define the flow constant
$\nu(\vec{v_0})$ when $\vec{v_0}\in\RR^d$ is the vector
$(0,\ldots,0,1)$. It follows from the subadditive
ergodic theorems of \cite{Ackoglu}, \cite{KrengelPyke} and \cite{Smythe}. Let $\mathbf{k}=(k_1,\ldots,k_{d-1})\in (\NN^*)^n$, and define
$A_{\mathbf{k}}=\prod_{i=1}^{d-1}[0,k_i] \times \{0\}$.
\begin{thm}[\cite{Kesten:flows}]
Suppose that $h(n)$ goes to infinity when $n$ goes to infinity, and
  that:
$$\int_0^\infty x\;dF(x)<\infty\;.$$
 Then,
$\tau(nA_{\mathbf{k}},h(n))/(n^{d-1}\prod_{i=1}^{d-1}k_i)$ converges almost surely  and
in $L^1$, when $n$ goes to infinity, to a non-negative,
finite constant $\nu(\vec{v_0})$ which does not depend on $\mathbf{k}$.
\end{thm}
An important problem is to know when $\nu(\vec{v_0})$ equals
zero. It is proved in \cite{Theret:small} (see also \cite{Chayes}
for capacities equal to zero or one) that $F(0)<1-p_c(d)$
implies $\nu(\vec{v_0}) >0$. Conversely, Zhang proved in \cite{Zhang}, Theorem 1.10 that
$\nu(\vec{v_0}) =0$ if $F(0) = 1-p_c(d)$, and so by a simple coupling of
probability if $F(0) \geq 1-p_c(d)$. Actually, he wrote the proof for $d=3$
but said himself that the argument works for $d\geq 3$ (see Remark 1 of
\cite{Zhang}). This property is also satisfied in dimension $d=2$ where we
can use duality arguments (see \cite{Kesten:StFlour} Theorem (6.1) and
Remark (6.2)). We gather these results in the following theorem:
\begin{thm}
\label{thm:criticalzhang}
Suppose that $\int_0^\infty x\;dF(x)$ is finite. Then, $\nu(\vec{v_0})=0$ if and only if $F(0)\geq 1-p_c(d)$.
\end{thm}
Finally, a crucial result is the following theorem of
Zhang, which allows to control the number of edges in a cut of minimal
capacity. Let $\mathbf{k}=(k_1,\ldots,k_{d-1})\in (\NN^*)^n$,
$m\in\NN^*$ and define:
$$\mathbf{B}(\mathbf{k},m)=\prod_{i=1}^{d-1}[0,k_i]\times[0,m]\;.$$
Let $N(\mathbf{k},m)$ be the number of edges of a cut $E$ between
$\mathbf{B}(\mathbf{k},m)$ and $\infty$ which achieves the minimal capacity
$V(E)=\sum_{e\in E}t(e)$ among all these cuts. If there are
more than one cut achieving the minimum, we use a deterministic method to
select a unique one with the minimum number of edges among these.
\begin{thm}[\cite{Zhang07}, Theorem 1]
\label{thmzhang}
If $F(0)<1-p_c(d)$, then there exists constants $\beta=\beta(F,d)$, $m_0(F,d)$ and
$C_i=C_i(F,d)$, for $i=1,2$ such that for all $n\geq \beta
\prod_{i=1}^{d-1}k_i$ and $m_0\leq m\leq \min_{i=1,\ldots ,d-1} k_i$,
$$\PP(N(\mathbf{k},m)\geq n)\leq C_1e^{-C_2n}\;.$$
\end{thm}
An analogue result is obtained in \cite{Zhang07}, Theorem 2, for the minimal cut
between the top and the bottom of $\mathbf{B}(\mathbf{k},m)$ inside
$\mathbf{B}(\mathbf{k},m)$. We shall make use of Theorem
\ref{thmzhang} through a slight modification, Proposition
\ref{prop:5.8Zhangbis} below.

Finally, Kesten proved in 1987 the law of large numbers for
$\phi$ in vertical boxes in dimension $3$ under the additional assumption
that $F(0)$ is sufficiently small and $h(n)$ not too large, plus an
assumption of finite exponential moment (see
Theorem 2.12 in \cite{Kesten:flows}). In a remarkable
paper, Zhang recently
improved Kesten's result by relaxing the assumption on $F(0)$ to the
relevant one $F(0)<1-p_c(d)$, and extended it to any dimension $d\geq 3$ (see
\cite{Zhang07}). Zhang proved the following result:
\begin{thm}[\cite{Zhang07}]
\label{cvphi}
Suppose $F(0)<1-p_c(d)$, and there exists $\gamma>0$ such that:
$$\int e^{\gamma x}\;dF(x)<\infty\;.$$
If $k_1$, ..., $k_{d-1}$, $m$ go to infinity in such a way that for some
$0<\eta\leq 1$, we have
$$ \log m \,\leq\, \max_{1\leq i\leq d-1} k_i^{1-\eta} \,,$$
then
$$ \lim_{k_1 ,..., k_{d-1}, m \rightarrow \infty} \frac{\phi(A_{\mathbf{k}}, m)}{
  k_1 \cdots k_{d-1}} \,=\, \nu(\vec{v_0}) \qquad a.s. \,\, and\,\, in
  \,\, L^1 \,.$$
\end{thm}
\subsection{Hypotheses on the distribution $F$ and the height $h$}

Here we gather and present the main hypotheses that we shall do on $F$
and on the height $h$. Notice that $\mathbf{
  (F5)}\Rightarrow\mathbf{ (F4)}\Rightarrow\mathbf{
  (F3)}\Rightarrow\mathbf{ (F2)}$ and $\mathbf{
  (H3)}\Rightarrow\mathbf{ (H2)}$.
{\begin{center}
\renewcommand{\arraystretch}{1.5}
\begin{tabular}{|l|l||l|}
\hline
\multicolumn{2}{|c||}{Hypotheses on the distribution} &
\multicolumn{1}{|c|}{Hypotheses on the height}\\
\hline
$\mathbf{(F1)}\quad F(0)<1-p_c(d)$&   & $\mathbf{ (H1)}\quad \lim_{n \rightarrow \infty}h(n)=+\infty$\\
$\mathbf{ (F2)}\quad \int_0^\infty
x\;dF(x)<\infty$&  $\mathbf{ (F4)}\quad \exists \gamma>0,\;\int_0^\infty e^{\gamma
  x}\;dF(x)<\infty $ & $\mathbf{ (H2)}\quad \lim_{n \rightarrow
  \infty}\frac{\log h(n)}{n^{d-1}}=0$ \\ 
$\mathbf{ (F3)}\quad \int_0^\infty x^{1+\frac{1}{d-1}}\;dF(x)<\infty$&$\mathbf{ (F5)}\quad \forall \gamma>0,\;\int_0^\infty e^{\gamma
  x}\;dF(x)<\infty$ & $\mathbf{
  (H3)}\quad \lim_{n \rightarrow \infty}\frac{h(n)}{n}=0$ \\
\hline 
\end{tabular}
\end{center}}
\vspace*{2mm}

The following table summarizes the needed hypotheses for the main
results presented in the next sections.  SLLN stands
for Strong Law of Large numbers, LDP
for Large Deviation Principle and R.F. for Rate Function (of the Large
deviation principles). 
{\begin{center}
\renewcommand{\arraystretch}{1.1}
\begin{tabular}{|c|c|c|c|c|c|c|c|c|}
\hline
&Existence &Positivity  & \multicolumn{2}{|c|}{SLLN
  for $\tau$}&\multicolumn{2}{|c|}{SLLN
  for }&LDP for $\tau$& LDP for \\
&of the R.F.  &of the R.F.  & \multicolumn{2}{|c|}{(and flat
  $\phi$)}&\multicolumn{2}{|c|}{straight $\phi$}& (and flat $\phi$)& 
straight $\phi$ \\
\cline{4-7}
& & &$0\in A$&$0\not\in A$&$0\in A$&$0\not\in A$& &\\
\hline
$\mathbf{(F1)}$& &$\times$  &&&&&$\times$&$\times$\\
$\mathbf{ (F2)}$& &$\times$  &$\times$ &$\times$ &$\times$ &$\times$ &$\times$&$\times$\\
$\mathbf{ (F3)}$& & &&$\times$ &&$\times$ &$\times$&$\times$ \\
$\mathbf{ (F4)}$& & &&&&&$\times$&$\times$\\
$\mathbf{ (F5)}$& & &&&&&$\times$&\\
$\mathbf{ (H1)}$& $\times$  &$\times$  &$\times$ &$\times$ &$\times$ &$\times$ &$\times$&$\times$\\
$\mathbf{ (H2)}$& & &&&$\times$ &$\times$ &&$\times$ \\
$\mathbf{ (H3)}$& & &($\times$) &($\times$) &&&($\times$) & \\
\hline 
\end{tabular}
\end{center}}

\vspace*{2mm}
Let us comment this table a bit. First, assumption $\mathbf{(H1)}$ is
not necessary to study the flows $\tau$ and $\phi$, but it is necessary to obtain a flow constant $\nu(\vec{v})$
which does not depend on the height $h$, and moreover it is natural when we
interpret our system as a model for porous media. All the other
assumptions are  optimal concerning  $\tau$ (cf. Remarks \ref{rem:condF2},
\ref{rem:F1PGDtau} and \ref{remcondF5}) except
perhaps $\mathbf{(F3)}$ (see Remark \ref{rem:condF3}). In addition, our
assumptions are also essentially optimal concerning $\phi$ (cf. Remarks
\ref{rem:F1PGDtau} and \ref{remcondH2}) except perhaps assumption
$\mathbf{(F4)}$ (see Remark \ref{rem:condF4phi}) and $\mathbf{(F3)}$,
for the same reason as for $\tau$. Finally, assumption
$\mathbf{(H3)}$, used to obtain results for cylinders which are not
straight, is certainly far from optimality (see Remark
\ref{rem:phiechoue}). This assumption gives results only for ``flat cylinders''.

\subsection{Results concerning $\tau$}

First, we will extend the definition of $\nu(.)$ in all directions. 
\begin{prop}[Definition of $\nu$]
\label{defnu} We suppose that $\mathbf{(F2)}$ and $\mathbf{(H1)}$
hold.  For every non-degenerate hyperrectangle $A$,
the limit
$$ \lim_{n\rightarrow \infty} \frac{\mathbb{E} (\tau(nA,h(n)))}{
  \mathcal{H}^{d-1}(nA)}  $$
exists and depends on the direction of $\vec{v}$, one of the two unit
  vectors orthogonal to $\hyp(A)$, and not on $A$ itself. We denote it by
  $\nu(\vec{v})$ (the dependence in $F$ and $d$ is implicit).
\end{prop}
\begin{rem}
\label{rem:condF2}
We chose to define simply the flow constant $\nu$ from the
convergence of the rescaled expectations. Having made this choice,
condition $\mathbf{(F2)}$ is necessary for the limit to be
finite. Indeed, for most orientations, there exists two vertices $x\in
A_1^h$ and $y\in A_2^h$ which are neighbours in $\ZZ^d$. Thus, the
corresponding edge must belong to any cutsets, and this implies that
the mean of $\tau(nA,h(n))$ is finite only if $\mathbf{(F2)}$
holds. Notice however that with some extra work, one could probably define a
flow constant without any moment condition as in
\cite{Kesten:StFlour}, section 2.
\end{rem}

The following Proposition states some basic properties of $\nu$, and
notably settles the question of its positivity.

\begin{prop}[Properties of $\nu$]
\label{thmnu}
Suppose that $\mathbf{(F2)}$ and $\mathbf{(H1)}$
hold. Let $\delta = \inf \{ \lambda \,|\, \mathbb{P} (t(e)
\leq \lambda) >0 \}$. Then,
\begin{itemize}
\item[(i)] for every unit vector $\vec{v}$, $\nu(\vec{v})\geq
\delta \|\vec{v}\|_1$. 
\item[(ii)] if $F(\delta) <1-p_c(d)$, then $\nu(\vec{v}) > \delta \|\vec{v}\|_1$ for
all unit vector $\vec{v}$. In the case $\delta =0$, the previous
implication is in fact an equivalence.
\item[(iii)] for every unit vector $\vec{v}$, and every non-degenerate
  hyperrectangle $A$ orthogonal to $\vec{v}$,
$$\nu(\vec{v})\leq
\inf_{n\in\NN}\left\{\frac{\EE(t(e))K(d,A)}{n}+\frac{\EE(\tau(nA,h(n)))}{\mathcal{H}^{d-1}(nA)}\right\}\;,$$
where $K(d,A)=c(d) \mathcal{H}^{d-2}(\partial A)/\mathcal{H}^{d-1}(A)$, and
$c(d)$ is a constant depending only on the dimension $d$.
\end{itemize}
\end{prop}

We will derive the law of large numbers for $\tau(A,h)$ in big
cylinders $\cyl (A,h)$ as a consequence of an almost subadditive argument:

\begin{thm}[LLN for $\tau$]
\label{thm:llntau}
We suppose that $\mathbf{(F2)}$ and $\mathbf{(H1)}$
hold. Then,
$$ \lim_{n\rightarrow \infty}\frac{\tau(nA,h(n))}{ \mathcal{H}^{d-1}(nA)}
\,=\, \nu(\vec{v})  \qquad \mbox{ in }L^1\,. $$
Moreover, if $0\in A$, where $0$ denotes the origin of $\ZZ^d$, or if $\mathbf{(F3)}$ holds,
$$ \lim_{n\rightarrow \infty}\frac{\tau(nA,h(n))}{ \mathcal{H}^{d-1}(nA)}
\,=\, \nu(\vec{v})  \qquad a.s. $$
\end{thm}

Propositions \ref{defnu}, \ref{thmnu} and Theorem \ref{thm:llntau} will be proven in
section \ref{convergence}. Concerning large deviations, we will show two results: the first gives the speed of decay of the
probability that the rescaled flow $\tau$ is abnormally small, and
the second one states a large deviation principle for the rescaled variable
$\tau$.

The estimate of lower large deviations is the following. Notice that Theorem \ref{thmzhang} is the key to obtain the relevant condition
$F(0)<1-p_c(d)$.
\begin{thm}[Lower deviations for $\tau$]
\label{devinf}
Suppose that $\mathbf{(F1)}$, $\mathbf{(F2)}$ and $\mathbf{(H1)}$ hold. Then for every $\varepsilon >0$ there exists a
positive 
constant $C(d,F,\varepsilon)$ such that for
every unit vector $\vec{v}$ and every non-degenerate hyperrectangle $A$ orthogonal to
$\vec{v}$, there exists a constant $\widetilde{C} (d,F,A,\varepsilon)$ (possibly depending on
all the parameters $d$, $F$, $A$, $\varepsilon$) such that:
$$ \mathbb{P} \left( \frac{\tau(nA, h(n))}{\mathcal{H}^{d-1}(nA)} \leq 
\nu(\vec{v}) -\varepsilon \right) \,\leq\,\widetilde{C}(d,F,A,\varepsilon) \exp \left(-C(d,F,\varepsilon)
  \mathcal{H}^{d-1}(A) n^{d-1}  \right) \,.$$
\end{thm}
 Now we can state a large deviation principle:
\begin{thm}[LDP for $\tau$]
\label{thmpgd}
Suppose that $\mathbf{(F1)}$, $\mathbf{(F5)}$ and $\mathbf{(H1)}$  hold. Then for
every unit vector $\vec{v}$ and every non-degenerate hyperrectangle $A$ orthogonal to
$\vec{v}$, the sequence
$$ \left( \frac{\tau(nA, h(n))}{\mathcal{H}^{d-1}(nA)}, n\in \mathbb{N}
\right) $$
satisfies a large deviation principle of speed $\mathcal{H}^{d-1}(nA)$ with
the good rate function $\mathcal{J}_{\vec{v}}$. Moreover we know that
$\mathcal{J}_{\vec{v}}$ is convex on $\mathbb{R}^+$, infinite on $[0,
  \delta \|\vec{v}\|_1[\cup ]\nu(\vec{v}), +\infty[$,
    where $\delta = \inf \{ \lambda \,|\, \mathbb{P} (t(e)
\leq \lambda) >0 \}$, equal to $0$ at $\nu(\vec{v})$, and if $\delta
\|\vec{v}\|_1 < \nu(\vec{v})$ we also know that $\mathcal{J}_{\vec{v}}$ is
finite on $]\delta \|\vec{v}\|_1,
  \nu(\vec{v})]$, continuous and strictly decreasing on $[\delta \|\vec{v}\|_1,
  \nu(\vec{v})]$ and strictlypositive on $[\delta \|\vec{v}\|_1,
  \nu(\vec{v})[$.
\end{thm}
\begin{rem}
\label{rem:F1PGDtau}
Notice that, from Proposition \ref{thmnu}, assumption $\mathbf{(F1)}$
is necessary to have positive asymptotic rescaled maximal flow, and thus to give a sense to the study of lower large
deviations. Moreover, Theorem \ref{thmpgd} is interesting only if $\nu(\vec{v}) > \delta
\|\vec{v}\|_1$. Proposition \ref{thmnu} states that it is the case at least
if $F(\delta) <1-p_c(d)$, and in the case $\delta =0$, this condition
is optimal. We do not know the optimal condition on $F(\delta)$ when
$\delta\not =0$.
\end{rem}

\begin{rem}
In his PhD-thesis \cite{Wouts}, section 2, Wouts shows a similar lower large deviations
result in the context of the dilute Ising model. More precisely, for every temperature $T$, a Gibbs measure $\Phi_{n,T}$ with i.i.d. nonnegative, bounded random
interactions $(J_e)_{e\in\EE^d}$ is constructed on the set of configurations $\{0,1\}^{E_n}$, where $E_n$ is the set of edges of a cube $B_n$ of length $n$, and 0 (resp. 1) means the edge is closed (resp. open). Wouts defines the quenched surface tension in this box as the normalized logarithm of the $\Phi_{n,T}$-probability of the event that there is a disconnection between the upper and lower parts of the boundary of $B_n$. Then, Wouts shows that for Lebesgue-almost every
temperature $T$, the  quenched surface tension satisfies a large deviation principle at surface order. A remarkable feature of this work is that the proof, quite simple, relies on a concentration property that avoids the use of any estimate like that of Theorem \ref{thmzhang}. A similar treatment could be done in our setting, with the value of $F(0)$ playing the role of the inverse temperature. Of course, this is quite artificial and unsatisfactory for our purpose, since one would not obtain any information for a precise distribution function $F$, but rather for almost all distributions of the form $p\delta_0+(1-p)dF$, $p\in [0,1]$. Still, it seems to us that Wouts' method deserves further investigation.
\end{rem}


\subsection{Results concerning $\phi$ in flat cylinders}
Under the additional assumption that the cylinder we study is sufficiently
flat, in the sense that we suppose $\lim_{n\rightarrow \infty}
h(n)/n =0$, we can transport results from $\tau$ to $\phi$ even in
non-straight boxes,
because the behaviour of these two variables are very similar in that
case. We obtain the following two results:
\begin{thm}[Lower deviations for flat $\phi$]
\label{cordevinf}
Suppose $\mathbf{(F1)}$, $\mathbf{(F2)}$, $\mathbf{(H1)}$ and $\mathbf{(H3)}$  hold. Then for every $\varepsilon >0$ there exists a positive
constant $C'(d,F,\varepsilon)$ such that for every unit vector
$\vec{v}$, every non-degenerate hyperrectangle $A$ orthogonal to $\vec{v}$, there exists a constant $\widetilde{C}'
(d,F,A,h,\varepsilon)$ (possibly depending on
all the parameters $d$, $F$, $A$, $h$, $\varepsilon$) such that
$$ \mathbb{P} \left( \frac{\phi(nA,
  h(n))}{\mathcal{H}^{d-1}(nA)} \leq 
\nu(\vec{v}) -\varepsilon \right) \,\leq\,\widetilde{C}'
(d,F,A,h,\varepsilon)  \exp \left(-C'(d,F,\varepsilon)
  \mathcal{H}^{d-1}(A) n^{d-1}  \right) \,.$$
\end{thm}
\begin{cor}[of Theorem \ref{thmpgd}, LDP for flat $\phi$]
\label{corthmpgd}
Suppose  $\mathbf{(F1)}$, $\mathbf{(F5)}$, $\mathbf{(H1)}$ and $\mathbf{(H3)}$ hold. Then, for
every unit vector $\vec{v}$ and every non-degenerate hyperrectangle $A$ orthogonal to
$\vec{v}$, the sequence
$$ \left( \frac{\phi(nA, h(n))}{\mathcal{H}^{d-1}(nA)}, n\in \mathbb{N}
\right) $$
satisfies a large deviation principle of speed $\mathcal{H}^{d-1}(nA)$ with
the good rate function $\mathcal{J}_{\vec{v}}$ (the same as in Theorem
\ref{thmpgd}).
\end{cor}
\begin{rem}
Theorem \ref{cordevinf} will be proven exactly as Theorem \ref{devinf},
using the fact that the convergence of $\EE[\tau(nA, h(n))]/\H^{d-1}(nA)$
implies the convergence of $\EE [\phi(nA, h(n))]/\H^{d-1}(nA)$ under the
hypotheses $\mathbf{(F2)}$ and $\mathbf{(H3)}$. Corollary \ref{corthmpgd}
will be proven using the exponential equivalence of the rescaled variables
$\tau(nA, h(n))$ and $\phi(nA, h(n)$ under hypotheses $\mathbf{(F5)}$ and
$\mathbf{(H3)}$.
\end{rem}


\subsection[Results concerning $\phi$ in straight cylinders]{Results
    concerning $\phi$ in straight but high cylinders}
We shall say that a hyperrectangle $A$ is
 \emph{straight} if it is of the
 form $\prod_{i=1}^{d-1} [0,a_i] \times \{0\}$ ($a_i \in \mathbb{R}^+_*$ for
 all $i$, so a straight hyperrectangle is non-degenerate). In particular,
 Theorem \ref{cvphi} implies that for a straight hyperrectangle $A$, 
for every function $h:\mathbb{N}\rightarrow
\mathbb{R}^+$ satisfying $\lim_{n\rightarrow \infty } h(n) =+\infty$ and
$\log h(n) \leq n^{1-\eta}$ for some $0<\eta\leq 1$, we have
$$ \lim_{n\rightarrow \infty} \frac{\phi(nA,h(n))}{\mathcal{H}^{d-1}(nA)}
\,=\, \nu((0,...,0,1))\qquad a.s.\,\, and\,\, in \,\, L^1 \,.  $$
We obtain three results for the rescaled variable $\phi$ in straight
cylinders. Using subadditivity and symmetry arguments, we can prove the law of large
 numbers for $\phi$ in straight boxes under a minimal moment
 condition, and the hypothesis $(\mathbf{H2})$ on $h$:
\begin{thm}[LLN for straight $\phi$]
\label{thm:llnphistraight}
Suppose that  $\mathbf{(F2)}$, $\mathbf{(H1)}$ and
$\mathbf{(H2)}$ hold, and that $A$ is a straight hyperrectangle. Then,
$$\lim_{n\rightarrow\infty}\frac{\phi
  (nA,h(n))}{\H^{d-1}(nA)}= \nu(\vec{v_0})\quad a.s.\mbox{ and in }L^1$$
where $\vec{v_0}=(0,\ldots,0,1)$.  
\end{thm} 
Under the additional condition of an exponential moment for
$F$, we can prove a large deviation
principle for $\phi$ in straight boxes.
\begin{thm}[LDP for straight $\phi$]
\label{thmpgdphi}
Suppose $\mathbf{(F1)}$, $\mathbf{(F4)}$, $\mathbf{(H1)}$ and $\mathbf{(H2)}$ hold. Then for
 every straight hyperrectangle $A$, the sequence
$$ \left( \frac{\phi(nA, h(n))}{\mathcal{H}^{d-1}(nA)}, n\in \mathbb{N}
\right) $$
satisfies a large deviation principle of speed $\mathcal{H}^{d-1}(nA)$ with
the good rate function $\mathcal{J}_{\vec{v}}$ with $\vec{v}=(0,...,0,1)$ (the same as in Theorem
\ref{thmpgd}).
\end{thm}
We also obtain a result similar to Theorem \ref{devinf} for $\phi$:
\begin{thm}[Lower deviations for straight $\phi$]
\label{devinfphi}
Suppose $\mathbf{(F1)}$, $\mathbf{(F2)}$, $\mathbf{(H1)}$ and $\mathbf{(H2)}$ hold.  Then, for every $\varepsilon >0$ there exists a positive
constant $C''(d,F,\varepsilon)$ such that for every straight hyperrectangle $A$, there exists a
strictly positive constant $\widetilde{C}''(d,F,A,h,\varepsilon)$ (possibly depending on
all the parameters $d$, $F$, $A$,$h$ and $\varepsilon$) such that:
$$
\mathbb{P} \left( \frac{\phi(nA,
  h(n))}{\mathcal{H}^{d-1}(nA)} \leq 
\nu((0,...,0,1)) -\varepsilon \right) \,\leq\,
\widetilde{C}''(d,F,A,h,\varepsilon) \exp \left(-C''(d,F,\varepsilon)   \mathcal{H}^{d-1}(A) n^{d-1}  \right) \,.
 $$
\end{thm}
This result answers question (2.25) in \cite{Kesten:flows}. We have to comment these three theorems by some remarks.
\begin{rem}
We decided to state the law of large numbers (Theorem
\ref{thm:llnphistraight}) in the case were the origin of the graph belongs
to the straight hyperrectangle $A$ since it is the case in the literature
(see \cite{Kesten:flows}, \cite{Zhang07}). We also could state the same
result for a hyperrectangle $A$ of the form
$\prod_{i=1}^{d-1}[a_i,b_i]\times \{c\}$ for real numbers $a_i<b_i$ and
$c$. In this case, exactly as in Theorem \ref{thm:llntau}, the same hypotheses
$\mathbf{(F2)}$, $\mathbf{(H1)}$ and $\mathbf{(H2)}$ are required to obtain
the convergence  of $\phi(nA, h(n))/\H^{d-1}(nA)$ in $L^1$, but
we need moreover the stronger hypothesis $\mathbf{(F3)}$ to obtain the
a.s. convergence of the variable if the origin of the graph does not belong
to $A$. 
\end{rem}

\begin{rem}
The proofs of these three theorems are a little bit tangled. It comes from
our willingness to obtain the best hypotheses on $F$ each time. Indeed, we
stress the fact that Theorem \ref{devinfphi} is not a simple consequence of
Theorem \ref{thmpgdphi} when $\mathbf{(F4)}$ does not hold. In fact, we will prove first a proposition,
Proposition \ref{prop:fntauxphi}, that will
lead to Theorem \ref{thm:llnphistraight} and Theorem \ref{thmpgdphi}
independently. Theorem \ref{devinfphi} will be proven exactly as Theorems
\ref{devinf} and \ref{cordevinf}, using Theorem \ref{thm:llnphistraight}.
\end{rem}

\begin{rem}
\label{remcondH2}
Actually the condition $\mathbf{(H2)}$, i.e. $\lim_{n\rightarrow
  \infty} \log h(n) / n^{d-1} =0$, is essentially the good one. For instance, if
$A=[0,1]^{d-1}\times\{0\}$, $h(n) \geq \exp (k
n^{d-1})$ for a constant $k$ sufficiently large and $F(0)>0$, then the
maximal flow $\phi(nA,h(n))$ eventually equals $0$, almost surely. Indeed if the $n^{d-1}$ vertical edges of the cylinder that intersect one fixed horizontal plane have all $0$ for capacity then $\phi(nA,h(n))=0$. By independence and translation invariance of the model, we obtain:
$$ \mathbb{P} \left[ \phi(nA,h(n)) \neq 0 \right] \,\leq \, \left[
  1-F(0)^{n^{d-1}} \right]^{2\exp (k
n^{d-1})} \, $$
which is summable for $k$ large enough, and so we conclude by the
Borel-Cantelli lemma.
\end{rem}

\begin{rem}
Notice that our setting in Theorem \ref{thm:llnphistraight} is not entirely similar to the one of
\cite{Zhang07} since each side of $nA$ grows at the
same speed, whereas Zhang considers $A=\prod_{i=1}^{d-1}[0,k_i]\times\{0\}$ and
lets all the $k_i$ go to infinity, possibly at different speeds. In
  the case we consider, we improve the height and moment conditions in Theorem 3 of
  \cite{Zhang07} to the relevant one, and so partially answer the question
  contained in Remark (2.17) and question (2.24) in
  \cite{Kesten:flows}. See also Remark \ref{rem:improveZhang}.
\end{rem}





\section{Lower large deviations for $\tau$ and $\phi$ and law of large
  numbers for $\tau$}
\label{sec:lowerdev}
In section \ref{sec:concentration}, we derive the crucial deviation
inequalities from their means of the flows $\tau$ and $\phi$. This will lead to
the law of large numbers for $\tau$ rescaled in section \ref{subsec:llntau}, and the
deviations from $\nu$ of $\tau$ and flat $\phi$ rescaled in section
\ref{sectiondevinf}. Of course, we need to define properly $\nu$ in
any direction, and this is done in section \ref{convergence}, whereas properties
of $\nu$ are proven in section \ref{subsec:thmnu}, using a combinatorial
result stated in section \ref{sec:minimalsize}.

\subsection{Minimal size of a cutset}
\label{sec:minimalsize}

For every hyperrectangle $A$, we denote by $\mathcal{N}(A,h)$ the minimal
number of edges in $A$ that can disconnect $A_1^h$ from $A_2^h$ in
$\cyl(A,h)$, if $A_1^h$ and $A_2^h$ are non-empty. The following lemma gives the asymptotic order of
$\mathcal{N}(nA,h(n))$ when $n$ goes to infinity.
\begin{lem}
\label{plaquettes}
Let $\vec{v}$ be a unitary vector. Then for all hyperrectangle $A$ 
orthogonal to $\vec{v}$, for all function $h: \mathbb{N} \rightarrow
]2\sqrt{d},+\infty[$, and for every $n\in \NN$ such that $l_{min}(nA)>\sqrt{d}$,
$$\left\arrowvert\, \frac{\mathcal{N}(nA,h(n))}{\mathcal{H}^{d-1}(nA)} - \|
\vec{v}  \|_1 \,\right\arrowvert \,\leq\,  \frac{d \mathcal{H}^{d-2}(\partial A)}{n \mathcal{H}^{d-1}(A)} \,.$$
\end{lem}

\begin{dem}
We introduce some definitions. For $A$ a hyperrectangle
orthogonal to $\vec{v}$, we
denote by $P_i(A)$ the orthogonal projection of $A$ on the $i$-th
hyperplane of
coordinates, i.e., the hyperplane $\{(x_1,...,x_d) \in \mathbb{R}^d \,|\,
x_i =0 \}$. We have the property
$$ \frac{\sum_{i=1}^d \mathcal{H}^{d-1} (P_i(A))}{\mathcal{H}^{d-1}(A)}
\,=\, \| \vec{v} \|_1 \,.$$
Indeed, $\mathcal{H}^{d-1}(P_i(A)) = |v_i| \mathcal{H}^{d-1}(A)$, where
$\vec{v} = (v_1,...,v_d)$. We define now $E_i (nA)$ the set of edges orthogonal to the
$i$-th hyperplane of coordinates that `intersect' the
hyperrectangle $nA$ in the following sense:
$$ E_i(nA) \,=\, \{ e =\langle x,y \rangle \in \mathbb{E}^d \,|\, y_i = x_i
+1  \,\, and \,\, [x,y[ \,\cap\, nA \neq \emptyset\,\, and \,\, [x,y[
\not\subset nA\}\,.  $$
We exclude here the extremity $y$ in the segment $[x,y[$ to avoid problems
    of non uniqueness of such an edge intersecting $nA$ at a given point. On
    one hand, we have a straight path that goes from $(nA)_1^{h(n)}$ to
    $(nA)_2^{h(n)}$ through each edge of $E_i(nA)$, $i=1,...,d$,
    except maybe the edges that intersect $nA$ along $\partial (nA)$, and these
    paths are disjoint, so a set of edges that disconnect $(nA)_1^{h(n)}$
    from $(nA)_2^{h(n)}$ in $\cyl(nA,h(n))$ must cut each one of these
    paths, thus
$$ \mathcal{N}(nA,h(n)) \,\geq \, \sum_{i=1}^d  \left(
    \mathcal{H}^{d-1}(P_i(nA)) - \mathcal{H}^{d-2}(\partial P_i(nA)) \right) \,\geq\,
    \left( \|\vec{v} \|_1 - d
    \frac{\mathcal{H}^{d-2}(\partial(nA))}{\mathcal{H}^{d-1}(nA)} \right)
    \mathcal{H}^{d-1}(nA) \,.$$
On the other hand, each path from $(nA)_1^{h(n)}$ to $(nA)_2^{h(n)}$ in
    $\cyl(nA,h(n))$ must go through $nA$ and so contains an edge of one
    of the $E_i(nA)$, $i=1,...,d$. It suffices then to remove all the edges
    in the union of the sets  $E_i(nA),i=1,\ldots ,d$ to disconnect $(nA)_1^{h(n)}$
    from $(nA)_2^{h(n)}$ in $\cyl(nA,h(n))$, and so
$$\mathcal{N}(nA,h(n)) \,\leq\,  \sum_{i=1}^d  \left(
    \mathcal{H}^{d-1}(P_i(nA)) + \mathcal{H}^{d-2}(\partial P_i(nA)) \right)
    \,\leq\, \left( \|\vec{v} \|_1 + d
    \frac{\mathcal{H}^{d-2}(\partial(nA))}{\mathcal{H}^{d-1}(nA)} \right)
    \mathcal{H}^{d-1}(nA) \,.  $$
We conclude that
$$ \left\arrowvert\, \frac{\mathcal{N}(nA,h(n))}{\mathcal{H}^{d-1}(nA)} - \|
\vec{v}  \|_1 \,\right\arrowvert \,\leq\,
d\,\frac{\mathcal{H}^{d-2}(\partial(nA))}{\mathcal{H}^{d-1}(nA)} \,=\,
\frac{d \mathcal{H}^{d-2}(\partial A)}{n \mathcal{H}^{d-1}(A)} \,.  $$
\end{dem}
\subsection{Lower deviations of the maximal flows from their means}
\label{sec:concentration}

Let $A$ be a non-degenerate hyperrectangle. In this section, we obtain
deviation inequalities for $\phi(A,h)$ and $\tau(A,h)$ from their
means. These inequalities, stated below in Proposition \ref{prop:deviation}, give the right speed for the lower large deviation probabilities  as
  soon as the convergence of the rescaled expectation of the variables is
  known. This will be used in section \ref{convergence} to prove the law
of large numbers for $\tau$, but above all this will be essential to show
the positivity of the rate function for lower large deviations in section
\ref{positivity}. This positivity will be used to prove Theorem
\ref{thm:llnphistraight} in section \ref{sec:llnphistraight}.

To get this result, we state below in Proposition \ref{prop:5.8Zhangbis} a slight modification of Zhang's Theorem \ref{thmzhang}, which allows to control the number of edges in a cut of minimal
capacity.  Notice that in this precise form, Proposition \ref{prop:5.8Zhangbis} is almost a strict analogue for flow
problems of Proposition 5.8 in \cite{Kesten:StFlour}, the latter being of
utmost importance in the study of First Passage Percolation. 

We introduce the following notation: $E_{\tau(A,h)}$ (resp. $E_\phi(A,h)$)
is a cut whose capacity achieves
the minimum in the dual definition (\ref{eq:maxflowmincut}) of
$\tau(A,h)$ (resp. $\phi(A,h)$). If there are
more than one cut achieving the minimum, we use a deterministic method to
select a unique one with the minimum number of edges among
these. Recall also that for a hyperrectangle $A$, we defined
$l_{min}(A)$ as the ``smallest length of $A$'', i.e. the number $t$
such that $A$ is the isometric image of
$\prod_{i=1}^{d-1}[0,t_i]\times\{0\}$, with $t=t_1\leq\ldots\leq t_{d-1}$.
\begin{prop}
\label{prop:5.8Zhangbis}
Suppose that $\mathbf{(F1)}$ holds, i.e. $F(0)<1-p_c(d)$. Then, there are constants $\eps(F,d)$,
$C_1(F(0),d)$, $C_2(F(0),d)$ and $t_0(F(0),d)$, such that, for every
$s\in\RR$, every non-degenerate hyperrectangle $A$ such that
$l_{min}(A)\geq t_0$, and every $h>2\sqrt{d}$, we have:
$$\PP\left( \card (E_{\tau(A,h)}) \geq s\mbox{ and }  \tau(A,h)\leq \eps (F,d) s
\right) \leq C_1(F(0),d) e^{-C_2(F(0),d)s}\;,$$
and:
$$\PP\left( \card (E_{\phi(A,h)}) \geq s\mbox{ and }  \phi(A,h)\leq \eps (F,d) s \right) \leq
C_1 (F(0),d) h e^{-C_2 (F(0),d) s}\;.$$
Furthermore, the constant $\eps$ depends on $F$ only on the neighborhood of 0
in the sense that if $(F_n)_{n\in\NN}$ is a sequence of possible distribution
functions for $t(e)$, which coincide on $[0,\eta]$ for some $\eta>0$, then
one can take the same constants $\eps$, $C_1$ and $C_2$ for the whole
sequence in the above inequalities.
\end{prop}

\begin{dem}
First notice that when $d=2$, this is a consequence of Proposition~5.8
in \cite{Kesten:StFlour}, through duality. In fact, when $d=3$, all the hard work has been done by Zhang, giving Theorem \ref{thmzhang}, so we only stress the
minor differences for the reader who would like to check how one goes
from the proof of Theorem \ref{thmzhang} (i.e. Theorem 1 in \cite{Zhang07}) to
Proposition \ref{prop:5.8Zhangbis}, and we rely heavily on the proof
and notations of \cite{Zhang07}.

The first thing is to see that one can perform the renormalization
argument of section 2 of \cite{Zhang07}. To do this for $\tau(A,h)$, replace $\infty$ by $A_2^{h}$ and the box $\mathbf{B}(\mathbf{k},m)$ by
$A_1^{h}$. For $\phi(A,h)$, replace $\infty$ by $T(A,h)$ and
$\mathbf{B}(\mathbf{k},m)$ by $B(A,h)$. For both $\tau(A,h)$ and
$\phi(A,h)$ also, one requires that all the connectedness properties
happen "in $\cyl(A,h)$". Then, the construction of the linear cutset
is identical, except for one thing: when $B_t(u)$ is a block of the
"block cutset" such that
$\overline{B}_t(u)$ intersects $\partial\cyl(A,h)$, it has a property
slightly different than the "blocked property" of Zhang. Define $\overline{B}'_t(u)$ to be the set of $t$-cubes which
are $\mathbf{L}^d$-neighbours of the cubes in $\overline{B}_t(u)$. Let
us say that a set of vertices $V_0$ of $\ZZ^d$ is of \emph{smallest length} $t$ if there
is a hyperrectangle $H$ in $\RR^d$, isometric image of
$[0,t]^{d-1}\times\{0\}$, such that for each edge $e$ of $\ZZ^d$
intersecting $H$, there is an endpoint of $e$ which belongs to
$V$. 
Now, let us say that a block  $B_t(u)$ has a "blocking surface
property" if either one of the following holds:
\begin{itemize}
\item[(i)] there are two subsets of vertices $V_1$ and $V_2$ of
smallest length $t/2$  in $\overline{B}'_t(u)$ which cannot be
connected by an open path in $\overline{B}'_t(u)$,
\item[(ii)] or there are a subset of vertices $V_1$ of
smallest length $t/2$ and an open path $\gamma$ connecting $B_t(u)$
to  $\overline{B}_t(u)$ in $\overline{B}'_t(u)$ such that $\gamma$ and
$V_1$  cannot be
connected by an open path in $\overline{B}'_t(u)$.
\end{itemize}
Then, if $A$ is of \emph{smallest length} larger than $t$, and
if $B_t(u)$ is a block of the
"block cutset" such that $\overline{B}_t(u)$ intersects
$\partial\cyl(A,h)$, it has a "blocking surface property". Now, it is easy to see, using the same arguments as Zhang from \cite{Grimmett99}, section 7, that the
probability that $\overline{B}_t(u)$ has a "blocking surface property"
decays exponentially to zero as $t$ goes to infinity, when
$F(0)<1-p_c(d)$.  This shows that the renormalization works if $A$ is
of smallest length larger than some $t_0(F(0),d)$, see the choice of $t$ above (5.26) in
\cite{Zhang07}. Notice that to prove Lemma 8 in \cite{Zhang07}, Zhang appeals to Lemma 7.104 in \cite{Grimmett99} whereas it
seems better to see this as a direct consequence of the fact that
percolation in slabs occurs.

The rest of the proof is almost identical. Note however that when
considering $\tau(A,h)$, there is no need to put a sum over the
possible intersections of the cutset with $\mathbf{L}$ (in (5.4), and
before (5.26)), since we know there is a constant $R(d)$ such that
there is a set of $R(d)$ edges that a cut needs to intersect (it is essentially
``pinned'' at the border of $A$). This is why we do not have any
condition on the height in the first inequality of Proposition
\ref{prop:5.8Zhangbis}, and why on the contrary $h$ appears in our
second inequality: for $\phi(A,h)$, we only know a set of $h$ edges
that a cut needs to intersect.

Finally, notice that we do not have any condition of moment on $F$,
since we are bounding the probability that $\{\card (E_{\tau(A,h)})\geq k\}$
\emph{and} $\{\tau(A,h)\leq \eps k \}$ occur, not only $ \PP\left( \card
  (E_{\tau(A,h)}) \geq k\right)$,
  and Zhang uses the moment condition only to bound
  $\PP\left(\tau(A,h)\leq \eps k \right)$. Also, the last statement on
  the constants is easily seen by tracking the choice of $\eps$ (see
  (5.2) and below (5.10)).
\end{dem}

Thanks to Proposition \ref{prop:5.8Zhangbis} and general deviation
 inequalities due to \cite{Boucheronetal03}, we obtain the
following deviation result for the maximal flows $\tau(nA,h(n))$ and
$\phi(nA,h(n))$.

\begin{prop}
\label{prop:deviation}
Suppose that hypotheses $\mathbf{(F1)}$ and $\mathbf{(F2)}$ occur. Then,
for any $\eta\in]0,1]$, there are positive constants $C(\eta,F,d)$,
$C_3(F(0),d)$ and $t_0(F,d)$ such that, for every $n\in\NN^*$, every non-degenerate hyperrectangle
$A$ such that $nA$ has smallest length at least $t_0$:
\begin{equation}
\label{eq:dev_tau_ordre1}
\PP\left( \tau(nA,h(n)) \leq \EE(\tau(nA,h(n)))(1-\eta) \right) \leq
C_3(F(0),d)e^{-C(\eta,F,d)\EE(\tau(nA,h(n)))}\;,
\end{equation}
and:
\begin{equation}
\label{eq:dev_phi_ordre1}
\PP\left( \phi(nA,h(n)) \leq \EE(\phi(nA,h(n)))(1-\eta) \right) \leq C_3(F(0),d)h(n)e^{-C(\eta,F,d)\EE(\phi(nA,h(n)))}\;.
\end{equation}
\end{prop}
\begin{dem}
To shorten the notations, define $\tau_n=\tau(nA,h(n))$ and
$\phi_n=\phi(nA,h(n))$.
We prove the result for $\tau_n$, the variant for $\phi_n$ being
entirely similar. Since
$\PP\left( \tau_n\leq \EE(\tau_n)(1-\eta) \right)$ is a decreasing
function of $\eta$, it is enough to prove the result for all $\eta$
less or equal to some absolute $\eta_0\in]0,1[$. We use this remark to
exclude the case $\eta =1$ in our study, thus, from now on,
let $\eta$ be a fixed real number in $]0,1[$.

Fix $A$ a non-degenerate hyperrectangle, and $n$ such that $nA$ has smallest length
at least $t_0(F,d)$, with $t_0$ as in Proposition~\ref{prop:5.8Zhangbis}. We order the edges in $\cyl(nA,h(n))$ as
$e_1,\ldots,e_{m_n}$. For every hyperrectangle $A$, we denote by
$\mathcal{N}(A,h)$ the minimal number of edges in $A$ that can
disconnect $A_1^h$ from $A_2^h$ in $\cyl(A,h)$, as in section \ref{sec:minimalsize}. For any real number $r\geq \mathcal{N}(nA,h(n))$, we define:
$$ \tau_n^r\,=\, \min\left\{ \begin{array}{c} V(E)\mbox{ s.t. }\card(E)\leq
r\mbox{ and } E \mbox{ cuts }\\
(nA)_1^{h(n)}\mbox{
  from }(nA)_2^{h(n)}\mbox{ in }\cyl(nA,h(n)) \end{array} \right\}\;. $$
Now, suppose that $F(0)<1-p_c(d)$, let $\eps$, $C_1$ and $C_2$ be as
in Proposition~\ref{prop:5.8Zhangbis}, and define
$r=(1-\eta)\EE(\tau_n)/\eps$.  Suppose
first that $r<\mathcal{N}(nA,h(n))$. Then,
\begin{eqnarray*}
\PP(\tau_n\leq (1-\eta)\EE(\tau_n))&=&\PP(\tau_n\leq
(1-\eta)\EE(\tau_n)\mbox{ and } \card (E_{\tau_n}) \geq
(1-\eta)\EE(\tau_n)/\eps)\;,\\
&\leq&C_1e^{-C_2(1-\eta)\EE(\tau_n)/\eps}\;,
\end{eqnarray*}
from Proposition~\ref{prop:5.8Zhangbis}, and the desired inequality is
obtained. Suppose now that $r\geq\mathcal{N}(nA,h(n))$. Then,
\begin{eqnarray}
\nonumber \PP(\tau_n\leq (1-\eta)\EE(\tau_n))&=&\PP(\tau_n\leq
(1-\eta)\EE(\tau_n)\mbox{ and }\tau_n^r\not=\tau_n)+\PP(\tau_n^r\leq (1-\eta)\EE(\tau_n))\;,\\
\label{eq:tautaunr} &\leq&C_1e^{-C_2r}+\PP(\tau_n^r\leq (1-\eta)\EE(\tau_n^r))\;,
\end{eqnarray}
from Proposition \ref{prop:5.8Zhangbis} and the fact that $\tau_n^r\leq
\tau_n$. Now, we truncate our variables $t(e)$. Let $a$ be a positive
real number to be chosen later, and define $\tilde t(e)=t(e)\land
a$. Let:
$$ \tilde \tau_n^r\,=\, \min\left\{ \begin{array}{c} \sum_{e\in
      E}\tilde t(e)\mbox{ s.t. }\card(E)\leq r\mbox{ and }E \mbox{ cuts }\\
 (nA)_1^{h(n)}\mbox{
  from }(nA)_2^{h(n)}\mbox{ in }\cyl(nA,h(n)) \end{array} \right\}\;. $$
Notice that $\tilde \tau_n^r\leq \tau_n^r$. We shall denote by
$R_{\tilde\tau_n^r}$ the intersection of all the cuts whose capacity achieves the minimum in the
definition of $\tilde \tau_n^r$. Then,
\begin{eqnarray*}
0\leq \EE(\tau_n^r)-\EE(\tilde\tau_n^r)&\leq &\EE\left[\sum_{e\in
  R_{\tilde\tau_n^r}}t(e)-\sum_{e\in
  R_{\tilde\tau_n^r}}\tilde t(e)\right]\;,\\
&\leq &\EE\left[\sum_{e\in
  R_{\tilde\tau_n^r}}t(e)\II_{t(e)\geq a}\right]\;,\\
&=&\sum_{i=1}^{m_n}\EE(t(e_i)\II_{t(e_i)\geq a}\II_{e_i\in
  R_{\tilde\tau_n^r}})\;,\\
&=&\sum_{i=1}^{m_n}\EE\left\lbrack \EE\left(t(e_i)\II_{t(e_i)\geq a}\II_{e_i\in
  R_{\tilde\tau_n^r}}|(t(e_j))_{j\not= i}\right)\right\rbrack\;.
\end{eqnarray*}
Now, when $(t(e_j))_{j\not= i}$ is fixed, $t(e_i)\mapsto\II_{e_i\in
  R_{\tilde\tau_n^r}}$ is a non-increasing function and $t(e_i)\mapsto
t(e_i)\II_{t(e_i)\geq a}$ is of course non-decreasing. Furthermore,
since the variables $(t(e_i))$ are independent, the conditional
expectation $\EE\left(.|(t(e_j))_{j\not= i}\right)$ corresponds to
  expectation over $t(e_i)$, keeping $(t(e_j))_{j\not= i}$
  fixed. Thus, Chebyshev's association inequality (see \cite{HardyLittlewoodPolya52}, p.~43) implies:
\begin{eqnarray*}
\EE\left(t(e_i)\II_{t(e_i)\geq a}\II_{e_i\in
  R_{\tilde\tau_n^r}}|(t(e_j))_{j\not= i}\right)&\leq &\EE\left(t(e_i)\II_{t(e_i)\geq a}|(t(e_j))_{j\not= i}\right)\EE\left(\II_{e_i\in
  R_{\tilde\tau_n^r}}|(t(e_j))_{j\not= i}\right)\;,\\
&=&\EE\left(t(e_1)\II_{t(e_1)\geq a}\right)\EE\left(\II_{e_i\in
  R_{\tilde\tau_n^r}}|(t(e_j))_{j\not= i}\right)\;.
\end{eqnarray*}
Thus,
\begin{equation}
\label{eq:nucontinue}0\leq \EE(\tau_n^r)-\EE(\tilde\tau_n^r) \leq
\EE\left(t(e_1)\II_{t(e_1)\geq a}\right)\EE( \card (R_{\tilde\tau_n^r}))\leq
  r\EE\left(t(e_1)\II_{t(e_1)\geq a}\right)\;.
\end{equation}
Now, since $F$ has a finite moment of order 1, we can choose
$a(\eta,F,d)$ such that:
$$\frac{1-\eta}{\eps}\EE\left(t(e_1)\II_{t(e_1)\geq
    a}\right)\leq\frac{\eta}{2}\;,$$
to get:
\begin{eqnarray}
\nonumber
 \EE(\tau_n^r)-\EE(\tilde\tau_n^r)\leq\frac{\eta}{2}\EE(\tau_n)\leq\frac{\eta}{2}\EE(\tau_n^r)\;,\\
\label{eq:taunrtautilde}
\PP(\tau_n^r\leq (1-\eta)\EE(\tau_n^r))\leq\PP\left(\tilde \tau_n^r\leq\EE(\tilde \tau_n^r)-\frac{\eta}{2}\EE( \tau_n^r)\right)\;.
\end{eqnarray}
Now, we shall use Corollary 3 in \cite{Boucheronetal03}. To this end, we
need some notation. We take $\tilde t'$
an independent collection of capacities with the same law as $\tilde t=(\tilde t(e_i))_{i=1\ldots,m_n}$. For each
edge $e_i\in\cyl(A,h)$, we denote by $\tilde t^{(i)}$ the collection of
capacities obtained from $\tilde t$ by replacing $\tilde t(e_i)$ by
$\tilde t'(e_i)$, and leaving
all other coordinates unchanged. Define:
$$V_-:=\EE\left\lbrack\left.\sum_{i=1}^{m_n}(\tilde\tau_n^r(t)-\tilde\tau_n^r(t^{(i)}))_-^2\right|t\right\rbrack\;,$$
where $\tilde\tau_n^r(t)$ is the maximal flow through $\cyl(nA,h(n))$ when
capacities are given by $t$. Observe that:
$$\tilde\tau_n^r(t^{(i)})-\tilde\tau_n^r(t)\leq (\tilde t'(e_i)-\tilde
t(e_i))\II_{e_i\in R_{\tilde\tau_n^r}}\;,$$
and thus,
$$V_-\leq a^2\card(R_{\tilde\tau_n^r})\leq a^2r=a^2(1-\eta)\EE(\tau_n)/\eps\;.$$
Thus,
Corollary 3 in \cite{Boucheronetal03} implies that, for every
$\eta\in]0,1[$,
$$\PP \left(\tilde \tau_n^r\leq \EE(\tilde \tau_n^r)-\frac{\eta}{2}\EE( \tau_n^r)\right)\leq
e^{-\frac{\EE(\tau_n^r)^2\eta^2\eps}{16a^2(1-\eta)\EE(\tau_n)}}\leq e^{-\frac{\EE(\tau_n)\eta^2\eps}{16a^2(1-\eta)}}\;,$$
which, with inequalities (\ref{eq:taunrtautilde}) and (\ref{eq:tautaunr})
 finishes the proof of inequality (\ref{eq:dev_tau_ordre1}).


\end{dem}


\begin{rem}
If we suppose the existence of an exponential moment for $F$, then
one can get concentration inequalities: there are positive constants $D_1$ and $D_2$, depending only on $F$ and $d$ and such that, for every hyperrectangle $A$, every $h>0$ and every $u>0$,
$$\PP(|\tau(A,h)-\EE(\tau(A,h))|\geq u)\,\leq\, D_1 \exp
\left(-\frac{u^2}{D_2\mathcal{H}^{d-1}(A)}\right) +D_1 \exp
\left(-\frac{1}{D_2}\mathcal{H}^{d-1}(A)\right) \,.$$
Furthermore, for every $h\leq \exp(\mathcal{H}^{d-1}(A))$ and every $u>0$,
$$\PP(|\phi(A,h)-\EE(\phi(A,h))|\geq u)\,\leq\,
D_1 \exp \left(-\frac{u^2}{D_2\mathcal{H}^{d-1}(A)} \right)+D_1 \exp
\left(-\frac{1}{D_2}\mathcal{H}^{d-1}(A) \right)\,.$$
This can be proved much as in \cite{Zhang07}, section 9. It should be
noted that these results certainly do not give the right order of the
``typical fluctuations'', i.e., fluctuations that occur with a non negligible probability. Indeed, let $S_n$ be the square:
$$S_n\,=\,\partial\left(\left[-\frac{1}{2},n-\frac{1}{2}\right]^{d-1}\times\left\{\frac{1}{2}\right\}\right)\,.$$
We say that a set of edges $E$ ``is a cut based on $S_n$'' if it is finite,
and if every closed path in $\mathbb{Z}^d$ which is not contractible to one point in $\RR^d\smallsetminus S_n$ has to contain one edge of $E$. Let $\mathcal{E}_n$ be the set of all sets of edges which are a cut based on $S_n$ and define:
$$\tilde{\tau}_n\,=\,\inf \{V(E)|E\in\mathcal{E}_n\}\,.$$
Then, mimicking the work of \cite{Benjamini:variance}, one can prove
that the variance of $\tilde{\tau}_n$ is at most of order $C (n^{d-1}/\log
n)$ where $C$ is a constant (and there is no reason for this bound to be
optimal). It is then very
reasonable to think that $\tau(A,h)$ and $\phi(A,h)$ will inherit this
property to have ``submean'' variance, i.e. their typical fluctuations
should be small with respect to $(\mathcal{H}^{d-1}(A))^{1/2}$ when the
side lengths of $A$ tend to infinity.

Remark also that these concentration inequalities, while they reflect
the right order of lower large deviations, do not give the right asymptotic
of upper large deviations, which are of volume order. We do not know a
 simple route to reach that which would avoid the work of \cite{TheretUpper}.
\end{rem}


\subsection{Asymptotic of $\mathbb{E} (\tau)$}
\label{convergence}
\begin{figure}[ht!]
\centering

\begin{picture}(0,0)%
\includegraphics{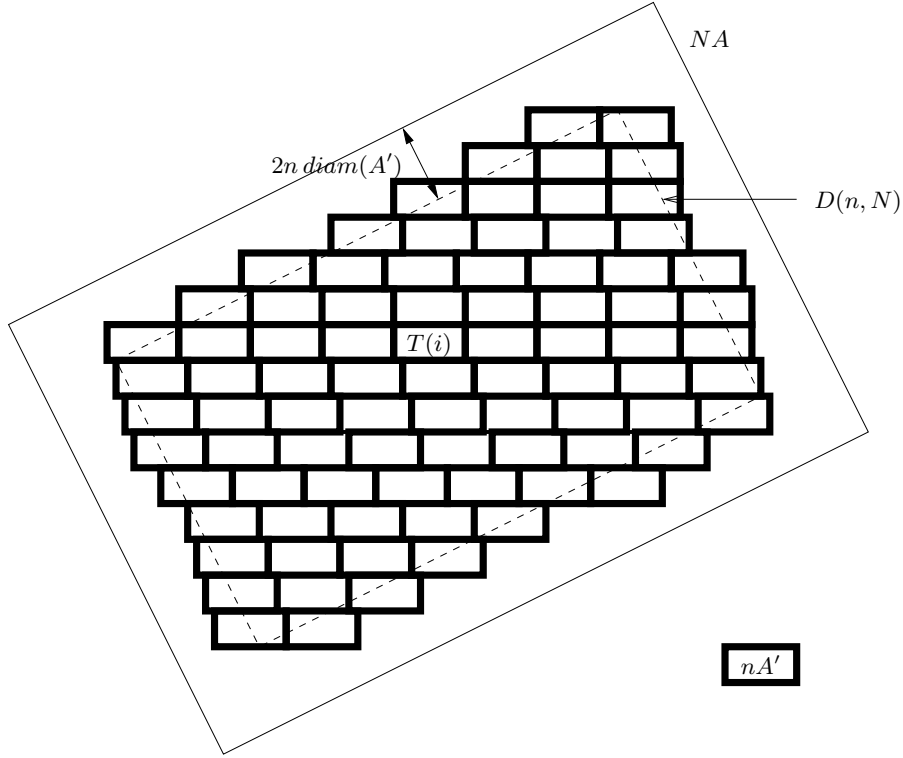}%
\end{picture}%
\setlength{\unitlength}{2960sp}%
\begingroup\makeatletter\ifx\SetFigFont\undefined%
\gdef\SetFigFont#1#2#3#4#5{%
  \reset@font\fontsize{#1}{#2pt}%
  \fontfamily{#3}\fontseries{#4}\fontshape{#5}%
  \selectfont}%
\fi\endgroup%
\begin{picture}(7515,6324)(1489,-6673)
\put(7201,-736){\makebox(0,0)[lb]{\smash{{\SetFigFont{9}{10.8}{\rmdefault}{\mddefault}{\updefault}{\color[rgb]{0,0,0}$NA$}%
}}}}
\put(5026,-3286){\makebox(0,0)[b]{\smash{{\SetFigFont{9}{10.8}{\rmdefault}{\mddefault}{\updefault}{\color[rgb]{0,0,0}$T(i)$}%
}}}}
\put(7801,-5986){\makebox(0,0)[b]{\smash{{\SetFigFont{9}{10.8}{\rmdefault}{\mddefault}{\updefault}{\color[rgb]{0,0,0}$nA'$}%
}}}}
\put(8251,-2086){\makebox(0,0)[lb]{\smash{{\SetFigFont{9}{10.8}{\rmdefault}{\mddefault}{\updefault}{\color[rgb]{0,0,0}$D(n,N)$}%
}}}}
\put(4801,-1786){\makebox(0,0)[rb]{\smash{{\SetFigFont{9}{10.8}{\rmdefault}{\mddefault}{\updefault}{\color[rgb]{0,0,0}$2n \, diam(A')$}%
}}}}
\end{picture}%

\caption{The hyperplane $\hyp(A)$.}
\label{hyperplan}
\end{figure}

Here, we prove Proposition \ref{defnu}, so we suppose that the capacity of the edges is in $L^1$.
Let us consider two
hyperrectangles $A$, $A'$ which have a common
orthogonal unit vector $\vec{v}$, and two functions $h,h':\mathbb{N}
\rightarrow \mathbb{R}^+$ such that $\lim _{n\rightarrow \infty} h(n) =\lim
_{n\rightarrow \infty} h'(n) = + \infty$. We take $n,N\in \mathbb{N}$ such that
$N\geq N_0(n)$ with $N_0(n)$ large enough to have 
$h(N)\geq h'(n) +1 $ and $N\, \diam (A) > n\, \diam (A')$ for all $N\geq N_0(n)$
(here $\diam (A) = \sup\{\|x-y\|_2 \,|\, x,y\in A  \}$). We define 
$$D(n,N)\,=\, \{x\in NA \,|\, d(x,\partial (NA)) > 2n\, \diam A'  \}\,.$$
There exists a finite
collection of sets $(T(i), i\in I)$ such that each
$T(i)$ is a translate of $nA'$ intersecting the set $D(n,N)$, the sets
$(T(i),i\in I)$ have pairwise disjoint interiors, and their union $\cup
_{i\in I} T(i)$ contains the set $D(n,N)$ (see Figure \ref{hyperplan}).
For all $i$, there exists a
vector $\vec{t_i}$ in $\mathbb{R}^d$ such that $\|\vec{t_i}\|_{\infty} <1$
and $T'(i)=T(i) +\vec{t_i}$ is the image of $nA'$ by an integer translation
(that leaves $\mathbb{Z}^d$ globally invariant). The cylinders $\cyl
(T'(i), h'(n))$ are still included in $\cyl(NA, h(N))$ for all $i\in I$,
and the family $(\tau(T'(i), h'(n)), i\in I)$ is identically distributed
(but not independent in general). For each $i$, by the max-flow min-cut theorem, we know that $\tau (T'(i),h'(n))$
is equal to the minimal capacity $\smash{V(E) = \sum_{e\in E} t(e)}$ of a set of edges
$E \subset \cyl(T'(i), h'(n))$ that cuts $\smash{T'(i)_1^{h'(n)}}$ from $\smash{T'(i)_2^{h'(n)}}$. For each $i\in I$, let $E_i$ be such a set
of edges of minimal capacity, i.e., $\smash{\tau
(T'(i), h'(n))= V(E_i)}$.

We fix $\zeta=4d $.
Let $E_0^{1}$ (resp. $E_0^{2}$, $E_0$) be the set of the edges included in
$\mathcal{E}_0^1$ (resp. $\mathcal{E}_0^2$, $\mathcal{E}_0$), where we define
$$ \mathcal{E}_0^1 \,=\, \bigcup_{i \in I} \left( \mathcal{V}( \cyl(\partial T'(i), +\infty) ,\zeta)
\cap \mathcal{V} (\hyp (NA),\zeta) \right)\,,  $$
$$  \mathcal{E}_0^2 \,=\,  \cyl(NA\smallsetminus D(n,N), \zeta) $$
and
$$ \mathcal{E}_0 \,=\, \mathcal{E}_0^1 \cup \mathcal{E}_0^2 \,.$$
The set of edges $E_0 \cup \bigcup_{i\in I}
E_i$ cuts $(NA)_1^{h(N)}$ from $(NA)_2^{h(N)}$ in $\cyl (NA,h(N))$, so 
\begin{eqnarray}
\label{sousadd}
\tau (NA, h(N))& \leq& V(E_0) + \sum_{i\in I} V(E_i) \nonumber \\
& \leq & V(E_0) + \sum_{i\in I} \tau (T'(i),h'(n)) \,.
\end{eqnarray}
Taking the expectation of (\ref{sousadd}), we obtain
\begin{align}
\nonumber\frac{\mathbb{E} \left( \tau(NA,h(N)) \right)}{\mathcal{H}^{d-1}(NA)}
& \,\leq\, \frac{\card(E_0)}{\mathcal{H}^{d-1}(NA)} \mathbb{E}(t) +
\frac{\card(I) \mathbb{E}(\tau(nA', h'(n)))}{\mathcal{H}^{d-1}(NA)} \\
\label{eq:sousaddmean}&  \,\leq\, \frac{\card(E_0)}{\mathcal{H}^{d-1}(NA)} \mathbb{E}(t) +
\frac{ \mathbb{E}(\tau(nA', h'(n)))}{\mathcal{H}^{d-1}(nA')} \,.
\end{align}
There exists a constant $c(d)$ such that:
\begin{equation}
\label{eq:borneE01et2}
 \card (E_0^1) \,\leq\, c(d)
\frac{\mathcal{H}^{d-1}(NA)}{\mathcal{H}^{d-1}(nA')}\mathcal{H}^{d-2}(\partial(nA'))
\quad \textrm{and} \quad  \card (E_0^2) \,\leq\, c(d)
\mathcal{H}^{d-2}(\partial(NA))diam(nA') \,,
\end{equation}
thus
\begin{equation}
\label{eq:borneE0}
\card (E_0) \,\leq\,
c(d)\left[\frac{\mathcal{H}^{d-1}(NA)}{\mathcal{H}^{d-1}(nA')}\mathcal{H}^{d-2}(\partial(nA'))+\mathcal{H}^{d-2}(\partial(NA))diam(nA')\right]\,,
\end{equation}
and so
$$ \lim_{n\rightarrow \infty} \lim_{N\rightarrow \infty}
\frac{\card(E_0)}{\mathcal{H}^{d-1}(NA)} \,=\, 0\,.$$
By sending $N$ to infinity, and then $n$ to infinity, we obtain that
$$ \limsup_{N\rightarrow \infty}\frac{\mathbb{E} \left( \tau(NA,h(N))
  \right)}{\mathcal{H}^{d-1}(NA)} \,\leq\, \liminf_{n\rightarrow \infty}
  \frac{\mathbb{E} \left( \tau(nA',h'(n))
  \right)}{\mathcal{H}^{d-1}(nA')}\,.  $$
For $A=A'$ and $h=h'$, we deduce from this inequality that
  $\lim_{n\rightarrow \infty}\mathbb{E} \left( \tau(nA,h(n))
  \right) / \mathcal{H}^{d-1}(nA) $ exists. For different $A,A',$ and
  $h,h'$, we conclude that this limit does not depend on $A$ and $h$, but
  only on the direction of $\vec{v}$ (and on $F$ and $d$ of course). We
  denote this limit by $\nu(\vec{v})$.


\subsection{Properties of $\nu$}
\label{subsec:thmnu}
Here we prove Proposition \ref{thmnu}. Lemma \ref{plaquettes} implies that $\nu(\vec{v})\geq
\delta \|\vec{v}\|_1$ for every unit vector $\vec{v}$, so we only need
to prove assertions $(ii)$ and $(iii)$ in Proposition \ref{thmnu}. First, let us show that $\nu(\vec{v})>0$ is equivalent to
$F(0)<1-p_c(d)$. We begin by stating the weak triangle inequality for $\nu(\vec{v})$:
\begin{prop}
\label{trig}
We suppose that $\mathbf{(F2)}$ holds. Let
$(ABC)$ be a non-degenerate triangle in $\mathbb{R}^d$ and let $\vec{v_A}$,
$\vec{v_B}$ and $\vec{v_C}$ be the exterior normal unit vectors to the sides
$[BC]$, $[AC]$, $[AB]$ in the plane spanned by $A$, $B$, $C$. Then
$$ \mathcal{H}^1 ([AB]) \nu(\vec{v_C}) \,\leq\, \mathcal{H}^1 ([AC])
\nu(\vec{v_B}) + \mathcal{H}^1 ([BC]) \nu(\vec{v_A}) \,. $$
\end{prop}

We do not prove Proposition \ref{trig} as it is the strict analogue of
Proposition 11.2 in \cite{Cerf:StFlour}. We stress the fact that it
uses only the definition of $\nu(\vec{v})$ as the limit of the expectation of
the rescaled variable $\tau$, i.e. Proposition \ref{defnu}. As in \cite{Kesten:StFlour}, one can extend $\nu$ as a function on
$\RR^d$ as follows:
$$\nu(\vec{0})=0,\quad\mbox{and}\quad\forall
\vec{u}\not=\vec{0},\;\nu(\vec{u}):=\|\vec{u}\|.\nu\left(\frac{\vec{u}}{\|\vec{u}\|}\right)\;.$$
Then, Proposition \ref{trig} shows that $\nu$ is convex (and even subadditive). Using this convexity, it is standard to obtain that
$$ \exists \vec{v} \neq \vec{0} \textrm{ s.t. } \nu(\vec{v})=0
\,\iff \, \forall \vec{v} \quad \nu(\vec{v}) \,=\,0 \,,  $$
see for example (3.15) in \cite{Kesten:StFlour}. We deduce that
\begin{equation}
\label{nu=0}
F(0) \geq 1-p_c(d)  \,\iff\, \exists \vec{v} \neq 0 \textrm{ s.t. }
\nu(\vec{v})=0 \,\iff\, \forall \vec{v} \quad \nu(\vec{v})=0\,.
\end{equation}

Now we study the case $\delta >0$. For a given realization of $(t(e), e\in
\mathbb{E}^d)$, we define the family of variables $(t'(e), e\in
\mathbb{E}^d)$ by $t'(e) = t(e) - \delta$ for all $e$. Then the variables
$(t'(e), e\in \mathbb{E}^d)$ are independent and identically distributed,
and if we denote by $F'$ their distribution function, we have $F'(\lambda)
= F(\lambda + \delta)$ for all $\lambda \in \mathbb{R}$. We compare the
variable $\tau(nA, h(n))$ and the corresponding variable $\tau'(nA, h(n))$
for the capacities $(t'(e))$, for a given hyperrectangle $A$ of normal unit
vector $\vec{v}$, and a given height function $h$ such that
$\lim_{n\rightarrow \infty} h(n) = +\infty$. We still denote by $\mathcal{N}(nA, h(n))$
the minimal number of edges that can disconnect $(nA)_1^{h(n)}$ from
$(nA)_2^{h(n)}$ in $\cyl(nA, h(n))$. By the max-flow min-cut theorem, we
easily obtain that
$$ \tau(nA, h(n)) \,\geq \, \tau'(nA, h(n)) + \delta \mathcal{N}(nA,
h(n))\,,$$
and so
$$ \frac{\mathbb{E}(\tau(nA, h(n)))}{\mathcal{H}^{d-1}(nA)} \,\geq \, \frac{\mathbb{E}(\tau'(nA, h(n)))}{\mathcal{H}^{d-1}(nA)} + \delta \frac{\mathcal{N}(nA,
h(n))}{\mathcal{H}^{d-1}(nA)}\,.  $$
Proposition \ref{defnu} and Lemma \ref{plaquettes} give us that
$$ \nu_F(\vec{v}) \,\geq \, \nu_{F'} (\vec{v}) + \delta \|\vec{v}\|_1
 $$
with trivial notations. Now $ F(\delta) =F'(0) < 1-p_c(d)$ implies that
$\nu_{F'} (\vec{v}) >0$, so $(ii)$  is proved. 

Finally, from inequalities (\ref{eq:sousaddmean}) and
(\ref{eq:borneE0}), with $A=A'$ and letting $N$ go to
infinity, we get, for every non-degenerate hyperrectangle $A$
orthogonal to some unit vector $\vec{v}$:
$$\nu(\vec{v})\leq \inf_{n\in\NN}\left\{\frac{\EE(t(e)) c(d) \mathcal{H}^{d-2}(\partial A)}{n\mathcal{H}^{d-1}(A)}+\frac{\EE(\tau(nA,h(n)))}{\mathcal{H}^{d-1}(nA)}\right\}\;.$$
Thus, Proposition \ref{thmnu} is proved.


\subsection{Law of large numbers for $\tau$}
\label{subsec:llntau}

Here, we prove Theorem \ref{thm:llntau}. We begin with the almost sure
convergence of $\tau(nA,h(n))/\mathcal{H}^{d-1}(nA)$. To deduce it from the convergence of
its expectation, we will use the following result:
\begin{lem}
\label{lem:deviation}
Suppose that hypotheses $\mathbf{(F1)}$ and $\mathbf{(F2)}$ occur. Then
$$
\liminf_{n\rightarrow \infty}\frac{\tau(nA,h(n))-\EE(\tau(nA,h(n)))}{\mathcal{H}^{d-1}(nA)}\geq 0\quad\mbox{ a.s. }$$
\end{lem}
\begin{dem}
It is a simple consequence of Proposition \ref{prop:deviation} and the fact that $\EE(\tau(nA,h(n)))$ is equivalent
to $\mathcal{H}^{d-1}(nA)\nu(\vec{v})$, using
Borel-Cantelli's lemma.
\end{dem}
We shall use (\ref{sousadd}) with  $h=h'$ and
$A=A'$, i.e. the sets $T'(i)$ are integer translates of $nA$. We
emphasize the dependence on $N$ and $n$ by writing $E_0^i=E_0^i(N,n)$ for $i\in\{1,2\}$, $I=I(N,n)$ and
$T'(i)=T'_{N,n}(i)$. Suppose first that $0\in A$. Then, we can
construct the sets $T'_{N,n}(i)$ in order to have:
$$\forall n\geq 1,\forall N'\geq N\geq N_0(n),\;(T'_{N,n}(i))_{i\in
  I(N,n)}\subset(T'_{N',n}(i))_{i\in I(N',n)} \;.$$
We obtain that
$$ \forall n\geq 1,\forall N'\geq N\geq N_0(n),\; E_0^1(N,n) \,\subset \,
E_0^1 (N',n) \,.$$
Thus, the strong law of large numbers for i.i.d. random variables
implies, using inequality (\ref{eq:borneE01et2}):
\begin{equation}
\label{eq:llnE01}
\limsup_{N\rightarrow \infty}\frac{V(E_0^1)}{\mathcal{H}^{d-1}(NA)}\leq
\EE(t(e))\limsup_{n \rightarrow \infty}\frac{\card
  (E_0^1)}{\mathcal{H}^{d-1}(NA)}\leq \frac{\EE(t(e))K(d,A)}{n}  \quad
\mbox{ a.s.}
\end{equation}
where $K(d,A)=c(d) \mathcal{H}^{d-2}(\partial A)/\mathcal{H}^{d-1}(A)$.
Moreover, we know (see (\ref{eq:borneE01et2})) that
$$ \card{E_0^2} \,\leq\, c(d) \mathcal{H}^{d-2} (\partial A) \diam (A)
N^{d-2}n \,. $$
Under the assumption $\mathbf{(F2)}$, Theorem 4.1 in \cite{Gut92} states that
$V(E_0^2(N,n))/\mathcal{H}^{d-1}(NA)$ converges completely to $0$, with the
definition of the complete convergence given by Gut (Definition (1.1) in
\cite{Gut92}). Complete convergence implies almost sure convergence
through Borel-Cantelli's lemma, thus
$$ \lim_{N\rightarrow \infty} \frac{V(E_0^2(N,n))}{\mathcal{H}^{d-1}(NA)}
\,=\, 0 \quad \mbox{a.s.} $$
Also, we claim that:
\begin{equation}
\label{eq:llntaudependent}
\limsup_{N \rightarrow \infty}\frac{\sum_{i\in I(N,n)} \tau
  (T'(i),h'(n))}{\mathcal{H}^{d-1}(NA)}\,=\,\frac{\EE(\tau(nA,h(n)))}{\mathcal{H}^{d-1}(nA)} \quad \mbox{ a.s.}
\end{equation}
Indeed, notice that (for $n$ large enough) $\tau (T'(i),h'(n))$ is independent of all the
other $\tau (T'(j),h'(n))$ except for at most $3^d-1$ values of $j$
corresponding to the $T'(j)$ that can intersect $T'(i)$. Thus,
(\ref{eq:llntaudependent}) follows by partitioning the sets $T'(j)$ into
$3^d-1$ classes of i.i.d. variables, and then applying  the strong law
of large numbers for i.i.d. random variables. Thus, for $n$ large enough,
$$\limsup_{N \rightarrow \infty}\frac{\tau (NA,
  h(N))}{\mathcal{H}^{d-1}(NA)}\,\leq \,\frac{\EE(t(e))K(d,A)}{n} +
\frac{\EE(\tau(nA,h(n)))}{\mathcal{H}^{d-1}(nA)} \quad \mbox{ a.s.}$$
and using Proposition \ref{defnu}:
\begin{equation}
\label{eq:upperboundlln}
\limsup_{N \rightarrow \infty}\frac{\tau (NA, h(N))}{\mathcal{H}^{d-1}(NA)}\,\leq\,
\nu(\vec{v})\quad \mbox{ a.s.}
\end{equation}
If $\nu(\vec{v})=0$, since $\tau$ is non-negative, we get the desired
result. We suppose that $\nu(\vec{v})>0$. From Proposition
\ref{thmnu}, we know that $\nu(\vec{v})>0$ is equivalent to $F(0)<1-p_c(d)$.
Then it follows from Lemma \ref{lem:deviation} and
the convergence of $\EE(\tau(nA,h(n)))/\mathcal{H}^{d-1}(nA)$ to
$\nu(\vec{v})$ that:
$$\nu(\vec{v})\leq\liminf_{N \rightarrow \infty}\frac{\tau (NA,
  h(N))}{\mathcal{H}^{d-1}(NA)}\quad \mbox{ a.s.}$$
which, together with (\ref{eq:upperboundlln}) gives the law of large
numbers for $\tau$.

Now, what happens if $0\not\in A$ ? Then, we suppose that
$\mathbf{(F3)}$ holds, and we can combine Borel-Cantelli's
Lemma with the complete convergence in the
law of large numbers for subsequences (Theorem 4.1 in \cite{Gut85}, or more
generally Theorem 4.1 in \cite{Gut92})  to
replace the classical law of large numbers to prove (\ref{eq:llnE01}) and
(\ref{eq:llntaudependent}).


This ends the proof of the almost sure convergence. Now, let us prove
the convergence in $L^1$. Suppose first that $0\in A$. Then, one can
find a sequence of sets of edges $(E(n))_{n\in\NN}$ such that for each
$n$, $E(n)$ is a cut between $(nA)_1^{h(n)}$ and $(nA)_2^{h(n)}$,
$E(n)\subset E(n+1)$ and:
$$\frac{\card (E(n)) }{\mathcal{H}^{d-1}(nA)}\xrightarrow[n\rightarrow
\infty]{}\|\vec{v}\|_1\;.$$
Now, define:
$$f_n=\frac{\tau_n}{\mathcal{H}^{d-1}(nA)}\mbox{ and }\quad g_n=\sum_{e\in
  E(n)}t(e)\;.$$
Then, we know the following:
\begin{itemize}
\item[(i)] $0\leq f_n\leq g_n$ for every $n$,
\item[(ii)] $(g_n)_{n\in\NN}$ converges almost surely and in $L^1$, thanks to the usual law of large numbers,
\item[(iii)] $(f_n)_{n\in\NN}$ converges almost surely to
  $\nu(\vec{v})$, thanks to the almost sure convergence for $0\in A$
  that we have just proven,
\item[(iv)] $(\EE(f_n))_{n\in\NN}$ converges  to   $\nu(\vec{v})$,
  thanks to Proposition \ref{defnu}.
\end{itemize}
It is then standard to show that $f_n$ converges in $L^1$ to
$\nu(\vec{v})$: apply the monotone convergence theorem to
$b_n=\inf_{m\geq n}(g_m-f_m)$, and then show that $(g-f-b_n)_{n\in\NN}$ and
$(g_n-f_n-b_n)_{n\in\NN}$ are positive sequences converging to zero in
$L^1$.

It remains to show the convergence in $L^1$ when we do not know
whether $0\in A$. Let $A''$ be the translate of $A$ such that $0\in A''$, and
$0$ is the center of $A''$. For
any fixed $n$, there exists a hyperrectangle $A_n'$ which is a translate
of $nA$ by an integer vector and such that $d_\infty(0,nA'_n)<1$ and
$d_\infty(nA'',A'_n)<1$, where $d_\infty$ denotes the distance induced by
$\|.\|_\infty$. We want to compare the maximal flow through
$\cyl(nA'',h(n))$ to the maximal flow through $\cyl (A_n',h(n))$. The
difficulty is that one of these cylinders is not included in the other. This is the reason why we will
construct bigger and smaller version of $\cyl(nA'', h(n))$. We recall that
$l_{\min}(A)$ is the smallest length of $A$, i.e.,
$$ l_{\min} (A) \,=\, \min_{i=1,...,d-1} k_i \,, $$
where $A$ is the image by an isometry of the set
$\prod_{i=1}^{d-1}[0,k_i]\times \{0\}$. We define the biggest length of $A$
as
$$ l_{\max} (A) \,=\, \max_{i=1,...,d-1} k_i $$
with the same notation. We only consider $n$ large enough such that $h(n) >1$. Thus the following inclusions holds:
$$ \cyl \left( \left( n- \left\lceil \frac{2}{l_{\min}(A)}  \right\rceil
  \right) A'', h(n) -1 \right)\,\subset\, \cyl(A_n',h(n)) \,\subset\, \cyl \left( \left( n+ \left\lceil \frac{2}{l_{\min}(A)}  \right\rceil
  \right) A'', h(n) +1\right) \,, $$
where $\lceil x \rceil$ is the smallest integer bigger than or equal to
$x$. For all $n$, we have
$$ \partial \left[ \left( n- \left\lceil \frac{2}{l_{\min}(A)}  \right\rceil
  \right) A'' \right] \,\subset\,  \mathcal{V} \left(\partial A_n' , l_{\max}(A)
  \left\lceil \frac{2}{l_{\min}(A)} \right\rceil + 1   \right)  $$
and
$$ \partial \left[ \left( n+ \left\lceil \frac{2}{l_{\min}(A)}  \right\rceil
  \right) A'' \right] \,\subset\,  \mathcal{V} \left(\partial A_n' , l_{\max}(A)
  \left\lceil \frac{2}{l_{\min}(A)} \right\rceil + 1   \right) \,. $$
Argueing as in section \ref{convergence}, let $F_n$ be the edges
included in $\mathcal{F}_n$ defined as
$$ \mathcal{F}_n \,=\, \mathcal{V} \left(\partial A_n' , l_{\max}(A)
  \left\lceil \frac{2}{l_{\min}(A)} \right\rceil + 1 + 4d  \right) \,.$$
We get, for $n$ large enough,
\begin{align*}
\tau \left( \left( n+ \left\lceil \frac{2}{l_{\min}(A)}  \right\rceil
  \right) A'', h(n) +1\right)& - V(F_n)\\ &\,\leq\,  \tau(A_n',h(n))
\,\leq\, \tau \left(  \left( n- \left\lceil
      \frac{2}{l_{\min}(A)}  \right\rceil \right) A'', h(n) -1 \right)
+V(F_n)\,.
\end{align*}
Using the convergence in $L^1$ for $A''$ which contains 0, we see that
$$\tau \left( \left( n+ \left\lceil \frac{2}{l_{\min}(A)}  \right\rceil
  \right) A'', h(n) +1\right) /\mathcal{H}^{d-1}(nA) \quad\textrm{and} \quad
\tau  \left(  \left( n- \left\lceil
      \frac{2}{l_{\min}(A)}  \right\rceil \right) A'', h(n) -1 \right)  /\mathcal{H}^{d-1}(nA)$$ 
converge to
$\nu(\vec{v})$ in $L^1$ as $n$ goes to infinity. Furthermore, since
$\card (F_n)$ is negligible compared to $n^{d-1}$, $V(F_n)/\mathcal{H}^{d-1}(nA)$ go to zero in $L^1$, and we get the convergence of
$\tau(A'_n,h(n))/\mathcal{H}^{d-1}(nA)$ to $\nu(\vec{v})$ in $L^1$. But since
$A_n'$ is an integer translate of $nA$, it implies the convergence of
$\tau(nA,h(n))/\mathcal{H}^{d-1}(nA)$ to $\nu(\vec{v})$ in $L^1$.

\begin{rem}
Most likely, the almost sure convergence of $(\tau(nA, h(n))/\mathcal{H}^{d-1}(nA), n
\in \mathbb{N})$ could also be obtained by adapting the proof of
\cite{Ackoglu}, and thus relaxing the independence hypothesis
on $(t(e))_e$ to stationarity. In any case, general subadditive
results existing in the literature are not well adapted to treat the
case of irrational directions, i.e. directions
$\vec{v}$ such that $\tau(nA,h(n))$ is not exactly subadditive and
stationary. Some authors circumvent this problem by proving that the
almost sure convergence is uniform with respect to rational
directions, which allows to extend the convergence to irrational
directions, see \cite{Kesten:StFlour} and \cite{Boivin} for
instance. But for flows like $\tau$, the uniform convergence requires a moment of order
strictly larger than 1, see for instance Theorems 1.3, 1.9 and section
4 in \cite{Boivin}. Notice also that Theorem 6.1 in \cite{Boivin}
shows directly the convergence in any direction for First Passage
Percolation in dimension 2, using techniques some of which are similar
to ours and others belong to the realm of ergodic theory. In this
paper, the
  strategy we adopt is to use the fact that our space $\RR^d$ has one
  dimension more than the hyperrectangles which are the
  indices of the almost subadditive family: we can move the
  hyperrectangles $T(i)$ out of the hypersurface spanned by $NA$ to obtain
  the hyperrectangles $T'(i)$ that have good properties. Moreover, the
  non-negativity of our variables $\tau$ implies that it is simpler to use
a concentration inequality than a maximal inequality as in the classical
subadditive ergodic theorems. 
\end{rem}
\begin{rem} We have obtained readily the independence of the
limit with regard to the precise form of the hyperrectangle we
consider and this is not surprising since it appears already in subadditive
ergodic theorems like in \cite{KrengelPyke}. 
\end{rem}

\begin{rem}
The almost sure convergence of $(\tau(nA, h(n))/\mathcal{H}^{d-1}(nA), n
\in \mathbb{N})$ is not necessary to prove Theorem \ref{devinf}, but we
need the convergence in probability to prove Theorem \ref{thmpgd}.
\end{rem}

\begin{rem}
\label{rem:condF3}
If $0\not\in A$, it is not clear to us whether condition $\mathbf{(F3)}$
is necessary or not: it is necessary for complete
convergence to hold, but complete convergence is stronger than the
a.s. convergence.
\end{rem}


\subsection{Lower deviations for $\tau$ and flat $\phi$: proofs of Theorem \ref{devinf} and Theorem \ref{cordevinf}}
\label{sectiondevinf}

Now, we can prove Theorem \ref{devinf}, and so we consider $F$, $h$,
$\vec{v}$ and $A$
as in the statement of this Theorem.
If $\nu(\vec{v})=0$, there is nothing to prove. Suppose now that
$\nu(\vec{v})>0$ and let $\eps\leq \nu(\vec{v})$ be a positive real
number. Let $u=\eps/(2\nu_{max})$, where
$\nu_{max} = \max \{ \nu(\vec{v})\,|\, \vec{v} \,\, unit\,\, vector
\}$. Then $u>0$ and we have
$$\frac{\nu(\vec{v})-\eps}{\nu(\vec{v})-\varepsilon /2}\,\leq\, 1-u\,.$$
Using assertion $(iii)$ in Proposition \ref{thmnu}, we know that there
exists a $n_0=n_0(A)$ (not depending on $h$) large enough to have
$$\forall n\geq n_0 \qquad
\frac{\EE(\tau(nA,h(n)))}{\mathcal{H}^{d-1}(nA)}\,\geq\,
\nu(\vec{v})-\frac{\varepsilon}{2}\,.$$
Then, for all $n\geq n_0$,
$$\mathbb{P} \left[ \tau (nA,h(n)) \leq \left(\nu(\vec{v})-\eps
  \right) \mathcal{H}^{d-1}(nA) \right] \leq  \mathbb{P} \left[\frac{\tau (nA,h(n))}{\EE(\tau(nA,h(n)))} \leq
  1-u \right]\,.$$
Now, the result follows easily from Proposition \ref{prop:deviation},
for $n$ larger than some $n_1=n_1(A)$. Adapting the constant for $n\leq n_1$
leads to $\widetilde{C}(d,F,A,\varepsilon)$.
\begin{rem}
Notice that for every hyperrectangle $A$:
$$\frac{2}{l_{min}(A)}\leq \frac{\mathcal{H}^{d-2}(\partial A)}{\mathcal{H}^{d-1}(A)}\leq \frac{2(d-1)}{l_{min}(A)}\;.$$
Thus, from the proof above, Proposition \ref{thmnu} $(iii)$ and Proposition
\ref{prop:deviation}, it can be seen that $n_1(A)$ and thus the constant
$\widetilde{C}(d,F,A,\varepsilon)$ depends on $A$ only through
$K(d,A)$, or equivalently, only through $l_{min}(A)$.
\end{rem}


We can do the same calculus for $\phi(nA, h(n))$ as soon as we know
that $\EE(\phi(nA, h(n)))/\H^{d-1}(nA)$ converges to
$\nu(\vec{v})$. To
prove Theorem \ref{cordevinf}, it is
sufficient to prove that it is the case under hypotheses $\mathbf{(F2)}$,
$\mathbf{(H1)}$ and $\mathbf{(H3)}$. We have to compare $\phi$ and
$\tau$. We suppose that $\lim_{n\rightarrow \infty} h(n) /n =0$, and
fix $\zeta\geq 2d$. We consider $n$ large enough such that the sides of
$nA$ have length bigger than $\zeta$, i.e., $l_{min}(A)\geq\zeta$.
Let $E_1^+ $ be the set of the edges that belong to $\mathcal{E}_1^+$,
defined as
$$ \mathcal{E}_1^+ \,=\, \mathcal{V}(\cyl(\partial (nA), h(n)), \zeta) \cap \cyl(nA, h(n)) \,. $$
We have, for all $n $ large enough,
$$\tau(nA, h(n)) \,\geq\, \phi(nA, h(n)) \,\geq\, \tau(nA, h(n)) - V(E_1^+) \,. $$
There exists a constant $C^+$ such that
$$\card(E_1^+) \leq C^+ n^{d-2} h(n) \,,$$
so we have
$$ \frac{|\EE[\phi(nA, h(n))] - \EE[\tau(nA, h(n))]|}{\H^{d-1}(nA)}
\,\leq\, \frac{C^+ n^{d-2} h(n)}{n^{d-1} \H^{d-1}(nA)}
\,\longrightarrow \,0 \quad \textrm{as }n\rightarrow \infty\,, $$
and this proves the convergence of $\EE[\phi(nA, h(n))]/\H^{d-1}(nA)$ to
$\nu(\vec{v})$. Notice that the speed of convergence depends on $h$.
Using Proposition \ref{prop:deviation}, we can find $n_1(d,F,A,h,\eps)$ such that
for all $n\geq n_1(d,F,A,h,\eps)$ we have
\begin{align*}
\mathbb{P} (\phi(nA,h(n)) \leq  & (\nu(\vec{v}) -
\eps) \mathcal{H}^{d-1}(nA) )\\
& \,\leq\, C_3(F(0),d) h(n) \exp \left(-C(\eps,
F, d) (\nu(\vec{v})-\eps/2) \mathcal{H}^{d-1}(nA)\right)\\
& \,\leq\,  C_3(F(0),d) \exp \left( \frac{\log h(n)}{n^{d-1}} n^{d-1} - C(\eps,
F, d) (\nu(\vec{v})-\eps/2) \mathcal{H}^{d-1}(nA)\right)\,.
\end{align*}
Using hypothesis $\mathbf{(H2)}$, which is implied by $\mathbf{(H3)}$,
Theorem \ref{cordevinf} is proved for $n\geq n_2(d,F, A, h, \eps)$, for
$n_2(d,F,A,h,\eps)$ large enough. Adapting the constant $C_3(F(0),d)$ for
the $n_2$ first terms, Theorem \ref{cordevinf} is proved for all $n$ with a
constant $\widetilde{C}'$ depending on $d$, $F$, $A$, $h$, $\eps$.

To prove Theorem \ref{devinfphi}, it remains to prove the convergence of
$\EE[\phi(nA, h(n))]/\H^{d-1}(nA)$ to $\nu(\vec{v})$ under the hypotheses
$\mathbf{(F1)}$, $\mathbf{(F2)}$, $\mathbf{(H1)}$ and $\mathbf{(H2)}$. This
will be done during the proof of Theorem \ref{thm:llnphistraight} in
section \ref{sec:llnphistraight}, so we postpone the end of the proof
of Theorem \ref{devinfphi} until section \ref{secfinal}.

\begin{rem}
Using Theorem \ref{cordevinf}, Theorem \ref{thm:llntau} and the fact that
$\phi(nA,h(n))\leq \tau(nA,h(n))$, we obtain the law of large numbers for
$\phi(nA,h(n))$ in flat cylinders (i.e., under hypothesis $\mathbf{(H3)}$)
under the same hypothesis as the one for $\tau(nA,h(n))$.
\end{rem}

\section[LDP for $\tau$]{Large deviation principle for $\tau$ and $\phi$ in
  flat cylinders}
\label{sec:pgdtau}
In this section, we show the large deviation principle for $\tau$. We
construct a precursor of the rate function in section \ref{existence}, and
then study its properties. Precisely, we show it is convex in section
\ref{convexity}, finite (and thus continuous) on
$]\delta\|\vec{v}\|_1,+\infty[$ in section \ref{cont}, and strictly positive on
$[0,\nu(\vec{v})[$ in section \ref{positivity}. After having shown in
section \ref{secdevsup} that upper large deviations occur at an order
bigger than the surface order, we can complete the proof of the full
large deviation principle for $\tau$ in section \ref{pgd} and deduce
the one for $\phi$ in flat cylinders in section \ref{subsec:ldpflatphi}.

\subsection{Construction of the rate function}
\label{existence}

We will prove the following lemma, for which no condition on $F$ is required.
\begin{lem}
\label{limitetaupert}
For every function $h:\mathbb{N} \rightarrow \mathbb{R}^+$ satisfying
$\mathbf{(H1)}$, for every non-degenerate hyperrectangle $A$, for
all $\lambda$ in $\mathbb{R}^+$, the limit
$$ \lim_{n\rightarrow \infty} \frac{-1}{\mathcal{H}^{d-1}(nA)} \log
\mathbb{P} \left[ \tau (nA,h(n)) \leq \left(\lambda-\frac{1}{\sqrt{n}}
  \right) \mathcal{H}^{d-1}(nA) \right] $$
exists in $[0,+\infty]$ and depends only on the direction of $\vec{v}$, one of the two unit
vectors orthogonal to $\hyp (A)$. We denote it by
$\mathcal{I}_{\vec{v}}(\lambda)$.
\end{lem}

We introduce a factor $1/\sqrt{n}$ in the definition of
$\mathcal{I}_{\vec{v}}(\lambda)$ because we want to work with subadditive
  objects, but $\tau(A,h)$ is not subadditive in $A$, except for straight cylinders. Indeed, if $A$ and
  $B$ are two hyperrectangles with a common orthogonal vector and with a
  common side, to glue
  together a set of edges in $\cyl(A,h)$ that cuts $A_1^h$ from $A_2^h$ and
  a set of edges in $\cyl(B,h)$ that cuts $B_1^h$ from $B_2^h$, we have to
  add edges at the common side of $A$ and $B$ (see the set of edges $E_0$
  defined in section \ref{convergence}). These edges may not have a capacity $0$, so they
  perturb the subadditivity of $\tau$. We add the factor $1/\sqrt{n}$ to
  compensate. 

\begin{rem}
It is natural to have no condition on $F$ in Lemma \ref{limitetaupert}
since it comes essentially from an almost subadditive property for a non-random
quantity. 
\end{rem}

\begin{dem}
For the proof of Lemma \ref{limitetaupert}, we consider the same
construction as in section \ref{convergence} (see Figure~\ref{hyperplan}). From (\ref{sousadd}) we deduce that for all
$\lambda \in \mathbb{R}^+_*$, we have
\begin{align*}
\mathbb{P} \bigg[\tau (NA,h(N))  & \leq \left(\lambda-
  \frac{1}{\sqrt{N}}\right) \mathcal{H}^{d-1} (NA)\bigg  ] \\
&  \,\geq\,
\mathbb{P} \left[ V(E_0) + \sum_{i\in I} \tau (T'(i),h'(n)) \leq
  \left(\lambda - \frac{1}{\sqrt{N}}\right)
  \mathcal{H}^{d-1} (NA)\right]\,.
\end{align*}

Let $\mathcal{D}=\{\lambda \,|\, \mathbb{P} (t(e)\leq \lambda) >0  \}$, and
$\delta=\inf \mathcal{D}$. We take $u=\delta+ \zeta$, so $p = \mathbb{P} (t(e)
\leq u) >0$. We use first the FKG inequality and then the fact that the
family $(\tau (T'(i), h'(n)), i \in I)$ is identically distributed to
obtain that
\begin{align*}
\mathbb{P} \bigg[\tau(NA,h(N)) & \leq
 \left(\lambda-\frac{1}{\sqrt{N}}\right) \mathcal{H}^{d-1}(NA)\bigg]\\
 & \,\geq \, \mathbb{P} \left[V(E_0) \leq u \,\card(E_0)\right] \\
& \qquad \times  \prod_{i\in I}
 \mathbb{P} \left[\tau (T'(i),h'(n)) \leq \frac{(\lambda - 1/\sqrt{N})
 \mathcal{H}^{d-1}(NA) - u \,\card(E_0)}{\card (I)}\right]\\
& \,\geq \, \mathbb{P} \left[t(e) \leq u \right]^{\card(E_0)} \\
& \qquad \times \mathbb{P} \left[\tau (nA',h'(n)) \leq \frac{(\lambda  - 1/\sqrt{N})
 \mathcal{H}^{d-1}(NA) - u\, \card(E_0)}{\card (I)}\right]^{\card(I)}\,.
\end{align*}
We have immediately that $\card (I) \leq \mathcal{H}^{d-1}(NA) /
\mathcal{H}^{d-1}(nA')$, so
\begin{align*}
\frac{-1}{\mathcal{H}^{d-1} (NA)} & \log \mathbb{P} \left[\tau (NA,h(N)) \leq
  \left(\lambda - \frac{1}{\sqrt{N}}\right) \mathcal{H}^{d-1}(NA)\right] \\
& \,\leq\, \frac{-1}{\mathcal{H}^{d-1}(nA')} \log \mathbb{P}
  \left[\tau(nA',h'(n)) \leq \beta\right] -
  \frac{\card(E_0)}{\mathcal{H}^{d-1}(NA)}  \log p \,,
\end{align*}
where 
$$\beta \,=\,\frac{\left(\lambda-1/\sqrt{N}\right) \mathcal{H}^{d-1}(NA) - u \,\card(E_0)}{\card (I)}\,.$$ 
As we saw in section \ref{convergence}, there exists a constant $c(d,\zeta,
A, A')$ such that
$$ \card (E_0) \,\leq\, c(d,\zeta,A,A') \left(N^{d-2}n + N^{d-1}/n
+1\right)\,.$$
On one hand, we obtain that
$$\lim_{n\rightarrow \infty} \lim_{N\rightarrow \infty}
\frac{\card(E_0)}{\mathcal{H}^{d-1}(NA)}  \log p \,=\, 0 \,. $$
On the other hand we want to compare $\beta$ with $(\lambda - 1/\sqrt{n})
\mathcal{H}^{d-1} (nA')$. Obviously we have
$$ \frac{\lambda \mathcal{H}^{d-1}(NA)}{\card(I)} \,\geq \, \lambda
\mathcal{H}^{d-1}(nA') \,. $$
We also know that
$$ \card(I) \,\geq\,
\frac{\mathcal{H}^{d-1}(D(n,N))}{\mathcal{H}^{d-1}(nA')}$$
so there exist a constant $c'(d,A,A')$ and an integer  $N_1(n)$ large enough to have,
for all $N\geq N_1(n)$,
$$ \card(I) \,\geq\, c'(d,A,A') \left( \frac{N}{n} \right)^{d-1} \,.$$
Thus, there exist constants $c_i(d,\zeta, A,A')$ such that for all $N\geq
N_1(n)$, we have
$$ \frac{\mathcal{H}^{d-1}(NA)}{\card(I) \sqrt{N}} \,\leq\,
\frac{c_1(d,\zeta,A,A')}{\sqrt{N}} \mathcal{H}^{d-1}(nA') $$
and
$$ \frac{u\, \card(E_0)}{\card(I)} \,\leq\, c_2(d,\zeta,A,A') \left( \frac{n}{N}
+ \frac{1}{n} \right) \mathcal{H}^{d-1}(nA') \,.$$
There exists $n_0$ such that for all $n\geq n_0$, $c_2/n \leq
1/(4\sqrt{n})$. Then there exists $N_2(n) \geq N_0(n) \vee N_1(n)$ such
that for all $N\geq N_2(n)$, $c_2n/N \leq 1/(4\sqrt{n})$ and $c_1/\sqrt{N}
\leq 1/(2\sqrt{n})$. Thus for a fixed $n\geq n_0$, for all $N\geq N_2(n)$,
we have
$$ \beta \,\geq \, \left(\lambda - \frac{1}{\sqrt{n}}  \right)
\mathcal{H}^{d-1}(nA') \,. $$
Now in the following inequality, obtained for $n\geq n_0$ and $N\geq N_2(n)$,
\begin{align*}
\frac{-1}{\mathcal{H}^{d-1} (NA)} & \log \mathbb{P} \left[\tau (NA,h(N)) \leq
  \left(\lambda - \frac{1}{\sqrt{N}}\right) \mathcal{H}^{d-1}(NA)\right] \\
& \,\leq\, \frac{-1}{\mathcal{H}^{d-1}(nA')} \log \mathbb{P}
  \left[\tau(nA',h'(n)) \leq \left(\lambda - \frac{1}{\sqrt{n}}  \right)
  \mathcal{H}^{d-1}(nA') \right] -
  \frac{\card(E_0)}{\mathcal{H}^{d-1}(NA)}  \log p \,,
\end{align*}
we send $N$ to infinity for a fixed $n\geq n_0$, and then we send $n$
to infinity. We thus obtain
\begin{align*}
\limsup_{N\rightarrow \infty} \frac{-1}{\mathcal{H}^{d-1} (NA)} & \log
\mathbb{P} \left[\tau (NA,h(N)) \leq \left(\lambda -
  \frac{1}{\sqrt{N}}\right) \mathcal{H}^{d-1}(NA)\right] \\
& \,\leq\, \liminf_{n\rightarrow \infty} \frac{-1}{\mathcal{H}^{d-1}(nA')}
\log \mathbb{P} \left[\tau(nA',h'(n)) \leq \left(\lambda -
  \frac{1}{\sqrt{n}}  \right) \mathcal{H}^{d-1}(nA') \right] \,.
\end{align*}
For $A=A'$ and $h=h'$, this gives us the existence of
$$ \lim_{n\rightarrow \infty}  \frac{-1}{\mathcal{H}^{d-1}(nA)}
\log \mathbb{P} \left[\tau(nA,h(n)) \leq \left(\lambda -
  \frac{1}{\sqrt{n}}  \right) \mathcal{H}^{d-1}(nA) \right]  $$
for all $\lambda \in \mathbb{R}^+_*$, and for different $A,A',h,h'$ this
shows that the limit is independent of $A$ and $h$. We denote this limit by
$\mathcal{I}_{\vec{v}} (\lambda)$.

For $\lambda = 0$, 
$$\mathbb{P} \left[\tau(nA,h(n)) \leq -\frac{\mathcal{H}^{d-1}(nA)}{
  \sqrt{n}}\right] \,=\, 0$$
for all $n\in \mathbb{N}$, so the previous limit equals $+\infty$,
  independently of $A$ and $\vec{v}$. This ends the proof of Lemma
\ref{limitetaupert}.
\end{dem}

\begin{rem}
The function $\mathcal{I}_{\vec{v}}$ is not exactly the rate function we will consider later: we will change its value from $0$ to $+\infty$ on $]\nu(\vec{v}), +\infty]$ and we will regularize it at $\|\vec{v}\|_1 \delta$.
\end{rem}


\subsection{Convexity of $\mathcal{I}_{\vec{v}}$}
\label{convexity}
We will prove that $\mathcal{I}_{\vec{v}}$ is convex, i.e., for all
$\lambda_1\geq \lambda_2 \in \mathbb{R}^+$ and $\alpha \in ]0,1[$, we have
$$ \mathcal{I}_{\vec{v}}(\alpha \lambda_1 + (1-\alpha) \lambda_2) \,\leq\,
    \alpha \mathcal{I}_{\vec{v}} (\lambda_1) + (1-\alpha)
    \mathcal{I}_{\vec{v}} (\lambda_2) \,. $$
For $\lambda_2 =0$, the result is obvious, so we suppose $\lambda_2 >0$.
We keep the same notations as in the previous section, for $D(n,N)$,
    $T(i)$, $E_i$, etc..., except that we take $A=A'$. We define
$$ \gamma \,=\, \lfloor \alpha \, \card(I) \rfloor \,. $$
If we have 
\begin{equation}
\label{eq:convex1}
\tau(T'(i), h(n)) \,\leq\, (\lambda_1 -1/\sqrt{n})
 \mathcal{H}^{d-1}(nA)\qquad for\,\, i=1,...,\gamma\,,
\end{equation}
\begin{equation}
\label{eq:convex2}
\tau(T'(i), h(n)) \,\leq\, (\lambda_2 -1/\sqrt{n})
 \mathcal{H}^{d-1}(nA)\qquad for\,\, i=\gamma +1,...,\card(I) \,,
\end{equation}
and 
$$V(E_0) \,\leq\, u\,\card(E_0)\,,$$
then we obtain that
\begin{align*}
\tau(NA,h(N)) &\,\leq \, \left( \gamma(\lambda_1 -\frac{1}{\sqrt{n}}) + (\card(I)
-\gamma) (\lambda_2 - \frac{1}{\sqrt{n}}) \right) \mathcal{H}^{d-1}(nA) +
u\, \card(E_0)\,,\\
& \,\leq\,(\alpha \lambda_1 + (1-\alpha) \lambda_2) \card(I)
\mathcal{H}^{d-1}(nA) - \frac{\card(I) \mathcal{H}^{d-1}(nA)}{\sqrt{n}} +
u\, \card(E_0)\,,\\ 
& \,\leq\, (\alpha \lambda_1 + (1-\alpha)\lambda_2) \mathcal{H}^{d-1}(NA) -\rho\,,
\end{align*}
where
$$ \rho \,=\, \frac{\card(I) \mathcal{H}^{d-1}(nA)}{\sqrt{n}} -u\, \card(E_0)
\,.$$
We want to prove that $\rho \geq \mathcal{H}^{d-1}(NA)/\sqrt{N}$ for $N$
large enough. We have seen in the previous section that there exists a
constant $c(d,\zeta, A)$ such that
$$ \card(E_0) \,\leq\, c(d,\zeta,A) N^{d-1} \left( \frac{n}{N}+\frac{1}{n}  \right)\,,
$$
and that there exists a constant $c'(d,A)$ and a $N_1(n)$ large enough to have,
for all $N\geq N_1(n)$,
$$ \card(I) \geq c'(d,A) \left( \frac{N}{n} \right)^{d-1} \,.$$
There exists $n_1$ such that for all $n\geq n_1$, $2c/n \leq c'/(2\sqrt{n})$.
For a fixed $n\geq n_1$, there exists constants $c_i(d,\zeta,A)$ and $N_3(n)$
such that for all $N\geq N_3(n)$ we have
$$ \frac{u \,\card(E_0)}{\mathcal{H}^{d-1}(NA)} \,\leq \, \frac{2c}{n}
\,\leq\,\frac{c'}{2\sqrt{n}} \,,\qquad \frac{\card(I)  \mathcal{H}^{d-1}(nA)}{\mathcal{H}^{d-1}(NA)\sqrt{n}}
\,\geq\,\frac{c'}{\sqrt{n}} \qquad and \qquad \frac{c'}{2\sqrt{n}}\,\geq\,
\frac{1}{\sqrt{N}}\,. $$
We conclude that for $n\geq n_1$ and $N\geq N_3(n)$, $\gamma \geq
\mathcal{H}^{d-1}(NA)/\sqrt{N}$ and then 
$$ \tau(NA,h(N)) \,\leq\, \left(\alpha \lambda_1 + (1-\alpha) \lambda_2
-\frac{1}{\sqrt{N}}\right)  \mathcal{H}^{d-1}(NA) \,, $$
as long as (\ref{eq:convex1}) and (\ref{eq:convex2}) hold.
Then, for all $n\geq n_1$ and $N\geq N_3(n)$, we have, by the FKG inequality:
\begin{align*}
\mathbb{P} \bigg(\tau(NA,h(N)) & \leq \left(\alpha \lambda_1 + (1-\alpha) \lambda_2 -
\frac{1}{\sqrt{N}}\right) \mathcal{H}^{d-1}(NA)\bigg) \\
& \,\geq\, \mathbb{P} \left(\tau(nA,h(n)) \leq (\lambda_1 - \frac{1}{\sqrt{n}})
\mathcal{H}^{d-1}(nA)\right)^{\gamma} \\
& \qquad \times  \mathbb{P} \left(\tau(nA,h(n)) \leq (\lambda_2 -
\frac{1}{\sqrt{n}}) \mathcal{H}^{d-1}(nA)\right)^{\card(I)-\gamma} p^{\card(E_0)}\,.
\end{align*}
We take the logarithm of this expression, we divide it by
$\mathcal{H}^{d-1}(NA)$, we send $N$ to infinity and then $n$ to infinity
to obtain
$$ \mathcal{I}_{\vec{v}}(\alpha \lambda_1 + (1-\alpha) \lambda_2) \,\leq\,
    \alpha \mathcal{I}_{\vec{v}} (\lambda_1) + (1-\alpha)
    \mathcal{I}_{\vec{v}} (\lambda_2) \,. $$
The convexity of $\mathcal{I}_{\vec{v}}$ is so proved.


\subsection{Continuity of $\mathcal{I}_{\vec{v}}$}
\label{cont}

Now we come back to the problem of the continuity of
$\mathcal{I}_{\vec{v}}$. Since $\mathcal{I}_{\vec{v}}$ is convex, we first
try to determine its domain. Recall that $\delta = \delta(F)=\inf \{ \lambda \,|\, \mathbb{P} (t(e)
\leq \lambda) >0 \}$. 

\noindent
\underline{$\bullet \lambda > \|\vec{v}\|_1 \delta$:} there exists $\varepsilon
>0$ such that $\lambda > (\|\vec{v}\|_1 + \varepsilon) (\delta +
2\varepsilon)$. Then there exists $n_0$ such that, for all $n\geq n_0$,
there exists a set of edges $E_0(n)$ that disconnects $\smash{(nA)_1^{h(n)}}$
    from $\smash{(nA)_2^{h(n)}}$ in $\cyl(nA,h(n))$ and such that $\smash{\card (E_0(n))
    \leq (\|\vec{v}\|_1 + \varepsilon) \mathcal{H}^{d-1}(nA)}$. We obtain
    for $n\geq n_0$
\begin{align*}
\mathbb{P} \left( \tau(nA,h(n)) \leq \left(\lambda - \frac{1}{\sqrt{n}}\right)
\mathcal{H}^{d-1}(nA) \right) & \,\geq\, \mathbb{P} \left( V(E_0(n)) \leq
\left(\lambda - \frac{1}{\sqrt{n}}\right) \mathcal{H}^{d-1}(nA) \right)\\
& \,\geq\, \mathbb{P} \left( t(e) \leq \frac{\lambda -
  1/\sqrt{n}}{\|\vec{v}\|_1 + \varepsilon}  \right)^{\lfloor (\|\vec{v}\|_1
  + \varepsilon) \mathcal{H}^{d-1}(nA) \rfloor}\,.
\end{align*}
But there exists $n_1$ large enough to have for all $n\geq n_1$, $\lambda -
1/\sqrt{n} \geq (\|\vec{v}\|_1 + \varepsilon) (\delta + \varepsilon)$, so
for all $n\geq n_0 \vee n_1$, we have
$$ \mathbb{P} \left( \tau(nA,h(n)) \leq \left(\lambda -
  \frac{1}{\sqrt{n}}\right) 
  \mathcal{H}^{d-1}(nA) \right) \,\geq\, \mathbb{P} (t(e) \leq \delta +
  \varepsilon) ^{\lfloor (\|\vec{v}\|_1 + \varepsilon) \mathcal{H}^{d-1}(nA)
  \rfloor} \,,$$
and finally
$$\mathcal{I}_{\vec{v}} (\lambda) \,\leq\, - (\|\vec{v}\|_1 + \varepsilon)
\log \mathbb{P} (t(e) \leq \delta + \varepsilon) \,<\, \infty \,. $$

\noindent
\underline{$\bullet \lambda \leq \|\vec{v}\|_1 \delta$:} for $\lambda >0$,
there exists $n_0$
such that for all $n\geq n_0$, 
$$ \frac{\tau (nA,h(n))}{\mathcal{H}^{d-1}(nA)} \,\geq\, \delta \frac{
  \mathcal{N}(nA,h(n))}{\mathcal{H}^{d-1}(nA,h(n))}  \,\geq\, \delta
  \|\vec{v}\|_1  - \frac{1}{2\sqrt{n}} \,>\, \lambda - \frac{1}{\sqrt{n}}
  \,,$$
and so for all $n\geq n_0$,
$$ \mathbb{P} \left( \tau(nA,h(n)) \leq \left(\lambda -
\frac{1}{\sqrt{n}}\right) 
\mathcal{H}^{d-1}(nA) \right) \,=\,0 \,.$$
The same result is true for $\lambda = 0$.
We obtain that $\mathcal{I}_{\vec{v}} (\lambda) = +\infty$.

Now, we know that $\mathcal{I}_{\vec{v}}$ is convex and finite on $]\delta
\|\vec{v}\|_1, +\infty[$ so it is continuous on $]\delta
\|\vec{v}\|_1, +\infty[$, and it is infinite on $[0,\delta \|\vec{v}\|_1]$.

\begin{rem}
\label{rem:Inu}
The only point we didn't study is the behaviour of
the function 
near $\delta \|\vec{v}\|_1$. In fact, we will eventually change the value
of $\mathcal{I}_{\vec{v}} (\delta \|\vec{v}\|_1)$ to obtain a lower
semicontinuous function. Moreover, the fact that $\mathcal{I}_{\vec{v}}
(\delta \|\vec{v}\|_1) = +\infty$ even if there exists an atom of the law
of $t(e)$ at $\delta$ is linked with the fact that we added a term
$1/\sqrt{n}$ and not with the behaviour of $\mathbb{P}
(\tau(nA,h(n)) \leq \delta \|\vec{v}\|_1 \mathcal{H}^{d-1}(nA))$. This
remark can be illustrated by an example in dimension $2$: let $A=
[-1/2,1/2]\times \{1/2\} $. Here $\vec{v} = (0,1)$ so $\|\vec{v}\|_1 =1$. We
consider a law of capacities with an atom at $\delta$.
\begin{figure}[ht!]
\centering

\begin{picture}(0,0)%
\includegraphics{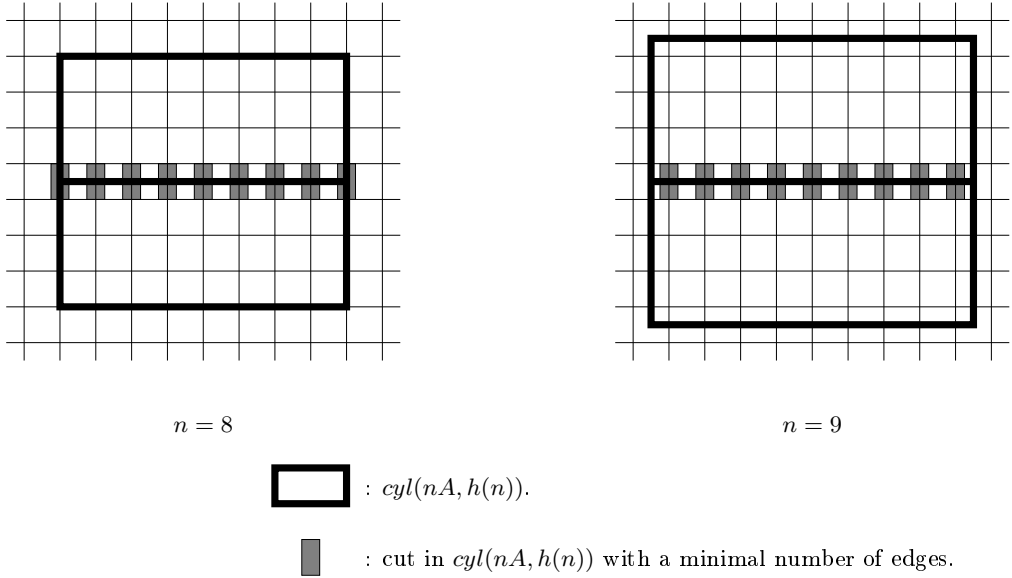}%
\end{picture}%
\setlength{\unitlength}{2960sp}%
\begingroup\makeatletter\ifx\SetFigFont\undefined%
\gdef\SetFigFont#1#2#3#4#5{%
  \reset@font\fontsize{#1}{#2pt}%
  \fontfamily{#3}\fontseries{#4}\fontshape{#5}%
  \selectfont}%
\fi\endgroup%
\begin{picture}(8424,4824)(1339,-5173)
\put(3001,-3961){\makebox(0,0)[b]{\smash{{\SetFigFont{9}{10.8}{\rmdefault}{\mddefault}{\updefault}{\color[rgb]{0,0,0}$n=8$}%
}}}}
\put(8101,-3961){\makebox(0,0)[b]{\smash{{\SetFigFont{9}{10.8}{\rmdefault}{\mddefault}{\updefault}{\color[rgb]{0,0,0}$n=9$}%
}}}}
\put(4351,-4486){\makebox(0,0)[lb]{\smash{{\SetFigFont{9}{10.8}{\rmdefault}{\mddefault}{\updefault}{\color[rgb]{0,0,0}: $cyl(nA, h(n))$.}%
}}}}
\put(4351,-5086){\makebox(0,0)[lb]{\smash{{\SetFigFont{9}{10.8}{\rmdefault}{\mddefault}{\updefault}{\color[rgb]{0,0,0}: cut in $cyl(nA, h(n))$ with a minimal number of edges.}%
}}}}
\end{picture}%

\caption{Examples of cuts.}
\label{deltav}
\end{figure}
We remark (see Figure \ref{deltav}) that $\mathcal{N} ((2n+1)A,2n+1) =
2n+1$. Moreover, there exists a unique cut $E_0(2n+1)$ in $\cyl ((2n+1)A, 2n+1)$
composed by $2n+1$ edges (see it on the Figure). So we have
$$ \mathbb{P} (\tau((2n+1)A, 2n+1) \leq (2n+1) \delta) \,=\, \mathbb{P}
(V(E_0(2n+1))=(2n+1) \delta ) \,=\, \mathbb{P} (t(e)=\delta )^{2n+1} $$
and
$$ \lim_{n\rightarrow \infty} \frac{-1}{2n+1} \log \mathbb{P}
(\tau((2n+1)A,2n+1 ) \leq (2n+1) \delta ) \,=\, - \log \mathbb{P} (t(e) =
\delta ) \,<\, \infty \,.$$
We also remark that $\mathcal{N} (2nA, 2n) \,=\, 2n+1$ because a cut in $\cyl
(2nA, 2n)$ must contain a vertical edge of first coordinate $i$ for
$i=0,...,2n$. Then we have
$$ \mathbb{P} (\tau(2nA,2n) \leq 2n \delta) \,=\, 0 $$
and
$$ \lim_{n\rightarrow \infty} \frac{-1}{2n} \log \mathbb{P}
(\tau(2nA,2n ) \leq 2n \delta ) \,=\, + \infty \,. $$
This example shows that the behaviour of $\mathbb{P}
(\tau(nA,h(n)) \leq \delta \|\vec{v}\|_1 \mathcal{H}^{d-1}(nA))$ is not
clear, and we will avoid the problem by taking later at $\|\vec{v}\|_1
\delta$ the value of the limit
$$ \lim_{\lambda > \|\vec{v}\|_1 \delta, \lambda
  \rightarrow \|\vec{v}\|_1 \delta} \mathcal{I}_{\vec{v}} (\lambda)$$
instead of $\mathcal{I}_{\vec{v}} (\|\vec{v}\|_1 \delta)$.
\end{rem}


\subsection{Positivity of $\mathcal{I}_{\vec{v}}$}
\label{positivity}

>From now on we need the assumptions $\mathbf{(F1)}$ and
$\mathbf{(F2)}$, i.e  $F(0) <1-p_c(d)$, and $F$ admits a moment of
order 1. It is an immediate consequence of Theorem \ref{thm:llntau} that
$\mathcal{I}_{\vec{v}}$ is equal to zero on $]\nu(\vec{v}), +\infty[$, and
Theorem \ref{devinf} implies immediately too that $\mathcal{I}_{\vec{v}}$ is
strictly positive on $[0, \nu(\vec{v})[$ if $\nu(\vec{v}) >0$.

\begin{rem}
We did not study the function $\mathcal{I}_{\vec{v}}$ at $\nu(\vec{v})$, i.e., if
$\mathcal{I}_{\vec{v}} (\nu(\vec{v})) =0$ or not. If $\nu(\vec{v}) > \delta \|\vec{v}\|_1$, then $\mathcal{I}_{\vec{v}}$ is
continuous at $\nu(\vec{v})$ and so $\mathcal{I}_{\vec{v}}(\nu(\vec{v}))
=0$. If $\nu(\vec{v}) = \delta \|\vec{v}\|_1$, the value of
$\mathcal{I}_{\vec{v}} (\nu(\vec{v}))$ is not relevant for the
understanding of the system as explained in Remark
\ref{rem:Inu}. Finally, Proposition
\ref{thmnu} gives a sufficient condition to have $\nu(\vec{v}) >\delta
\|\vec{v}\|_1$, and this condition is also necessary when
$\delta=0$.
\end{rem}








\subsection{Upper large deviations for $\tau$}
\label{secdevsup}

We will need the following result to prove the large deviation principle
for $\tau$ in the next section:
\begin{lem}
\label{devsup}
Suppose that $\mathbf{(H1)}$ and $\mathbf{(F5)}$ hold. Then we have, for all $\lambda > \nu(\vec{v})$,
\begin{equation}
\label{decroissance}
 \lim_{n\rightarrow \infty} \frac{1}{\mathcal{H}^{d-1}(nA)} \log
\mathbb{P} \left[ \frac{\tau(nA, h(n))}{\mathcal{H}^{d-1}(nA)} \geq \lambda
  \right] \,=\, -\infty \,.
\end{equation}
\end{lem}

We do not prove Lemma \ref{devsup} here. The proof is an adaptation of
 section 3.7 in \cite{TheretUpper}, that proves that the upper large
deviations for $\phi(nA, h(n))/\H^{d-1}(nA)$ in straight boxes are of
volume order. It is written completely in \cite{Theret:uppertau}, where
other assumptions on $F$ are also considered. We describe
here only the two adaptations required to get Lemma \ref{devsup} from
the proof in \cite{TheretUpper}. The proof for $\phi$ is based on a
comparison between the variable $\phi(NA, h(N))$ in a big cylinder, and
the minimum over $h(N)/h(n)$ possible choices of sums of
$\mathcal{H}^{d-1}(NA)/\H^{d-1}(nA)$ independent variables equal in law with
$\tau(nA, h(n))$, where $n$ is small compared to $N$. This comparison is
obtained by dividing the big cylinder $\cyl(NA, h(N))$ into $h(N)/h(n)$
slabs, and diving each slab in $\mathcal{H}^{d-1}(NA)/\H^{d-1}(nA)$
translates of $\cyl(nA, h(n))$. Then, in any fixed slab, if we glue
together cutsets in the small cylinder of size $n$, we can construct a
cutset in $\cyl(NA,h(N))$. There are two difficulties to replace $\phi(NA,
h(N))$ by $\tau(NA,h(N))$ in this construction, and to consider potentially
tilted cylinders. First, the fact that the cylinders we consider may be
tilted implies a default of subadditivity of the variable $\tau$, so we have to
add edges between the small cylinders of size $n$ to glue together the
different cutsets, and we have to control the number of the edges we
must consider. Then, when the small cutsets are glued together, they form a
set of edges that cuts the top from the bottom of $\cyl(NA, h(N))$. It
remains to link this cutset to the boundary of $NA$ to obtain a cutset
corresponding to the variable $\tau(NA,h(N))$. To obtain a control on the
number of edges we must add at this step, we have to consider only slabs
whose distance to $NA$ is negligible compared to $N$. Using
Cram\'er Theorem for each possible sum of independent variables in a slab, and
optimizing over the possible choices of slab, we obtain the desired result.


\begin{rem}
\label{remcondF5}
For the variable $\phi$, it suffices to have one exponential moment for
the law $F$ to obtain this speed of decay (see \cite{TheretUpper}). For
$\tau$, one exponential moment is not a sufficiently strong
condition. Consider for example an exponential law of parameter $1$ for the
capacities of the edges. We know that $\mathbb{E} (\exp (\gamma t))
<\infty$ for all $\gamma <1$. Let $x_0$ be a fixed point of the
boundary $\partial (nA)$. There are, at distance at most $4d$ of
$x_0$, one vertex of  $(nA)_1^{h(n)}$ and another of
$(nA)_2^{h(n)}$. Let $\gamma$ be some smallest path in $\cyl(nA,h(n))$
joining those two vertices. Its length is at most some constant
$R(d)$, and we know that every set of edges that
cuts $(nA)_1^{h(n)}$ from $(nA)_2^{h(n)}$ in $\cyl (nA, h(n))$ must contain
one of the edges of $\gamma$. The probability that all of them  have a capacity
bigger than $\lambda \mathcal{H}^{d-1}(nA)$ for some $\lambda > \nu(\vec{v})$,
and therefore that $\tau (nA, h(n))$ is bigger than $\lambda
\mathcal{H}^{d-1}(nA)$, is greater than $\exp (-R(d) \lambda
\mathcal{H}^{d-1}(nA) )$. Then the property (\ref{decroissance}) cannot hold.
\end{rem}

\begin{rem}
It is also proved in \cite{Theret:uppertau} that if the capacity of the edges
is bounded, the upper large deviations are of order $n^{d-1}\min (n,
h(n))$, and this is the right order of the upper large deviations in this case.
\end{rem}

\subsection{Proof of Theorem \ref{thmpgd}}
\label{pgd}

We define the function $\mathcal{J}_{\vec{v}}$ on
$\mathbb{R}^+$ by
$$ \mathcal{J}_{\vec{v}} (\lambda) \,=\, \left\{ \begin{array}{ll}
\mathcal{I}_{\vec{v}} (\lambda) & if \,\, \lambda \leq \nu(\vec{v}) \,\,
and\,\, \lambda \neq \|\vec{v}\|_1 \delta \,,\\
 \lim_{\mu > \|\vec{v}\|_1 \delta, \mu \rightarrow \|\vec{v}\|_1 \delta}
 \mathcal{I}_{\vec{v}} (\mu) & if\,\,\lambda = \|\vec{v}\|_1 \delta \,, \\
 +\infty & if\,\, \lambda > \nu(\vec{v})\,. \end{array} \right.$$
The study of the function $\mathcal{I}_{\vec{v}}$ made previously and the
construction of $\mathcal{J}_{\vec{v}}$ gives us immediately that the
function $\mathcal{J}_{\vec{v}}$ is a good rate function.  As soon as
  we know that the upper large deviations are of order bigger than the
  lower large deviations, the techniques we will use to prove the large
  deviation principle are standard (see for example \cite{Cerf:StFlour}).

\noindent
{\bf $\bullet$ Lower bound}

We have to prove that for all open subset $\mathcal{O}$ of $\mathbb{R}^+$,
$$ \liminf_{n\rightarrow \infty}
\frac{1}{\mathcal{H}^{d-1}(nA)} \log \mathbb{P} \left[
  \frac{\tau(nA,h(n))}{\mathcal{H}^{d-1}(nA)} \in \mathcal{O} \right] \,\geq \, -
\inf_{\mathcal{O}} \mathcal{J}_{\vec{v}} \,.$$
Classically, it suffices to prove the local lower bound:
$$\forall \alpha \in \mathbb{R}^+\,,\,\, \forall \varepsilon>0 \qquad \liminf_{n\rightarrow \infty}
\frac{1}{\mathcal{H}^{d-1}(nA)} \log \mathbb{P} \left[
  \frac{\tau(nA,h(n))}{\mathcal{H}^{d-1}(nA)} \in ]\alpha -
\varepsilon, \alpha + \varepsilon[ \right] \,\geq \, -
\mathcal{J}_{\vec{v}}(\alpha)  \,.$$
If
$\mathcal{J}_{\vec{v}}(\alpha) = +\infty$, the result is
trivial. Otherwise, suppose
$\mathcal{J}_{\vec{v}}(\alpha) < +\infty$. The function
$\mathcal{I}_{\vec{v}}$ is convex, equal to zero on
$[\nu(\vec{v}), +\infty[$, positive on $[0,\nu(\vec{v})[$ and finite on
         $]\|\vec{v}\|_1 \delta,
      +\infty] $. Then $\mathcal{I}_{\vec{v}} $ is strictly decreasing on
      $]\|\vec{v}\|_1 \delta,
      \nu(\vec{v})] $, and so is $\mathcal{J}_{\vec{v}}$ (because $\mathcal{I}_{\vec{v}} =\mathcal{J}_{\vec{v}}$ on $]\|\vec{v}\|_1 \delta,
      \nu(\vec{v})] $). Yet $\mathcal{J}_{\vec{v}}(\alpha) < +\infty$ implies that $\alpha \in ]\|\vec{v}\|_1 \delta,
      \nu(\vec{v})]$ or $\alpha = \|\vec{v}\|_1 \delta$ and $\mathcal{J}_{\vec{v}}(\|\vec{v}\|_1 \delta) < +\infty$. In both cases, we so obtain that $\mathcal{J}_{\vec{v}}(\alpha) <\mathcal{J}_{\vec{v}}(\alpha -
            \varepsilon/2) $. Then the following inequality, true
            for $n> 4/\varepsilon^2$,
$$ \mathbb{P} \left[  \frac{\tau(nA,h(n))}{\mathcal{H}^{d-1}(nA)} \in
              ]\alpha - \varepsilon,\alpha+\varepsilon [ \right] \,\geq \,
              \mathbb{P} \left[
              \frac{\tau(nA,h(n))}{\mathcal{H}^{d-1}(nA)} \leq \alpha -
              \frac{1}{\sqrt{n}} \right] -
              \mathbb{P} \left[
              \frac{\tau(nA,h(n))}{\mathcal{H}^{d-1}(nA)} \leq \alpha
              -\frac{\varepsilon}{2} - \frac{1}{\sqrt{n}} \right] $$
leads to
$$ \liminf_{n\rightarrow \infty}
\frac{1}{\mathcal{H}^{d-1}(nA)} \log \mathbb{P} \left[
  \frac{\tau(nA,h(n))}{\mathcal{H}^{d-1}(nA)} \in ]\alpha -
\varepsilon, \alpha + \varepsilon[ \right] \,\geq \, -
\mathcal{J}_{\vec{v}}(\alpha)  \,. $$

\noindent
{\bf $\bullet$ Upper bound}

We have to prove that for all closed subset $\mathcal{F}$ of $\mathbb{R}^+$
$$ \limsup_{n\rightarrow \infty}
\frac{1}{\mathcal{H}^{d-1}(nA)} \log \mathbb{P} \left[
  \frac{\tau(nA,h(n))}{\mathcal{H}^{d-1}(nA)} \in \mathcal{F} \right] \,\leq \, -
\inf_{\mathcal{F}} \mathcal{J}_{\vec{v}} \,.$$
Let $\mathcal{F}$ be a closed subset of $\mathbb{R}^+$. If $\nu(\vec{v})
\in \mathcal{F}$, the result is obvious. We suppose now that $\nu(\vec{v})
\notin \mathcal{F}$. We consider $\mathcal{F}_1 = \mathcal{F} \cap
       [0,\nu(\vec{v})]$ and $\mathcal{F}_2 = \mathcal{F} \cap
          ]\nu(\vec{v}),+\infty[$. Let $f_1 = \sup \mathcal{F}_1$ ($f_1
            < \nu(\vec{v})$ because $\mathcal{F}$ is closed) and $f_2 = \inf \mathcal{F}_2$ ($f_2 >
            \nu(\vec{v})$ for the same reason). Then,
\begin{align*}
\limsup_{n\rightarrow \infty} \frac{1}{\mathcal{H}^{d-1}(nA)} &  \log
\mathbb{P} \left[ \frac{\tau(nA,h(n))}{\mathcal{H}^{d-1}(nA)} \in
  \mathcal{F} \right]\\ & \,\leq\,
\limsup_{n\rightarrow \infty} \frac{1}{\mathcal{H}^{d-1}(nA)} \log
\left( \mathbb{P} \left[ \frac{\tau(nA,h(n))}{\mathcal{H}^{d-1}(nA)}
  \leq f_1 \right] + \mathbb{P} \left[
  \frac{\tau(nA,h(n))}{\mathcal{H}^{d-1}(nA)} \geq f_2 \right]
\right)\,.\\
\end{align*}
 We know that:
\begin{align*}
\limsup_{n\rightarrow \infty}
  \frac{1}{\mathcal{H}^{d-1}(nA)} &\log \mathbb{P} \left(
    \frac{\tau(nA,h(n))}{\mathcal{H}^{d-1}(nA)} \leq f_1 \right) \\ & \,\leq\,
  \lim_{\eta \rightarrow 0}\limsup_{n\rightarrow \infty}
  \frac{1}{\mathcal{H}^{d-1}(nA)} \log \mathbb{P} \left(
    \frac{\tau(nA,h(n))}{\mathcal{H}^{d-1}(nA)} \leq f_1 + \eta -
    \frac{1}{\sqrt{n}} \right) \\
 & \,=\, -\lim_{\eta \rightarrow 0} \mathcal{I}_{\vec{v}} (f_1 + \eta) 
\,=\, -\mathcal{J}_{\vec{v}} (f_1) \,,
\end{align*}
and since $\mathcal{J}_{\vec{v}}$ is non-increasing on $[0,\nu_{\vec{v}}]$
and the upper large deviations of $\tau(nA, h(n))$ are of order bigger than
$n^{d-1}$, we
obtain:
$$ \limsup_{n\rightarrow \infty} \frac{1}{\mathcal{H}^{d-1}(nA)} \log
\mathbb{P} \left[ \frac{\tau(nA,h(n))}{\mathcal{H}^{d-1}(nA)} \in
  \mathcal{F} \right] \,\leq\, -\mathcal{J}_{\vec{v}} (f_1) \,=\, -\inf
_{\mathcal{F}} \mathcal{J}_{\vec{v}}\,.$$

\subsection{Large deviation principle for $\phi$ in small boxes}
\label{subsec:ldpflatphi}
In this section, we shall prove Corollary \ref{corthmpgd}, i.e., under the  assumption that $\lim_{n\rightarrow
  \infty}h(n)/n = 0$, the sequence 
$$\left(\frac{\phi(nA, h(n))}{\mathcal{H}^{d-1}(nA)}, n\in \mathbb{N}\right)$$
satisfies the same large deviation principle as $(\tau(nA,
h(n))/\mathcal{H}^{d-1}(nA), n\in \mathbb{N})$.

We will use a result of exponential equivalence. For $(X_n)$ and $(Y_n)$
two sequences of random
variables defined on the same probability space $(\Omega, \mathcal{A},
\mathbb{P})$, and for a given speed function $v(n)$ which goes to infinity
with $n$, we say that $(X_n)$ and $(Y_n)$ are exponentially equivalent with
regard to $v(n)$ if and only if for all positive $\varepsilon$ we have
$$ \limsup_{n\rightarrow \infty} \frac{1}{v(n)} \log \mathbb{P} \left( |X_n
- Y_n | \geq \varepsilon \right) \,=\, -\infty \,.  $$
The following result is classical in large deviations theory (see
\cite{Dembo-Zeitouni}, Theorem 4.2.13):
\begin{thm}
\label{expeq}
Let $(X_n)$ and $(Y_n)$ be two sequences of random
variables defined on the same probability space $(\Omega, \mathcal{A},
\mathbb{P})$. If $(X_n)$ satisfies a large deviation principle of speed
$v(n)$ with a good rate function, and if $(X_n)$ and $(Y_n)$ are
exponentially equivalent with regard
to $v(n)$, then $(Y_n)$ satisfies the same large deviation principle as
$(X_n)$.
\end{thm}
We will prove that the sequences $(\phi(nA,h(n))/\mathcal{H}^{d-1}(nA))$ and
$(\tau(nA,h(n))/\mathcal{H}^{d-1}(nA))$ are exponentially equivalent with
regard to $\mathcal{H}^{d-1}(nA)$ under the assumptions that there exist
exponential moment of the law of capacity of all orders and for height
functions $h$ satisfying $\lim_{n\rightarrow \infty} h(n)/n =0$.

We take a hyperrectangle $A$ and use the same notations as in section \ref{sectiondevinf}.
Let $\zeta\geq 2d$, and $n$ large enough such that the sides of $nA$ have
length bigger than $\zeta$.
Let $E_1^+ $ be the set of the edges that belong to $\mathcal{E}_1^+$
defined as
$$ \mathcal{E}_1^+ \,=\, \mathcal{V}(\cyl( \partial (nA) ,h(n)), \zeta) \cap\cyl(nA,h(n)) \,. $$
We have for all $n \geq p$
$$\phi(nA,h(n)) \,\leq\, \tau (nA, h(n))
\,\leq\, \phi(nA,h(n)) + V(E_1^+) \,. $$
Thus for all $\varepsilon >0$, for all $n\geq p$, we obtain
$$ \mathbb{P} \left(\Big\arrowvert
\frac{\phi(nA,h(n))}{\mathcal{H}^{d-1}(nA)} - \frac{\tau(nA,
  h(n))}{\mathcal{H}^{d-1}(nA)} \Big\arrowvert \geq \varepsilon \right)
\,\leq\, \mathbb{P}\left( V(E_1^+) \geq \varepsilon \mathcal{H}^{d-1}(nA)
\right) \,.$$
We know that there exists a constant $C^+$ such that
$$\card(E_1^+) \leq C^+ n^{d-2} h(n) \,,$$
so for all $\varepsilon>0$, for all $\gamma >0$, for a family $(t_k)$ of
independent variables with the same law as the capacities of the edges, we have
\begin{align*}
\mathbb{P} \left[ V(E_1^+) \geq \varepsilon \mathcal{H}^{d-1}(nA) \right] &
\,\leq\, \mathbb{P} \left[ \sum_{k=1}^{C^+ n^{d-2} h(n)} t_k \geq \varepsilon
  \mathcal{H}^{d-1}(nA) \right] \\
& \,\leq\, \mathbb{E} (e^{\gamma t})^{C^+ n^{d-2}h(n)} \exp \left(-\gamma
\varepsilon \mathcal{H}^{d-1}(nA) \right) \\
& \,\leq\, \exp \left(-\mathcal{H}^{d-1}(nA) \left( \gamma \varepsilon -
C^+ \frac{n^{d-2} h(n)}{\mathcal{H}^{d-1}(nA)} \log
  \mathbb{E} (e^{\gamma t})\right)\right) \,.\\
\end{align*}
For a fixed $R>0$, we can choose $\gamma$ large enough to have $\gamma \varepsilon
\geq 2R$, and also there exists $n_2$ such that for all $n\geq n_2$ we have 
$$ C^+  \frac{n^{d-2} h(n)}{\mathcal{H}^{d-1}(nA)} \log
  \mathbb{E} (e^{\gamma t}) \,\leq\, R  \,,$$
so for all $R>0$
$$ \limsup_{n\rightarrow \infty} \frac{1}{\mathcal{H}^{d-1}(nA)} \log
\mathbb{P } \left[V(E_1^+) \geq \varepsilon \mathcal{H}^{d-1}(nA)\right] \,\leq
\,- R $$ 
and then
$$ \limsup_{n\rightarrow \infty} \frac{1}{\mathcal{H}^{d-1}(nA)} \log
\mathbb{P } \left[V(E_1^+) \geq \varepsilon \mathcal{H}^{d-1}(nA)\right] \,=\,
-\infty \,.$$
We obtain immediately that $(\phi(nA,h(n))/\mathcal{H}^{d-1}(nA))$ and
$(\tau(n A, h(n))/\mathcal{H}^{d-1}(nA))$ are exponentially equivalent
with regard to $\mathcal{H}^{d-1}(nA)$, and so by Theorem \ref{expeq}, 
$(\phi(nA,h(n))/\mathcal{H}^{d-1}(nA))$ satisfies the same large deviation
principle as $(\tau(nA,h(n))/\mathcal{H}^{d-1}(nA))$.

\section[LLN, LDP and lower large deviations for $\phi$]{Law of large
  numbers, large deviation principle and lower large deviations
  for $\phi$ in straight boxes}
\label{sec:pgdphi}

The main work is done in section \ref{subsec:pgdphistraight}, where
one proves that $\phi$ and $\tau$ share the same rate function in
straight boxes. Then, the law of large numbers is proven in section
\ref{sec:llnphistraight}. The large deviation principle is proven in section \ref{secfinal}, as
well as the deviation inequality from $\nu$ (Theorem \ref{devinfphi}).


\subsection{Comparison between $\phi$ and $\tau$}
\label{subsec:pgdphistraight}

We prove in this section that under hypotheses $\mathbf{(F1)}$,
$\mathbf{(F4)}$, $\mathbf{(H1)}$ and $\mathbf{(H2)}$, the
lower large deviations of $\phi(nA, h(n))$ and $\tau(nA, h(n))$ are of
the same exponential order. The
following proposition is the key to prove both Theorem \ref{thm:llnphistraight}
and Theorem \ref{thmpgdphi}.
\begin{prop}
\label{prop:fntauxphi}
Suppose that $\mathbf{(F1)}$, $\mathbf{(F4)}$, $\mathbf{(H1)}$ and
$\mathbf{(H2)}$ hold. Let $A$ be a non-degenerate straight hyperrectangle. Then, for every $\lambda$ in  $\RR^+$,
$$\lim_{n\rightarrow \infty}\frac{-1}{\mathcal{H}^{d-1}(nA)} \log
\mathbb{P} \left[ \phi (nA,h(n)) \leq \left(\lambda-\frac{1}{\sqrt{n}}\right)
  \mathcal{H}^{d-1}(nA) \right]\,=\,\mathcal{I}_{\vec{v}}(\lambda)\;, $$
where $\vec{v}=(0,\ldots,0,1)$.
\end{prop}
\begin{dem}
Since
$\phi (nA,h(n)) \leq \tau (nA,h(n))$, we only need to show that:
$$\liminf _{n\rightarrow\infty}\frac{-1}{\mathcal{H}^{d-1}(nA)} \log
\mathbb{P} \left[ \phi (nA,h(n)) \leq \left(\lambda-\frac{1}{\sqrt{n}}\right) \mathcal{H}^{d-1}(nA) \right]\geq\mathcal{I}_{\vec{v}}(\lambda)\;.$$
To shorten the notations, we shall suppose that
$$A\,=\,[0,1]^{d-1}\times\{0\}\,,$$
the general case of a straight hyperrectangle being handled exactly along
the same lines. Notice that $\mathcal{H}^{d-1}(nA)=n^{d-1}$. As in
section \ref{sec:concentration}, we shall write $\phi_n$ instead of
$\phi(nA,h(n))$, and  denote by $E_{\phi_n}$ a cut whose capacity achieves
the minimum in the dual definition (\ref{eq:maxflowmincut}) of $\phi_n$.

\emph{The idea of the proof is the following}. The minimal cut $E_{\phi_n}$ has a certain intersection with the sides of the cylinder
$\mbox{cyl}(A,h)$. Thanks to Zhang's result, Theorem
\ref{thmpgd}, and after having eventually reduced a
little the cylinder, one can prove that the intersection of $E_{\phi_{n}}$ with
the sides of this reduced cylinder has less than $Cn^{d-1}/n^{1/3}$ edges with
very high probability (here $C$ is a constant). This shows that (with very high probability) $\phi_n$ is larger than the minimum of a
collection of random variables $(\tau_F)_{F\in I_n}$, where $F$ designs a
possible \emph{trace} of $E_{\phi_{n}}$, i.e., its intersection  with the
sides of the reduced cylinder, and where $I_n$ is the set of all the
possible choices for $F$. Since $E_{\phi_{n}}$ itself has less than
$Cn^{d-1}$ edges, and since it is connected (in the dual sense), a trivial bound
for the cardinal of $I_n$ is roughly:
$$ \card (I_n) \,\leq\, h(n)(C'n^{2d-3})^{C n^{d-1}/n^{1/3}}\,.$$
The important point here is that $\log \card (I_n)$ is small
compared to $n^{d-1}$. Having done this, a subadditive argument using
symmetries can be
performed to show that in fact the smallest $\tau_F$ (in distribution) behaves
essentially like $\tau(nA,h(n))$, which has $\mathcal{I}_{\vec{v}}$ as a rate
function.

\emph{Now, we turn to a formal proof}.
In the sequel, we shall suppose that $n$ is large enough to ensure that
$$\log h(n)\,\leq\, n^{d-1} \,, \,\, nl_{\min}(A) \,\geq\, t_0 \,\,
\textrm{and} \,\, h(n) >2\sqrt{d}\,,$$
where $t_0$ is defined in Proposition~\ref{prop:5.8Zhangbis}. We consider
$\gamma >0$ such that $\mathbb{E} (\exp(\gamma t(e)))<\infty$. Let $E_n$ be
the cutset defined by
$E_n \,=\, \{ e=\langle x,y \rangle \in \mathbb{E}^d  \,|\,
x\in nA\,\, and \,\, y_d=1\} \,. $ Notice that for $n$ large enough,
$$ \card(E_n) \,\leq\, 2n^{d-1} \,. $$
For a fixed $L$, and for $\eps$, $C_1$ and $C_2$ as in
Proposition~\ref{prop:5.8Zhangbis}, using this proposition we obtain:
\begin{align*}
\mathbb{P} \left( \card (E_{\phi_n}) \geq Ln^{d-1} \right) & \,\leq\,
\mathbb{P} \left( \card (E_{\phi_n}) \geq Ln^{d-1} \mbox{ and } \phi_n
  \leq \eps L n^{d-1} \right) + \mathbb{P} \left( \phi_n \geq \eps L
  n^{d-1} \right)\\ 
&\,\leq\, C_1 h(n) \exp \left[-C_2 L n^{d-1}\right] + \mathbb{P}
\left( V(E_n) \geq \eps L n^{d-1}  \right)\\
& \,\leq\,C_1 \exp \left[-(C_2 L -1)n^{d-1}\right] + \mathbb{P}
\left(\sum_{j=1}^{2n^{d-1}} t(e_j) \geq \eps L n^{d-1}
\right) \\
&\,\leq\, C_1 \exp \left[-(C_2 L -1)n^{d-1}\right]  + \exp \left[
-\left( \gamma \eps L - 2\log \mathbb{E}
  \left(e^{\gamma t(e)} \right)\right)n^{d-1}\right]\,.
\end{align*}
Thus, there exist constants $\beta(F,d)$ and $C_i'(F,d)$ for $i=1,2$ such
that for all $L\geq \beta$ and every $n$, we have
$$ \mathbb{P} \left( \card (E_{\phi_n}) \geq Ln^{d-1} \right) \,\leq\,
C_1'(F,d) e^{-C'_2(F,d) L n^{d-1}}\,. $$
We fix a real number $L\geq \beta$ to be chosen later. Define
\begin{align*}
\phi_{L,n}\,=\,\min \{V(E)\,| \, & E\mbox{ is  a
}(B(nA,h(n)),T(nA,h(n)))\mbox{-cut in }cyl(nA,h(n)) \\
& \,\,\,\, \mbox{ and }\card
(E)\leq L n^{d-1}\} \;.
\end{align*}
Thus,
\begin{equation}
\label{eq:moinsdaretes}
\PP(\phi_n\leq \lambda n^{d-1})\,\leq\, \PP(\phi_{L,n}\leq \lambda
n^{d-1})+C'_1e^{-C'_2Ln^{d-1}}\;.
\end{equation}
We shall now concentrate on the first summand in the right-hand side of the last inequality. Let $\psi(n)=\lceil n^{1/3}\rceil$. For any $k$ in $\{1,\ldots ,\psi(n)\}$, define
$$A_{n,k}=[k,n-k]^{d-1}\times \{0\}\;,$$
$$B_{n,k}=[k,n-k]^{d-1}\times [-h(n),h(n)]\;,$$
and
$$S_{n,k}= \partial\left([k,n-k]^{d-1}\right)\times [-h(n),h(n)]\;.$$
In order to perform the announced subadditive argument, we shall need to patch together
two cuts of neighbouring
boxes which share a same trace in the intersection of these boxes. It is not so trivial to show that one obtains a cut doing
so. This is why we shall impose a kind of ``connection trace'' on the sides of
the box, which remembers if a vertex of the side is connected to the top or
the bottom of the cylinder once the cut $E_{\phi_{n}}$ has been removed. Let
  us precise the needed definitions. If $X$ is a subset of vertices of a subgraph $G$ of $\ZZ^d$, we denote by $C_{G}(X)$ the union of all the connected components of $G$ intersecting $X$. If $v$ is a vertex of $G$, we write $C_{G}(v)$ instead of $C_{G}(\{v\})$. We shall say that a function $x$ from some $S_{n,k}$ to $\{0,1,2\}$ is a \emph{weak connection function} for $F$ in $S_{n,k}$ if for every $u$ and $v$ in $S_{n,k}$,
$$u\in C_{S_{n,k}\smallsetminus F}(v)\,\Rightarrow\, x(u)=x(v)\;.$$
If $E$ cuts $B(A_{n,k},h(n))$ from $T(A_{n,k},h(n))$ in $B_{k,n}$, we define $x_E$, \emph{the connection function of $E$ (in $B_{k,n}$)} as follows:
$$\forall u\in \mbox{cyl}(A_{n,k},h(n)) \qquad
x_E(u)=\left\lbrace\begin{array}{l}1\mbox{ if }u\in C_{ B_{k,n}\setminus
  E}(T(A_{n,k},h(n)))\;,\\ 0\mbox{ if }u\in C_{ B_{k,n}\setminus
  E}(B(A_{n,k},h(n)))\;,\\ 2 \mbox{ else .}\end{array}\right.\;$$
Clearly, $\tilde{x}_E$, the restriction of $x_E$ to $S_{n,k}$ is a weak connection function for $E\cap S_{n,k}$. Then, define the following set of ``good'' couples $(F,x)$ of a trace $F$ and a weak connection function $x$:
$$\begin{array}{rl}
\displaystyle I_n=\bigcup_{k=1}^{\psi(n)}\bigcup_{h=-h(n)}^{h(n)-L
  n^{d-1}}&\left\{(F,x)\,|\, F \subset \EE^d \cap\partial\left([k,n-k]^{d-1}\right)\times [h,h+Ln^{d-1}] ,\;\card (F)\leq \frac{L
    n^{d-1}}{\psi(n)},\right.\;\\
& \left.\phantom{\bigcup_{k=1}^{\psi(n)}} x\mbox{ is a weak connection function  for }F\mbox{ in }S_{n,k}\right\}\;.
\end{array}$$
If $F$ satisfies the conditions in the above definition, then there are at most  $2^dL n^{d-1}/\psi(n)$ distinct connected components in $S_{n,k}\smallsetminus F$. Thus, for a fixed $F$, there are at most $3^{2^dL n^{d-1}/\psi(n)}$ distinct weak configuration functions $x$ such that $(F,x)$ belongs to $I_n$. Thus, there is a constant $C_3$, which depends only on $d$, such that
\begin{equation}
\label{eqIn}
(2h(n)+1)\,\leq\, \card (I_n) \,\leq\,  2h(n)\psi(n)(C_3Ln^{2d-3})^{L n^{d-1}/\psi(n)}\;.
\end{equation}
On the other hand, define, for $(F,x)$ in $I_n$ and $k$ such that $F\subset S_{n,k}$,
\begin{eqnarray*}
\mathcal{C}_{F,x}&=&\left\{E\subset \EE^d \,|\,E\mbox{ is  a
}(B(A_{n,k},h(n)),T(A_{n,k},h(n)))\mbox{-cut in }B_{n,k}\, ,\right.\\ &&
\,\,\,\, \left. E\cap S_{n,k}=F ,\;\card (E)\leq L n^{d-1}\mbox{ and }\tilde{x}_E=x\right\}\;,
\end{eqnarray*}
and
$$\tau_{(F,x)}\,=\,\min \left\{V(E)\,|\, E\in\mathcal{C}_{F,x}\right\}\;.$$
We claim that
\begin{equation}
\label{eqphitau}\phi_{L,n}\,\geq\, \min_{(F,x)\in I_n}\tau_{(F,x)} \;.
\end{equation}
To see why (\ref{eqphitau}) is true, notice that for any $k$ in $\{1,\ldots
,\psi(n)\}$, $E_{\phi_{L,n}}\cap
B_{n,k}$ cuts $B(A_{n,k},h(n)$ from $T(A_{n,k},h(n))$ in $B_{n,k}$,
and has less than $L n^{d-1}$ edges. $E_{\phi_{L,n}}$ is connected in
the dual sense (see the proof of Lemma 12 in \cite{Zhang07}), and has less than $Ln^{d-1}$
edges. Then there is an $h$ such that $E_{\phi_{L,n}}$ is included in
$[0,n]^{d-1} \times [h,h+L n^{d-1}]$. Thus, there is an $h$ such that $E_{\phi_{L,n}}\cap B_{n,k}$ is included in $[k,n-k]^{d-1}\times [h,h+L n^{d-1}]$. Furthermore, since $S_{n,1},\ldots ,S_{n,\psi(n)}$ are pairwise disjoint, there is at least one $k$ in $\{1,\ldots ,\psi(n)\}$ such that 
$$\card (E_{\phi_{L,n}}\cap S_{n,k}) \,\leq\,\frac{ \card (E_{\phi_{L,n}}
  )}{\psi(n)} \,\leq\, \frac{L n^{d-1}}{\psi(n)}\;.$$
Thus, denoting $F=E_{\phi_{L,n}}\cap S_{n,k}$ and $x=\tilde{x}_{E_{\phi_{L,n}}\cap S_{n,k}}$, this shows that $\phi_{L,n}\geq \tau_{(F,x)}$, and claim (\ref{eqphitau}) is proved.

Now, we need to show that $\min_{(F,x)\in I_n}\tau_{(F,x)}/{n^{d-1}}$ has lower large deviations given by $\mathcal{I}_{\vec{v}}$. First, notice that
\begin{eqnarray}
\label{eqsubad}
\PP(\phi_{L,n}\leq \lambda n^{d-1})  & \leq & \PP(\min_{(F,x)\in
  I_n}\tau_{(F,x)}\leq \lambda n^{d-1}) \nonumber  \\
& \leq &\sum_{(F,x)\in I_n}\PP(\tau_{(F,x)}\leq \lambda n^{d-1}) \,.
\end{eqnarray}
Since, according to inequality (\ref{eqIn}), $\log \card (I_n)$ is small
compared to $n^{d-1}$, we shall be done if we can show that, uniformly in
$(F,x)\in I_n$, the probability of deviation $\PP(\tau_{(F,x)}\leq \lambda
n^{d-1})$ is asymptotically of order at most $\exp\left(-\mathcal{I}_{\vec{v}}(\lambda)n^{d-1}\right)$. We shall do this using a subadditivity 
argument. From now on, we fix $(F,x)$ in $I_n$ and $k$ such that $F\subset
S_{n,k}$. The notations and rigorous proofs are a little cumbersome, but
everything can be guessed in two stages, looking at Figures
\ref{fig:patch4} and \ref{fig:etape2}.

\begin{figure}[!ht]
\centering

\begin{picture}(0,0)%
\includegraphics{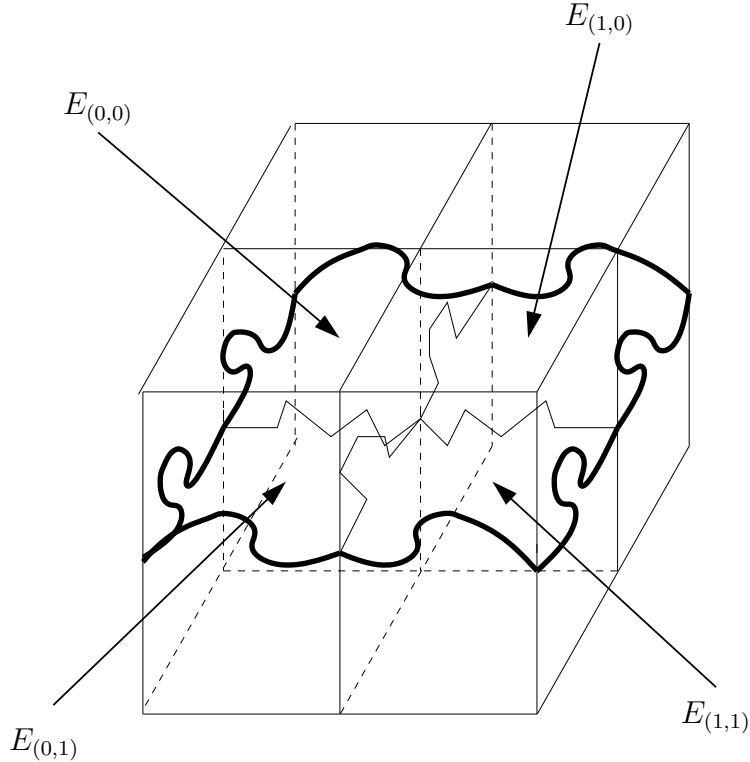}%
\end{picture}%
\setlength{\unitlength}{2960sp}%
\begingroup\makeatletter\ifx\SetFigFont\undefined%
\gdef\SetFigFont#1#2#3#4#5{%
  \reset@font\fontsize{#1}{#2pt}%
  \fontfamily{#3}\fontseries{#4}\fontshape{#5}%
  \selectfont}%
\fi\endgroup%
\begin{picture}(5662,6402)(1111,-7132)
\put(5776,-961){\makebox(0,0)[b]{\smash{{\SetFigFont{12}{14.4}{\rmdefault}{\mddefault}{\updefault}{\color[rgb]{0,0,0}$E_{(1,0)}$}%
}}}}
\put(6751,-6811){\makebox(0,0)[b]{\smash{{\SetFigFont{12}{14.4}{\rmdefault}{\mddefault}{\updefault}{\color[rgb]{0,0,0}$E_{(1,1)}$}%
}}}}
\put(1126,-7036){\makebox(0,0)[b]{\smash{{\SetFigFont{12}{14.4}{\rmdefault}{\mddefault}{\updefault}{\color[rgb]{0,0,0}$E_{(0,1)}$}%
}}}}
\put(1576,-1711){\makebox(0,0)[b]{\smash{{\SetFigFont{12}{14.4}{\rmdefault}{\mddefault}{\updefault}{\color[rgb]{0,0,0}$E_{(0,0)}$}%
}}}}
\end{picture}%

\caption{Patching $E_b$ for $b\in\{0,1\}^{d-1}$ when $d=3$.}
\label{fig:patch4}
\end{figure}

\begin{figure}[!ht]
\centering

\begin{picture}(0,0)%
\includegraphics{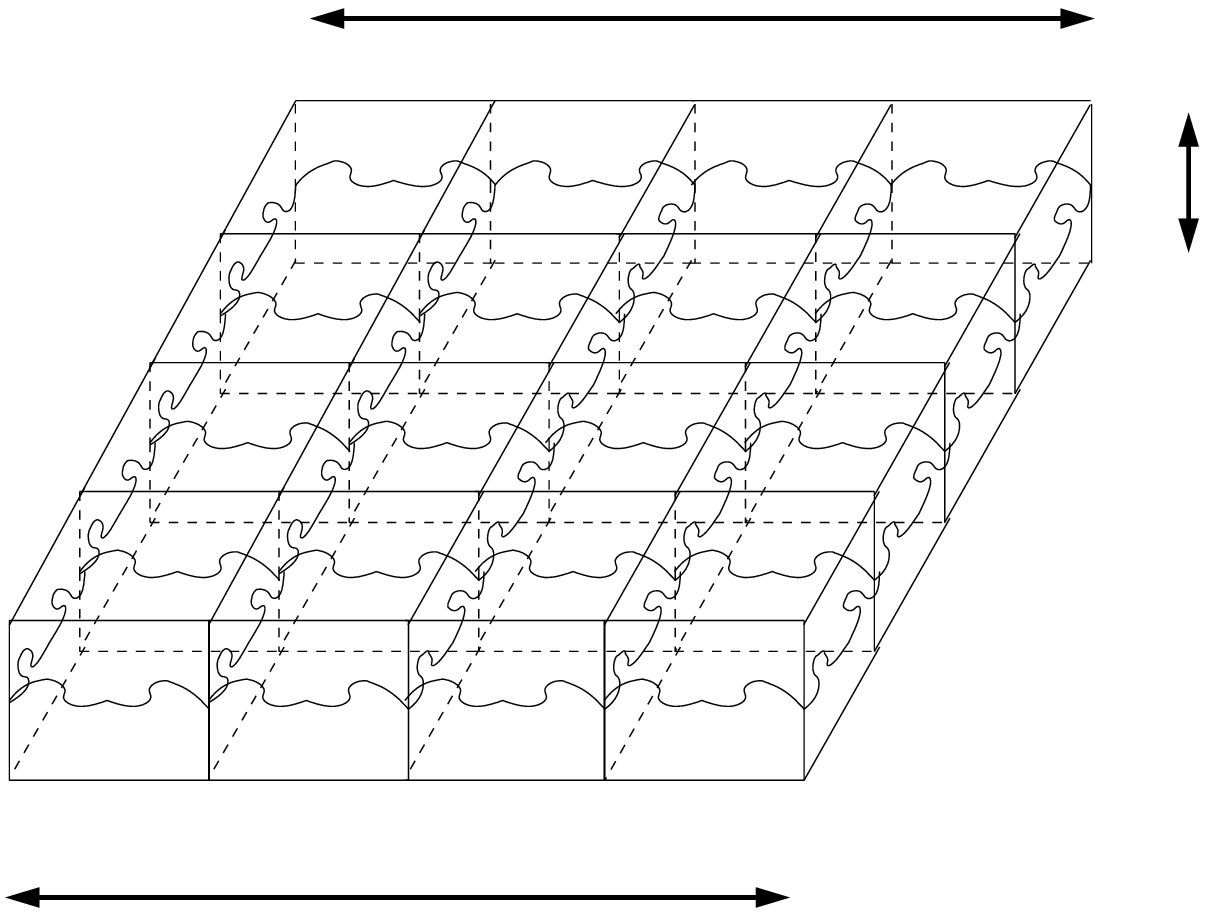}%
\end{picture}%
\setlength{\unitlength}{2960sp}%
\begingroup\makeatletter\ifx\SetFigFont\undefined%
\gdef\SetFigFont#1#2#3#4#5{%
  \reset@font\fontsize{#1}{#2pt}%
  \fontfamily{#3}\fontseries{#4}\fontshape{#5}%
  \selectfont}%
\fi\endgroup%
\begin{picture}(7923,6528)(418,-7357)
\put(8326,-2386){\makebox(0,0)[lb]{\smash{{\SetFigFont{12}{14.4}{\rmdefault}{\mddefault}{\updefault}{\color[rgb]{0,0,0}$2h(n)$}%
}}}}
\put(4951,-1036){\makebox(0,0)[b]{\smash{{\SetFigFont{12}{14.4}{\rmdefault}{\mddefault}{\updefault}{\color[rgb]{0,0,0}$2N$ cubes}%
}}}}
\put(2851,-7261){\makebox(0,0)[b]{\smash{{\SetFigFont{12}{14.4}{\rmdefault}{\mddefault}{\updefault}{\color[rgb]{0,0,0}$4N(n-2k)$}%
}}}}
\end{picture}%

\caption{Patching cuts with the same perimeter.}
\label{fig:etape2}
\end{figure}

Let $N$ be an integer such that for every $N'\geq N$, $h(2(n-2k)N')\geq h(n)$. Define, for $i=1,\ldots ,d-1$, the following hyperplanes:
$$H_i\,=\,\RR^{i-1}\times\{n-k\}\times \RR^{d-i}\;.$$
We define $\sigma_i$ to be the affine orthogonal reflection relative to $H_i$, and $\tr_i$ to be the following translation along coordinate $i$:
$$\tr_i(z)\,=\,z+2(n-2k)\mathbf{e}_{i}\;,$$
where $(\mathbf{e}_1,...,\mathbf{e}_d)$ is the canonical orthonormal basis of
$\mathbb{R}^d$. For any $b\in \{-2N,\ldots,2N-1\}^{d-1}$, we define
the map $\sigma_b$ as follows. For every $i$ in $\{1,\ldots,d-1\}$, let
$a_i=\lfloor b_i/2\rfloor$ and $c_i=b_i-2a_i$. Then, we denote by
$\sigma_b$ the (commutative) product of translations and reflections
$\prod_{i=1}^{d-1}\tr_i^{a_i}\circ\prod_{i=1}^{d-1}\sigma_i^{c_i}$, where
$\sigma_i^{c_i}$ (respectively $\tr_i^{a_i}$) is the $c_i$-th iterate of
$\sigma_i$ (respectively the $a_i$-th iterate of $\tr_i$). Finally, we define also, for any set of vertices or set of edges $X$,
$$\sigma_N(X)\,=\,\bigcup_{b\in \{-2N,\ldots,2N-1\}^{d-1}}\sigma_b(X)\;,$$
and
$$\tilde{\sigma}_N(X)\,=\,\sigma_N(X)\cap S_N\;,$$
where
$$S_N\,=\,\partial \left([k-2N(n-2k),k+2N(n-2k)]^{d-1}\right)\times [-h(n),h(n)]\;.$$
The following lemma should be intuitive looking at Figures \ref{fig:patch4}
and \ref{fig:etape2}. In words, the main message of this lemma
(assertion $(ii)$) is the
following.  Let  $E$ (resp. $E'$) be a 
cut between the top and the bottom in some box $B$
(resp. $B'$). Suppose that $B$ and $B'$ share exactly a face, and that
the connection functions of $E$ and $E'$ coincide on this face. Then,
$E\cup E'$ is a cut between the top and the bottom in $B\cup
B'$. Notice that assertion $(i)$ is just an obvious property of symmetry: if you
take a cut $E$ between the top and the bottom in a box $B$, then
$\sigma_b(E)$ is a cut between the top and the bottom in $\sigma_b(B)$,
for any $b$.
\begin{lem}
\label{lem:connection}
Let $(F,x)$ be fixed in $I_n$. Suppose that for every $b\in
\{-2N,\ldots,2N-1\}^{d-1}$, we are given a set $E_b$ of edges that cuts $B(\sigma_b(A_{n,k}),h(n))$
from $T(\sigma_b(A_{n,k}),h(n))$ in $\sigma_b(B_{n,k})$. Let
$\mathbf{0}$ denote ${(0,\ldots,0)}$ and define:
$$E=\bigcup_{b\in \{-2N,\ldots,2N-1\}^{d-1}}E_b\;.$$
\begin{itemize}
\item[(i)] If $E_{\mathbf{0}}\cap S_{n,k}=F$, and $\tilde{x}_{E_{\mathbf{0}}}=x$, then for every $b\in\{-2N,\ldots,2N-1\}^{d-1}$, the set of edges
  $\sigma_b(E_{(0,\ldots,0)})$ cuts $B(\sigma_b(A_{n,k}),h(n))$ from $T(\sigma_b(A_{n,k}),h(n))$ in
  $\sigma_b(B_{n,k})$, has configuration function $x\circ\sigma_b^{-1}$, and satisfies
$$\sigma_b(E_{(0,\ldots,0)})\cap\sigma_b(S_{n,k})\,=\,\sigma_b(F)\;.$$

\item[(ii)] If, for every $b\in
\{-2N,\ldots,2N-1\}^{d-1}$, 
$$\tilde{x}_{E_b}\circ\sigma_b\,=\,x\;,$$
then $E$ cuts
  $B(\sigma_N(A_{n,k}),h(n))$ from $T(\sigma_N(A_{n,k}),h(n))$ in
  $\sigma_N(B_{n,k})$.


\end{itemize}
\end{lem}
\begin{dem}
Assertion $(ii)$ is the only non-trivial point to show. Let $b$ and $b'$ be two members of
$\{-2N,\ldots,2N-1\}^{d-1}$. The hypotheses on the cuts $E_b$ and $E_b'$
ensure that $x_{E_b}$ and $x_{E_{b'}}$ coincide on
$\sigma_b(B_{n,k})\cap\sigma_{b'} (B_{n,k})$. Thus, we can extend all
the functions $(x_b)_{b\in \{-2N,\ldots,2N-1\}^{d-1}}$ in a single
function $x$ on $\sigma_N(B_{n,k})$. This implies that for every two
neighbours $u$ and $v$ in $\sigma_N(B_{n,k})$, if $\langle
u,v\rangle\not\in E$, then $x(u)=x(v)$. Thus, $x$ is constant on each
connected component of $\sigma_N(B_{n,k})\smallsetminus E$. Since in
each box $\sigma_b(B_{n,k})$, $E_b$ cuts $B(\sigma_b(A_{n,k}),h(n))$
from $T(\sigma_b(A_{n,k}),h(n))$, we have that $x$ takes the value $1$ on $B(\sigma_N(A_{n,k}),h(n))$, and $0$ on  $T(\sigma_N(A_{n,k}),h(n))$. Thus, these two sets are disconnected in $\sigma_N(B_{n,k})\smallsetminus E$.
\end{dem}
Now, for every $b\in \mathbb{Z}^{d-1}$, define
\begin{eqnarray*}
\mathcal{C}_{F,x,b}&=&\left\{E\subset \EE^d \,|\,E\mbox{ is  a
}(B(\sigma_b(A_{n,k}),h(n)),T(\sigma_b(A_{n,k}),h(n)))\mbox{-cut in
}\sigma_b(B_{n,k}) \,,\right.\\ &&\,\,\,\,\left. E\cap
\sigma_b(S_{n,k})=\sigma_b(F) ,\;\card (E)\leq L n^{d-1}\mbox{ and }\tilde{x}_E\circ\sigma_b=x\right\}\;,
\end{eqnarray*}
and
$$\tau_{(F,x,b)}=\min \left\{V(E)\,|\, E\in\mathcal{C}_{F,x,b}\right\}\;.$$
For every $N$, let $E_N$ denote the set of the edges $e$ in $\sigma_N(B_{n,k})$  such that at least one endpoint of $e$ belongs to $S_N$. Define
$M(N)=N+\psi(N)$ and, for $N'\in\{N,M(N)\}$,
$$\tau_{N'}\,=\,\tau(\sigma_{N'}(A_{n,k}),h(N'))\;.$$
If $E$ is a $(B(\sigma_N(A_{n,k}),h(n)),T(\sigma_N(A_{n,k}),h(n)))$-cut in 
$\sigma_N(A_{n,k})$, $E\cup E_{N'}$ clearly cuts $\sigma_N(A_{n,k})_1^{h(n)}$ from $\sigma_N(A_{n,k})_2^{h(n)}$. Thus, part $(ii)$ of Lemma \ref{lem:connection} gives us that
$$\sum_{b\in\{-2M(N),\ldots ,2M(N)-1\}^{d-1}}\tau_{(F,x,b)}+\min_{N'\in\{N,\ldots,M(N)\}}\sum_{e\in
  E_{N'}}t(e) \,\geq\, \min_{N'\in\{N,\ldots,M(N)\}}\tau_{N'}\;.$$
Notice that some edges are counted twice on the left-hand side of the
preceding inequality. From part $(i)$ of Lemma
  \ref{lem:connection}, we know that the random variables
  $(\tau_{(F,x,b)})_{b\in\{-2M(N),\ldots ,2M(N)-1\}^{d-1}}$ are identically distributed,
  with the same distribution as $\tau_{(F,x)}$. Using the FKG inequality,
\begin{align*}
\PP(\tau_{(F,x)}\leq \lambda &  n^{d-1})^{(4M(N))^{d-1}} \\
&\,=\, \prod_{b\in\{-2M(N),\ldots
  ,2M(N)-1\}^{d-1}}\PP(\tau_{(F,x,b)}\leq \lambda n^{d-1})\\
&\,\leq \,\PP(\forall b\in\{-2M(N),\ldots ,2M(N)-1\}^{d-1},\;\tau_{(F,x,b)}\leq \lambda
  n^{d-1})\\
&\,\leq \,\PP\left(\sum_{b\in\{-2M(N),\ldots ,2M(N)-1\}^{d-1}}\tau_{(F,x,b)}\leq \lambda
  n^{d-1}(4M(N))^{d-1}\right)\\
&\,\leq \,\PP\left(\min_{N'\in\{N,\ldots,M(N)\}}\tau_{N'}-\min_{N'\in\{N,\ldots,M(N)\}}\sum_{e\in
  E_{N'}}t(e)\leq \lambda
  n^{d-1}(4M(N))^{d-1}\right)\,.\\
\end{align*}
Let $\eps>0$ be a fixed positive real number.
\begin{eqnarray}
\nonumber \PP(\tau_{(F,x)}\leq \lambda n^{d-1})&\leq &\Bigg( \PP\bigg(\min_{N'\in\{N,\ldots,M(N)\}}\tau_{N'}\leq (\lambda+\eps)
  n^{d-1}(4M(N))^{d-1}\bigg)\\ 
\label{eq:majotauF1}&& +\PP\bigg(\min_{N'\in\{N,\ldots,M(N)\}}\sum_{e\in
  E_{N'}}t(e)\geq \eps n^{d-1}(4M(N))^{d-1}\bigg)\Bigg)^{\frac{1}{(4M(N))^{d-1}}}\;.
\end{eqnarray}
Now we shall let $N$ go to infinity. Using Lemma
\ref{limitetaupert},  the fact that $\lim_{N\rightarrow \infty}\psi(N)/N=
0$ and a union bound,

\begin{align}
\nonumber &\liminf_{N\rightarrow \infty} -\frac{1}{(4N(n-2k))^{d-1}}\log \PP\bigg(\min_{N'\in\{N,\ldots,M(N)\}}\tau_{N'}\leq (\lambda+\eps)
  n^{d-1}(4M(N))^{d-1}\bigg)\\
\nonumber &\,\geq\,\liminf_{N\rightarrow
  \infty}-\frac{1}{(4N(n-2k))^{d-1}} \times \\
\nonumber &\quad\max_{N'\in\{N,\ldots,M(N)\}}\log \PP\left(\tau_{N'}\leq \left(\lambda+2\eps-\frac{1}{\sqrt{4(n-2k)N'}}\right)
  n^{d-1}(4N')^{d-1}\right)\\
\label{eq:majotauF2}&\,\geq\, \mathcal{I}_{\vec{v}}\left((\lambda+2\eps )\left(\frac{n}{n-2k}\right)^{d-1}\right)\;.
\end{align}
Now, we use the fact that $F$ possesses an exponential moment, and that the
sets $E_{N'}$ are disjoint. Using
Chebyshev inequality, there are positive constants $C_4$ and $C_5$, depending only on $F$ and $d$, such that
\begin{align*}
\PP\bigg( \min_{N'\in\{N,\ldots,M(N)\}} &  \sum_{e\in
  E_{N'}}t(e)\geq \eps n^{d-1}(4M(N))^{d-1}\bigg)\\
&\,=\,\prod_{N'\in\{N,\ldots,M(N)\}}\PP\bigg(\sum_{e\in
  E_{N'}}t(e)\geq \eps n^{d-1}(4M(N))^{d-1}\bigg)\;,\\
&\,\leq\, \left(e^{C_4h(n)n^{d-2}M(N)^{d-2}-C_5n^{d-1}(4M(N))^{d-1}}\right)^{\psi(N)}\;.\\
\end{align*}
Thus,
$$\liminf_{N\rightarrow \infty}-\frac{1}{(4N(n-2k))^{d-1}}\log \PP\left(\min_{N'\in\{N,\ldots,M(N)\}}\sum_{e\in
  E_{N'}}t(e)\geq \eps n^{d-1}(4M(N))^{d-1}\right)\,=\,+\infty\;.$$
Therefore, inequalities (\ref{eq:majotauF1}) and (\ref{eq:majotauF2}) imply:
$$-\frac{1}{n^{d-1}}\log \PP(\tau_{(F,x)}\leq \lambda n^{d-1})\,\geq\,\frac{(n-2k)^{d-1}}{n^{d-1}}\mathcal{I}_{\vec{v}}\left((\lambda+2\eps )\left(\frac{n}{n-2k}\right)^{d-1}\right)\;.$$
We choose $\eps=\frac{1}{n}$, and replace $\lambda$ by
$\lambda-\frac{1}{\sqrt{n}}$ to get
$$-\frac{1}{n^{d-1}}\log \PP\left(\tau_{(F,x)}\leq \left(\lambda-\frac{1}{\sqrt{n}}\right)n^{d-1}\right)\,\geq\,\frac{(n-2k)^{d-1}}{n^{d-1}}\mathcal{I}_{\vec{v}}\left(\left(\lambda-\frac{1}{\sqrt{n}}+\frac{2}{n}\right)\left(\frac{n}{n-2k}\right)^{d-1}\right)\;.$$
Since $k\leq \psi(n)$ and $\psi(n)$ is small compared to $\sqrt{n}$,
and since $\mathcal{I}_{\vec{v}}$ is non-increasing, for $n$
large enough,
$$-\frac{1}{n^{d-1}}\log \PP\left(\tau_{(F,x)}\leq
\left(\lambda-\frac{1}{\sqrt{n}}\right)
n^{d-1}\right)\geq\frac{(n-2k)^{d-1}}{n^{d-1}}\mathcal{I}_{\vec{v}}\left(\lambda\right)\;.$$
Using inequalities (\ref{eqIn}) and  (\ref{eqsubad}),
$$\liminf_{n\rightarrow \infty}-\frac{1}{n^{d-1}}\log \PP(\phi_{L,n}\leq \lambda n^{d-1})\geq \mathcal{I}_{\vec{v}}(\lambda)\;.$$
And thus, from inequality (\ref{eq:moinsdaretes}),
$$\liminf_{n\rightarrow \infty}-\frac{1}{n^{d-1}}\log \PP(\phi_{n}\leq \lambda
n^{d-1})\geq
\min\{\mathcal{I}_{\vec{v}}(\lambda),C'_2L\}\;.$$
Letting $L$ tend to infinity finishes the proof of Proposition \ref{prop:fntauxphi}.
\end{dem}
\begin{rem}
\label{rem:phiechoue}
This ``symmetric-subadditive'' argument does not work in the ``non-straight'' 
case. It is
perhaps important to note that in this case, it is not obvious
at all to know in advance for which $F$ the random variable $\tau_F$ has the ``minimal'' distribution. It is
natural to conjecture that this ``minimal'' $F$ is a hyperrectangle, but we do
not know how to prove this for all dimensions. When $d=2$, though, we are
able to solve this problem and to show that if $h(n) /n$ converges towards
$\tan(\alpha)$ for some $\alpha$ in $[0, \pi /2]$, and if $\vec{v} = (\cos
\theta, \sin \theta) = \vec{v_\theta}$ is orthogonal to $A = A_\theta$,
then $\phi(nA_\theta, h(n))/n$ converges towards $\min
\{\nu(\vec{v_{\widetilde{\theta}}}) / \cos (\widetilde{\theta}
-\theta) \,\mbox{ s.t. }\,
|\widetilde{\theta} -\theta | \leq \alpha \}$. A similar method gives
an analog result for the lower large deviations. This will be done
rigorously in a forthcoming paper.
\end{rem}


\subsection[LLN for $\phi$]{Law of large numbers}
\label{sec:llnphistraight}
In this section, we prove Theorem \ref{thm:llnphistraight}. So we
suppose that $\mathbf{(F2)}$, $\mathbf{(H1)}$ and
$\mathbf{(H2)}$ hold, and that $A$ is a straight (so non-degenerate)
hyperrectangle. Notice first that if $\mathbf{(F1)}$ does not hold,
then $\nu$ always equals zero (cf. Proposition \ref{thmnu}) and the
law of large numbers for $\phi$ is a consequence of the one for
$\tau$, Theorem \ref{thm:llntau}. Thus, we may suppose that
$\mathbf{(F1)}$ holds. We first prove the a.s. convergence of the
rescaled variable. Since
$\phi (nA,h(n)) \leq \tau (nA,h(n))$,  we only need to show that
\begin{equation}
\label{eq:llnphiinf}
\liminf _{n\rightarrow\infty}\frac{\phi
  (nA,h(n))}{\H^{d-1}(nA)}\geq \nu(\vec{v_0})\quad a.s.
\end{equation}
where $\vec{v_0}=(0,\ldots,0,1)$.  

\emph{Suppose first that $F$ has bounded support.} Then, $\mathbf{(F4)}$ is
obviously satisfied, and we deduce from Proposition
\ref{prop:fntauxphi}, the positivity of $\mathcal{I}_{\vec{v}}$ on
$[0,\nu(\vec{v})[$ and Borel-Cantelli's lemma that (\ref{eq:llnphiinf})
is true. 

\emph{Now, let $F$ be general}, i.e. satisfy $\mathbf{(F2)}$. We rely on the
ideas of Proposition \ref{prop:deviation}. Let $a>1/2$ be
a real number to be chosen later, define $\tilde
t(e)=t(e)\land a$ and let $F_a$ be the distribution function of $\tilde
t(e)$. We define:
$$ \tilde \tau_n\,=\, \min\left\{ \begin{array}{c} \sum_{e\in
      E}\tilde t(e)\mbox{ s.t. }E \mbox{ cuts }\\
 (nA)_1^{h(n)}\mbox{
  from }(nA)_2^{h(n)}\mbox{ in }\cyl(nA,h(n)) \end{array}
\right\}\;, $$
and define analogously $\tilde\phi_n$. We use the notations $\nu_{F}(\vec{v})$ (resp. $\nu_{F_a}(\vec{v})$)
  to denote the limit of the rescaled flow $\tau$ corresponding to
  capacities of distribution function $F$ (resp. $F_a$).
As we obtained (\ref{eq:nucontinue}), we get:
$$ \EE(\tau_n)-\EE(\tilde \tau_n)\leq \EE\left(t(e_1)\II_{t(e_1)\geq
    a}\right)\EE(\card E_{\tilde\tau_n})\;.$$
Proposition \ref{prop:5.8Zhangbis} implies that there are constants
$\eps$, $C_1$ and $C_2$ such that:
$$\EE(\card E_{\tilde \tau_n})\leq
\frac{C_1}{C_2}+\frac{1}{\eps}\EE(\tilde\tau_n)\leq \frac{C_1}{C_2}+\frac{1}{\eps}\EE(\tau_n)\;,$$
where the constants $\eps$, $C_1$ and $C_2$ depend only on $d$ and $F$, and
not on $a$, since $F$ and $ F_a$ coincide on $[0,1/2]$. Then, for any
$\eps'>0$, one can choose $a$ large enough so that:
$$ \EE(\tau_n)-\EE(\tilde \tau_n)\leq \eps \H^{d-1}(nA)\;,$$
leading to $\nu_F(\vec{v_0})-\nu_{F_a}(\vec{v_0})\leq \eps$. Since
$\phi_n\geq \tilde\phi_n$, and using the result for $F_a$ which has
bounded support, we get for every $\eps' >0$:
$$\liminf _{n\rightarrow\infty}\frac{\phi
  (nA,h(n))}{\H^{d-1}(nA)}\geq \nu_{F_a}(\vec{v_0})\geq
\nu_F(\vec{v_0})-\eps\quad a.s.$$
Which gives the desired result.

It remains to prove the convergence in $L^1$. It may be derived
exactly as in the proof of the convergence in $L^1$ of Theorem
\ref{thm:llntau} as soon as we have proved the convergence of the
expectation of the rescaled maximal flow. But this is immediate thanks to Fatou's lemma:
\begin{align*}
 \nu(\vec{v_0}) \,=\, & \EE \left[ \lim_{n\rightarrow \infty}
  \frac{\phi(nA, h(n))}{\H^{d-1}(nA)}  \right]  \,\leq\, \liminf_{n\rightarrow
  \infty} \EE \left[ \frac{\phi(nA, h(n))}{\H^{d-1}(nA)}  \right]
\,\leq\,\\
& \,\leq\,  \limsup_{n\rightarrow
  \infty} \EE \left[ \frac{\phi(nA, h(n))}{\H^{d-1}(nA)}  \right]  \,\leq\, \lim_{n\rightarrow
  \infty} \EE \left[ \frac{\tau(nA, h(n))}{\H^{d-1}(nA)}  \right] \,=\,
\nu(\vec{v_0}) \,. 
\end{align*}
This ends the proof of Theorem \ref{thm:llnphistraight}.
\begin{rem}
\label{rem:improveZhang}
Our proof of Proposition \ref{prop:fntauxphi} can be carried out in
Kesten and Zhang's setting \cite{Zhang07}, who
consider $A_{\mathbf{k}}=\prod_{i=1}^{d-1}[0,k_i]\times\{0\}$ with
$k_1\leq\ldots\leq k_{d-1}$ and
let all the $k_i$ go to infinity, possibly at different speeds. The
only obstacle to do this is when one reduces the sides of the box:
$\psi(n)$ has to be replaced by $\psi(k_1)$, and the set $I_n$ by a
set $I_{\mathbf{k}}$ satisfying:
$$\card(I_{\mathbf{k}})\,\leq \,
C_3h(\mathbf{k})\psi(k_1)\left(C_4L\prod_{i=1}^{d-1}k_i\right)^{LC_5
  \prod_{i=1}^{d-1}k_i/\psi(k_1)}\;,$$
where $C_3$, $C_4$ and $C_5$ are constants depending on $d$ and $h(\mathbf{k})$ is
the height of the box. Then, the proof works as long as
$\log\card(I_{\mathbf{k}})$ is small with respect to
$\prod_{i=1}^{d-1}k_i$, which is the case if $\log h(\mathbf{k})$ is
small with respect to $\prod_{i=1}^{d-1}k_i$ and $\log k_{d-1}$ is small
with respect to $k_1$. Thus,  we obtain the law of large numbers (and
also Proposition \ref{prop:fntauxphi}) under the conditions:
$$\left\{\begin{array}{l}(\log
      h(\mathbf{k}))/\prod_{i=1}^{d-1}k_i\xrightarrow[\mathbf{k}\rightarrow\infty]{}0\\
    (\log
  k_{d-1})/k_1\xrightarrow[\mathbf{k}\rightarrow\infty]{}0 \end{array}\right.\;$$
and condition $\mathbf{(F1)}$. So, in a sense, the height condition is better than in Theorem
\ref{cvphi} (and essentially optimal), however we are not able to get
rid of the second ugly condition -  which imposes that the sides of $A$ do
not have too different asymptotic behaviours - without requiring a stronger condition
on $h$, similar to the one of Kesten and Zhang.
\end{rem}

\subsection{Final steps of the proofs of Theorem \ref{devinfphi} and
  Theorem \ref{thmpgd}}
\label{secfinal}

The proof of Theorem \ref{devinfphi} is exactly the same as the one of
Theorem \ref{devinf}, using Theorem \ref{thm:llnphistraight} and
Proposition \ref{prop:deviation}. It remains
to end the proof of Theorem \ref{thmpgdphi}. Proposition
\ref{prop:fntauxphi} states in a sense that $\phi$ and $\tau$ share the
same rate function. Since this function has already been studied, and since
the upper large deviations of $\phi$ have been studied in \cite{TheretUpper}, the
construction of the rate function of $\phi$ was the main work to do in
order to show the large deviation principle for $\phi$ in straight
boxes. Indeed, the only thing we have to prove is that for all $\lambda > \nu(\vec{v})$,
\begin{equation}
\label{devsupphi}
\lim_{n\rightarrow \infty} \frac{1}{\mathcal{H}^{d-1}(nA)} \log
\mathbb{P} \left[ \frac{\phi(nA,h(n))}{\mathcal{H}^{d-1}(nA)} \geq \lambda \right]\,=\, -\infty  \,.
\end{equation}
As soon as we have (\ref{devsupphi}), we can write exactly the same proof
for Theorem \ref{thmpgdphi} as for Theorem \ref{thmpgd} (see section
\ref{pgd}), since we
have proved that $\phi(nA, h(n))/\H^{d-1}(nA)$ converges a.s. to
$\nu(\vec{v_0})$. To obtain (\ref{devsupphi}), we can refer to section 3.7 in
\cite{TheretUpper} (here only the existence of one exponential moment is
required).

\begin{rem}
\label{rem:condF4phi}
We leave the following questions open: is condition $\mathbf{(F4)}$
necessary to obtain the existence of a rate function for $\phi$ ? If this
rate function exists under weaker hypothesis than $\mathbf{(F4)}$, is it
necessarily the same as the one for $\tau$ ? When the rate function exists, do
we necessarily obtain the corresponding large deviation principle ?
\end{rem}


\section*{Acknowledgements}
We are very indebted to Rapha\"el Cerf for the time he devoted to us on
helpful discussions and his numerous valuable comments and advises.
We also would like to express our gratitude to Marc Wouts for stimulating
conversations. We thank an anonymous referee for his many
  relevant questions that incited us notably to improve our moment conditions in
  the entire article.



\def\cprime{$'$}

\end{document}